\documentclass[10pt]{amsart}
\usepackage[leqno]{amsmath}
\usepackage{amsthm}
\usepackage{amsfonts}
\usepackage{amssymb}
\usepackage{eucal}
\usepackage[all]{xy}

\DeclareFontFamily{OT1}{rsfs}{}
\DeclareFontShape{OT1}{rsfs}{n}{it}{<-> rsfs10}{}
\DeclareMathAlphabet{\mathscr}{OT1}{rsfs}{n}{it}

\DeclareMathOperator{\Hom}{Hom}

\DeclareMathOperator{\Spec}{Spec}

\DeclareMathOperator{\et}{\acute{e}t}
\newcommand*{\invlim}{\varprojlim}                               % invlimit
\newcommand*{\Z}{\ensuremath{\mathbf{Z}}}                        % integers
\renewcommand*{\P}{\ensuremath{\mathbf{P}}}                        % proj space
\renewcommand*{\setminus}{-}
\newcommand*{\ms}{\mathscr}                                      % font love
\newcommand*{\calA}{\mathscr{A}}
\theoremstyle{plain}
  \newtheorem{theorem}{Theorem}[subsection]
  \newtheorem{proposition}[theorem]{Proposition}
  \newtheorem{lemma}[theorem]{Lemma}
  \newtheorem{corollary}[theorem]{Corollary}
  
\theoremstyle{remark} 
  \newtheorem{remark}[theorem]{Remark}
  \newtheorem{example}[theorem]{Example}
  \newtheorem{notation}[theorem]{Notation}
  
\theoremstyle{definition}
  \newtheorem{definition}[theorem]{Definition} 

\numberwithin{equation}{theorem}
\renewcommand{\theequation}{\ifnum\value{theorem}>0{\thetheorem.\arabic{equation}}\else{\thesubsection.\arabic{equation}}\fi}

\title{Nagata compactification for algebraic spaces}
\author{Brian Conrad}
\address{Department of Mathematics\\
Stanford University\\
Stanford, CA 48109, USA}
\email{conrad@math.stanford.edu}
\author{Max Lieblich}
\address{Department of Mathematics\\
University of Washington\\
Seattle, WA  98195 USA}
\email{lieblich@math.washington.edu}
\author{Martin Olsson}
\address{Department of Mathematics\\
Univ. of California\\
Berkeley, CA 94720, USA}
\email{molsson@math.berkeley.edu}
\date{October 25, 2010}
\thanks{BC was partially supported by NSF grant DMS-0917686, ML was
  partially supported by NSF grant DMS-0758391 and the Sloan
  Foundation, and MO was partially supported by NSF grants DMS-0714086 and DMS-0748718
and the Sloan Foundation.  The authors are grateful to the
Mathematical Institute at Oberwolfach and MSRI for their hospitality
and atmosphere, and to Johan de Jong and especially Ofer Gabber for
helpful suggestions at the start of this project.  We also received numerous helpful comments from David
Rydh on an early draft of this paper. We are also grateful to  Rydh for obtaining and forwarding a copy of the unpublished thesis of Raoult \cite{raoulthesis}, after the first version of this paper was written.}
\begin{document}

\begin{abstract}
We prove the Nagata compactification theorem for any separated map of finite type between quasi-compact and quasi-separated algebraic spaces, generalizing earlier results of Raoult.  Along the way we also prove (and use) 
absolute noetherian 
approximation for such algebraic spaces, generalizing earlier results in the 
case of schemes. 

\end{abstract}

\maketitle

\centerline{{\em To the memory of Masayoshi Nagata}}

\section{Introduction}

\subsection{Motivation}

The Nagata compactification theorem for schemes is a very useful and fundamental result.  It says that if $S$ is a quasi-compact and quasi-separated scheme (e.g., any noetherian scheme)
and if $f:X \rightarrow S$ is a separated map of finite type from a scheme $X$ then $f$ fits into a commutative diagram of schemes
\begin{equation}\label{factorf}
\xymatrix{
{X} \ar[r]^-{j} \ar[dr]_-{f} & {\overline{X}} \ar[d]^-{\overline{f}} \\
& {S}}
\end{equation}
with $j$ an open immersion and $\overline{f}$ proper;
we call such an $\overline{X}$ an {\em $S$-compactification of $X$}.

Nagata's papers (\cite{n1}, \cite{n2}) focused on the case of
noetherian schemes and unfortunately are difficult to read nowadays (due to the use of
an older style of algebraic geometry), but there are several available proofs 
in modern language.  The proof by L\"utkebohmert \cite{l} applies in the noetherian case, and 
the proof of Deligne (\cite{d}, \cite{conrad}) is a modern interpretation of Nagata's method
which applies in the general scheme case.  The preprint \cite{vojta} by Vojta 
gives an exposition of Deligne's approach in the noetherian case.  Temkin has 
recently introduced some new valuation-theoretic ideas that give yet another
proof in the general scheme case. 
The noetherian case is the essential one for proving the theorem
because it implies the general case via approximation
arguments \cite[Thm.~4.3]{conrad}.

An important application of the Nagata compactification theorem for
schemes is in the definition of \'etale cohomology with proper
supports for any separated map of finite type $f:X \rightarrow S$
between arbitrary schemes. Since any algebraic space is
\'etale-locally a scheme, the main obstacle to having a similar
construction of such a theory for \'etale cohomology of algebraic
spaces is the availability of a version of Nagata's theorem for
algebraic spaces.  Strictly speaking, it is possible to develop the
full ``six operations'' formalism even for non-separated Artin stacks
(\cite{lo1}, \cite{lo2}) despite the lack of a compactification
theorem in such cases.  However, the availability of a form of
Nagata's theorem simplifies matters tremendously, and there are
cohomological applications for which the approach through
compactifications seems essential, such as the proof of Fujiwara's
theorem for algebraic spaces \cite{varshavsky} (from which one can
deduce the result for Deligne--Mumford stacks via the use of coarse
spaces).  The existence of compactifications is useful in many
non-cohomological contexts as well.

\subsection{Main results} 
In his (unpublished) thesis \cite{raoulthesis}, Raoult proved the
Nagata compactification theorem in the following situations:
\begin{enumerate}
\item $X$ is a normal algebraic space and $S$ is noetherian
  \cite[Proposition, page 2]{raoulthesis};
\item $X$ is arbitrary and $S$ is the spectrum of $\Z$, a field $k$,
  or a complete local noetherian ring \cite[Theorem, page
  11]{raoulthesis}.
\end{enumerate}
A sketch of the second case when $S$ is the spectrum of a field is
contained in \cite{raoult2}.  One can deduce from the second case that
compactifications exist when $S$ is separated and of finite type over
$\Z$, a field, or a complete local noetherian ring.  Raoult's argument
was conditional on an unpublished result of Deligne concerning the
existence and properties of quotients by finite group actions on
separated algebraic spaces.  This result of Deligne is a
consequence of subsequent general results on the existence of coarse
moduli spaces associated to Artin stacks, which we shall review later (and it
is also addressed in \cite[Thm.~5.4]{rydh-qt} in a stronger form).

Following a sketch suggested to us by Gabber, we handle an essentially arbitrary
base $S$:

\begin{theorem}\label{mainresult} 
  Let $f:X \rightarrow S$ be a separated map of finite type between
  algebraic spaces, with $S$ quasi-compact and quasi-separated.  There
  exists an open immersion $j:X \hookrightarrow \overline{X}$ over $S$
  into an algebraic space $\overline{X}$ that is proper over $S$.  If
  $f$ is finitely presented then $\overline{X}$ may be taken to be
  finitely presented over $S$.
\end{theorem}

The proof of Theorem \ref{mainresult} consists of two rather separate
parts: technical approximation arguments reduce the proof to the case
when $S$ is of finite presentation over $\Spec \Z$ (or an excellent
noetherian scheme), and geometric arguments handle this special case
by reducing to the known normal case.  (We include a proof in 
the normal case, so there is no logical dependence on the unpublished
\cite{raoulthesis}.)  We became aware of Raoult's
work only after working out our proof in general, and our basic
strategy (after reducing to $S$ of finite presentation over $\Z$) is
similar to Raoult's.  The reader who is only interested in the case
that $S$ is of finite presentation over $\Z$ or a field (or more
generally, an excellent noetherian scheme) may go immediately to
\S\ref{exccase} and can ignore the rest of this paper.  Theorem
\ref{mainresult} has also been announced by K. Fujiwara and F. Kato,
to appear in a book in progress, as well as by D. Rydh (who has also
announced progress in the case of Deligne--Mumford stacks).

The approximation part of our proof contains some results which are useful for eliminating noetherian hypotheses more generally, so we now make some remarks on this feature of
the present work.   Limit methods of Grothendieck 
\cite[IV$_3$,~\S8--\S11]{ega} are a standard technique for
avoiding noetherian hypotheses, and a very useful supplement to these methods
is the remarkable \cite[App.~C]{tt}.  The key innovation there going beyond
\cite{ega} is an absolute noetherian approximation property \cite[Thm.~C.9]{tt}:  any 
quasi-compact and quasi-separated scheme $S$ admits the form
$S \simeq \invlim S_{\lambda}$ where $\{S_{\lambda}\}$ is an inverse system of $\Z$-schemes of
finite type, with affine transition maps $S_{\lambda'} \rightarrow S_\lambda$ for $\lambda' \ge \lambda$. 
(Conversely, any such limit scheme is obviously affine over any
$S_{\lambda_0}$ and so is necessarily quasi-compact and quasi-separated.) The crux of the
matter is that every quasi-compact and quasi-separated scheme 
$S$ is affine over a $\Z$-scheme of finite type. 
This appromxation result 
was used to avoid noetherian hypotheses in Nagata's theorem for schemes in
\cite{conrad}, and we likewise need a version of it
for algebraic spaces.  This is of interest in its own right, so we state it here (and prove it in
\S\ref{approxsec}).  We note that D. Rydh has a different proof of
this result that also works for certain algebraic stacks \cite{rydh-approx}.

\begin{theorem}\label{absapprox} Let $S$ be a quasi-compact and quasi-separated algebraic
space.  There exists an inverse system $\{S_\lambda\}$ of algebraic spaces of finite presentation 
over $\Z$ such that the transition maps $S_{\lambda'} \rightarrow S_{\lambda}$
are affine for $\lambda' \ge \lambda$ and $S \simeq \invlim S_{\lambda}$.
Moreover, $S$ is separated if and only if $S_{\lambda}$ is separated for sufficiently large $\lambda$.  
\end{theorem}

\begin{remark} The limit algebraic space $\invlim S_{\lambda}$ is defined
\'etale-locally over any single $S_{\lambda_0}$ by using the analogous well-known
limit construction in the case of schemes.  By working \'etale-locally it is easy to check
that such an algebraic space has the universal property of an inverse limit in the category of
algebraic spaces.  Also, since $\invlim S_{\lambda}$ is affine over any
$S_{\lambda_0}$, it is quasi-compact and quasi-separated.  
\end{remark}

We now briefly outline the paper.  
We first consider Theorem \ref{mainresult}
when $S$ is of finite presentation over an excellent noetherian scheme
(such as $\Spec \Z$).  This case is the focus of our efforts
in  \S\ref{exccase}, and the base $S$ is fixed throughout most this section
but we progressively simplify $X$. In \S\ref{redredcase} we use the cotangent complex 
to reduce to the case when $X$ is reduced.
%  (This is 
% an improvement of \cite[Prop.~3]{raoult2}, which treats the case $S = \Spec k$ for a field $k$;
% arguing via the cotangent complex is also better-suited to generalization
% to Artin stacks.)  
Then 
in \S\ref{normalred} we use a contraction result of Artin \cite[Thm.~6.1]{artin}
and various results of Raynaud--Gruson \cite[I,~\S5.7]{rg} 
to reduce to the case when $X$ is normal.  
% (The proof of this part
% is broadly similar to the proof of \cite[Prop.~2]{raoult2}.) 
The case of normal $X$ is handled in \S\ref{normalproof} 
by using a group quotient argument to reduce to the known 
case (whose proof we also provide, for the
convenience of the reader) when $X$ is normal and $S$ is an excellent noetherian scheme.
Note that this settles Theorem \ref{mainresult} for all ``practical'' cases, but not yet the general noetherian case. 

The passage to the general case of Theorem \ref{mainresult} is the aim of \S\ref{approxsec},
via various approximation methods.  
In \S \ref{absnoetherian}
we use stratification techniques of Raynaud--Gruson for algebraic
spaces to prove Theorem \ref{absapprox} by reducing it to the known case of schemes, and
then in \S\ref{fpresec} we reduce to the proof of Theorem \ref{mainresult} in general to the case when 
$f:X \rightarrow S$ is finitely presented (not just finite type). 
An application of Theorem \ref{absapprox} then allows us
to reduce Theorem \ref{mainresult} to the proved case when $S$ is of finite presentation over $\Z$.  
It is only at this point that the general noetherian case is also settled. 

The appendices 
provide some foundational facts we need concerning algebraic spaces; in \S\ref{notation}
we offer some ``justification'' for the appendices.  Much of what is in the appendices
may be known to some experts, but we did not know of a reference
in the literature for the results discussed there.  The reader who is content with taking $S$ to be finitely presented
over an excellent scheme  in Theorem \ref{mainresult} can ignore
\S\ref{approxsec} and almost all of the appendices, and other readers should probably only consult the
appendices when they are cited in the main text.  

New ideas are needed  to prove a general version of Nagata's theorem 
for Deligne--Mumford stacks. 

\subsection{Terminology and conventions}\label{notation}
 We write {\em qcqs} as shorthand for ``quasi-compact
and quasi-separated'' (for schemes, algebraic spaces, or morphisms between them),
and we freely identify any scheme with the corresponding sheaf of sets that it represents
on the \'etale site of the category of schemes.  

The reader who wishes to understand the proof of Theorem \ref{mainresult} for general
noetherian $S$ (or anything beyond the case of $S$ of finite presentation over an excellent
noetherian scheme) will need to read \S\ref{approxsec}, for which the following comments should be helpful.  
Although quasi-separatedness is required in the definition of an algebraic space in
\cite{knutson}, there are natural reasons for wanting to avoid such foundational
restrictions.  We need to use several kinds of 
pushout and gluing constructions with algebraic spaces, 
and the construction and study of these pushouts becomes unpleasant 
if we cannot postpone the consideration of quasi-separatedness properties until
{\em after} the construction has been carried out.  It is a remarkable fact that quasi-separatedness is not necessary
in the foundations of the theory of algebraic spaces; this was known to some experts long ago, but seems to not be 
as widely known as it should be.  

We define an algebraic space $X$ to be an algebraic
space over $\Spec \Z$ as in \cite[I,~5.7.1]{rg}:  it is an \'etale sheaf on the 
category of schemes such that it is isomorphic to a quotient sheaf
$U/R$ for an \'etale equivalence relation in schemes $R \rightrightarrows U$; there are no
quasi-compactness hypotheses in this definition.
The key point is that by using the method of proof of \cite[I,~5.7.2]{rg}, 
it can be proved that for any such $X = U/R$, the fiber product 
$V \times_X W$ is a scheme for any pair of maps $V \rightarrow X$ and
$W \rightarrow X$ with schemes $V$ and $W$. Such representability
was proved in \cite{knutson} under 
quasi-separatedness hypotheses, and is one of the main reasons that
quasi-separatedness pervades that work.  For the convenience of the reader,
we include a proof of this general representability result in \S\ref{repgeneral}, where
we also show (without quasi-separatedness hypotheses)
that quotients by \'etale equivalence relations in algebraic spaces
are always algebraic spaces.  

The avoidance of
quasi-separatedness simplifies the discussion of a number of gluing constructions. 
% In Example \ref{galex} and Example \ref{nongenex} 
% we illustrate some of the subtleties of non-quasi-separated algebraic spaces. 
% We require noetherian algebraic spaces to be quasi-separated by definition; see
% Definition \ref{noethdef}ff. 
Beware that if one removes quasi-separatedness from the definition of an algebraic space
 then some strange things can happen,
such as non-quasi-separated algebraic spaces admitting an \'etale cover by the affine line over a field
%(Example \ref{galex}) 
and unusual behavior for generic points. 
For this reason, when working with algebraic spaces over a noetherian
scheme it is stronger to say ``finite presentation''
(i.e., finite type and quasi-separated) 
than ``finite type'' (even though for schemes
there is no distinction over a noetherian base).  

Whenever we use a result from \cite{knutson} we
are either already working with quasi-separated algebraic spaces or it is trivial
to reduce the desired assertion to the case of quasi-separated algebraic spaces
(such as by working \'etale-locally on the base).  
Note also that the concept of ``algebraic space over a scheme $S$'' in the sense
defined in \cite[I,~5.7.1]{rg} is the same thing as an algebraic space (as defined above) 
equipped with a map to $S$.  

% \section*{Acknowledgments}
% \label{sec:acknowledgments}

% We thank Ofer Gabber for explaining the outline of the proof and
% Mathematisches Forschungsinsitut Oberwolfach for providing a
% delightful working environment.  We also 

\section{The excellent case}\label{exccase}

In this section we prove Theorem \ref{mainresult} when the algebraic
space $S$ is of finite presentation over an excellent noetherian
scheme (such as $\Spec \Z$).  
This case will be used to deduce the general case in \S\ref{approxsec}.

We will proceed by noetherian induction 
on $X$, so first we use deformation theory to show that the result for
$X_{\rm{red}} \rightarrow S$ implies the result for $X \rightarrow S$. Then 
we will be in position to begin the induction.  The base
space $S$ will remain fixed throughout the induction process.

\subsection{Reduction to the reduced case}\label{redredcase}

Suppose that Theorem \ref{mainresult} is proved for $X_{\rm{red}}
\rightarrow S$ with a fixed noetherian algebraic
space $S$.  Let us deduce the result for $X \rightarrow S$. 
We induct on the order of nilpotence
of the nilradical of $X$, so we may assume that there is a square-zero coherent ideal
sheaf $\mathscr{J}$ on $X$ such that the closed subspace $X_0 \hookrightarrow X$
defined by killing $\mathscr{J}$ admits an $S$-compactification, say 
$\sigma:X_0 \hookrightarrow \overline{X}_0$.  
Let $f_0:X_0 \rightarrow S$ and $f:X \rightarrow S$ be
the structure maps. 

By blowing up the noetherian $\overline{X}_0$ along a closed subspace
structure on $\overline X_0\setminus X_0$ (such as the reduced
structure) we can arrange that $\overline{X}_0 - X_0$ admits a
structure of effective Cartier divisor, so $\sigma$ is an affine
morphism.  Let us check that it suffices to construct a cartesian
diagram of algebraic spaces
\begin{equation}\label{squirt}
\xymatrix{
{X_0} \ar[r]^-{\sigma} \ar[d] & {\overline{X}_0} \ar[d] \\
{X} \ar[r] & {\overline{X}}}
\end{equation}
over $S$ in which the bottom arrow is an open immersion
and the right vertical arrow is a square-zero closed immersion
defined by a quasi-coherent ideal sheaf of ${\ms O}_{\overline{X}}$ whose
natural ${\ms O}_{\overline{X}_0}$-module structure is coherent.
In such a situation, since the square-zero ideal sheaf 
$\ker({\ms O}_{\overline{X}} \twoheadrightarrow {\ms O}_{\overline{X}_0})$
on $\overline{X}$ is coherent as an ${\ms O}_{\overline{X}_0}$-module,
$\overline{X}$ is necessarily of finite type over $S$ and thus is $S$-proper
(since $\overline{X}_0$ is $S$-proper).  We would therefore be done.

Rather than construct the geometric object $\overline{X}$ directly, we will construct
its structure sheaf ${\ms O}_{\overline{X}}$ as a square-zero thickening of
${\ms O}_{\overline{X}_0}$ in a manner that restricts over
the open $X_0 \subseteq \overline{X}_0$ to ${\ms O}_X$ viewed as a square-zero thickening
of ${\ms O}_{X_0}$.  The sufficiency of such a sheaf-theoretic approach is provided by an application
of Theorem \ref{sqzero} in the noetherian case (also see Remark \ref{remsqzero}), 
according to which the existence of diagram \eqref{squirt} is equivalent to extending 
the square-zero extension of
$f_0^{-1}({\ms O}_S)$-algebras
\begin{equation}\label{dagger}
0 \rightarrow \mathscr{J}  \rightarrow {\ms O}_X \rightarrow {\ms O}_{X_0} \rightarrow 0
\end{equation} 
on $(X_0)_{\et}$ to a square-zero extension of
$\overline{f}_0^{-1}({\ms O}_S)$-algebras
$$0 \rightarrow \overline{\mathscr{J}} \rightarrow \mathscr{A} \rightarrow
{\ms O}_{\overline{X}_0} \rightarrow 0$$
on $(\overline{X}_0)_{\et}$ in which the kernel $\overline{\mathscr{J}}$ is coherent as an ${\ms O}_{\overline X_0}$-module. 

Since $\sigma$ is affine we have ${\rm{R}}^1 \sigma_{\ast,\et}(\mathscr{J}) = 0$, so 
applying $\sigma_{\ast}$ to (\ref{dagger}) gives a  square-zero extension 
of $\overline{f}_0^{-1}({\ms O}_S)$-algebras 
$$0 \rightarrow \sigma_{\ast}(\mathscr{J}) \rightarrow
\sigma_{\ast}({\ms O}_X) \rightarrow \sigma_{\ast}({\ms O}_{X_0}) \rightarrow 0$$
on $(\overline{X}_0)_{\et}$ whose pullback along ${\ms O}_{\overline{X}_0}
\rightarrow \sigma_{\ast}({\ms O}_{X_0})$ is a square-zero extension  of
$\overline{f}_0^{-1}({\ms O}_S)$-algebras 
$$0 \rightarrow \sigma_{\ast}(\mathscr{J}) \rightarrow
\mathscr{B} \rightarrow {\ms O}_{\overline{X}_0} \rightarrow 0$$
in which the kernel $\sigma_\ast(\mathscr{J})$ is only quasi-coherent as an ${\ms O}_{\overline X_0}$-module.

By \cite[III,~Thm.~1.1,~Cor.~1.2]{knutson}, we know that $\sigma_{\ast}(\mathscr{J}) = \varinjlim \overline{\mathscr{J}}_{\alpha}$
where $\overline{\mathscr{J}}_{\alpha}$ ranges through the directed system of
${\ms O}_{\overline{X}_0}$-coherent subsheaves of $\sigma_{\ast}(\mathscr{J})$
satisfying $\overline{\mathscr{J}}_{\alpha}|_{(X_0)_{\et}} = \mathscr{J}$. 
Hence, 
our problem is reduced to proving bijectivity of the natural map
$$\varinjlim {\rm{Exal}}_{{\ms O}_S}({\ms O}_{\overline{X}_0},\mathscr{M}_i)
\rightarrow {\rm{Exal}}_{{\ms O}_S}({\ms O}_{\overline{X}_0}, \varinjlim
\mathscr{M}_i)$$ for any directed system $\mathscr{M}_i$ of
quasi-coherent ${\ms O}_{\overline{X}_0}$-modules. Here, as usual, we let
${\rm{Exal}}_{{\ms O}_S}(\mathscr{A}, \mathscr{M})$ denote the set of
isomorphism classes of square-zero extensions of an
$\overline{f}_0^{-1}({\ms O}_S)$-algebra $\mathscr{A}$ by an
$\mathscr{A}$-module $\mathscr{M}$ on $(\overline{X}_0)_{\et}$ (see for example 
\cite[$0_{\rm{IV}}$,~\S18]{ega} and \cite[III.1.1]{illusie}).

By \cite[Thm.~III.1.2.3]{illusie} 
(applied to the ringed topos $((\overline{X}_0)_{\et}, \overline{f}_0^{-1}({\ms O}_S))$), 
the cotangent complex ${\rm{L}}_{\overline{X}_0/S}$ of
${\ms O}_{\overline{X}_0}$-modules is bounded above and satisfies 
$${\rm{Exal}}_{{\ms O}_S}({\ms O}_{\overline{X}_0},\mathscr{M}) \simeq
{\rm{Ext}}^1_{{\ms O}_{\overline{X}_0}}({\rm{L}}_{\overline{X}_0/S}, \mathscr{M})$$
naturally in any ${\ms O}_{\overline{X}_0}$-module $\mathscr{M}$.  Moreover, 
by \cite[Cor.~II.2.3.7]{illusie} the complex ${\rm{L}}_{\overline{X}_0/S}$ has coherent homology modules. 
% (The formation of this complex is compatible with \'etale localization, by \cite[Prop.~II.2.30]{illusie},
% but we do not need this fact.)  

We are now reduced to showing that if $Z$ is any noetherian algebraic space (such as $\overline{X}_0$) and
$\mathscr{F}^{\bullet}$ is any bounded-above complex of ${\ms O}_{Z}$-modules with
coherent homology modules (such as ${\rm{L}}_{\overline{X}_0/S}$)  then the functor
${\rm{Ext}}^j_{{\ms O}_Z}(\mathscr{F}^{\bullet}, \cdot)$
on quasi-coherent ${\ms O}_{Z}$-modules commutes
with the formation of direct limits for every $j \in \Z$.  
This is a standard fact: one reduces first to the case
that $\mathscr{F}^{\bullet} = \mathscr{F}[0]$
for a coherent sheaf $\mathscr{F}$ on $Z$, and then uses the
the local-to-global Ext spectral sequence and the compatibility of
\'etale cohomology on the qcqs $Z$ with direct limits to reduce
to the case of affine $Z$ with the Zariski topology, which is handled by
degree-shifting in $j$. 

\subsection{Reduction to normal case}\label{normalred}

Now take $S$ to be an algebraic space of finite presentation over an excellent noetherian scheme.  In this section we prove the following result.

\begin{proposition}\label{P:normalred}
For our fixed $S$, if Theorem $\ref{mainresult}$ holds whenever $X$ is normal
  then it   holds in general. 
\end{proposition}
\begin{proof}
By noetherian induction on $X$ (with its fixed $S$-structure),
to prove Theorem \ref{mainresult} for $X \rightarrow S$ we may
assume that every proper closed subspace of $X$ admits an
$S$-compactification.  If $X$ is not reduced then \S\ref{redredcase}
may be used to conclude the argument, so we may assume that $X$ is
reduced.  

Let $\pi:\widetilde{X} \rightarrow X$ be the finite
surjective normalization, and let $Z \hookrightarrow X$ the 
closed subspace cut out by the coherent conductor ideal 
${\rm{Ann}}_X(\pi_{\ast}({\ms O}_{\widetilde{X}})/{\ms O}_X)$ of $\widetilde{X}$ over $X$. The open complement 
$X - Z$ is the maximal open subspace of $X$ over which $\pi$ is
an isomorphism.  We have $Z \ne X$ since $X$ is reduced (though $Z$ may be non-reduced).  By the
noetherian induction hypothesis, the separated finite type map $Z
\rightarrow S$ admits an $S$-compactification $\overline{Z}
\rightarrow S$.  Assuming that Theorem \ref{mainresult} is proved for
all normal $S$-separated algebraic spaces of finite type over $S$, so
$\widetilde{X}$ admits a compactification $\widetilde{X}^{-}$ over
$S$, let us see how to construct an $S$-compactification for $X$.

The idea is to reconstruct $X$ from $\widetilde{X}$ via a contraction
along the finite surjective map $\pi^{-1}(Z) \rightarrow Z$, and to
then apply an analogous such contraction to the $S$-compactification $\widetilde{X}^{-}$ 
of $\widetilde{X}$ (using the $S$-compactification 
$\overline{Z}$ in place of $Z$) to construct an $S$-compactification
of $X$.  We first record a refinement of a contraction theorem of
Artin.

\begin{theorem}[Artin]\label{artincontract}
Let $S$ be an algebraic space of finite presentation over an excellent 
noetherian scheme, and let
$$\xymatrix{
{Y'} \ar[r] \ar[d] & {X'} \ar[d] \\
{Y} \ar[r] & {S}}$$
be a commutative diagram of quasi-separated algebraic spaces locally of
finite type over $S$, with $Y' \rightarrow X'$ a closed immersion and
$Y' \rightarrow Y$ a finite surjective map. 

  \begin{enumerate}
  \item The pushout $X = Y \coprod_{Y'} X'$ exists in the category of
    algebraic spaces, it is quasi-separated and locally of finite type
    over $S$, and the pushout diagram
    \begin{equation}\label{yxyx}
      \xymatrix{
        Y' \ar[r] \ar[d] & X' \ar[d]^-{\pi} \\
        Y \ar[r] & X}
    \end{equation} 
    is cartesian, with $Y \rightarrow X$ a closed immersion and $X'
    \rightarrow X$ a finite surjection.  If $X'$ is $S$-separated
    $($resp. of finite presentation over $S$, resp. $S$-proper$)$ then
    so is $X$.
  \item The formation of this diagram $($as a pushout$)$ commutes with
    any flat base change on $X$ in the sense that if $X_1 \rightarrow
    X$ is a flat map of algebraic spaces then the cartesian diagram
    \begin{equation}\label{x1y1}
      \xymatrix{
        Y'_1 \ar[r] \ar[d] & X' _1\ar[d] \\
        Y_1 \ar[r] & X_1}
\end{equation}
obtained after base change is a pushout diagram in the category of
algebraic spaces.  In particular, the formation of $(\ref{yxyx})$ commutes with \'etale
base change on $X$.
\end{enumerate}
\end{theorem}

Before we prove Theorem \ref{artincontract}, we make some remarks.  

\begin{remark} 
By descent for morphisms it suffices to prove the result \'etale-locally
on $S$, so the case when $X'$ and $Y$ are $S$-separated (which is all
we will need in this paper) is easily reduced to the case when $X'$
and $Y$ are separated (over $\Spec \Z$).  
% In this case the result is
% asserted by Raoult for separated noetherian algebraic spaces (without
% any $S$ at all) in \cite[Prop.~1]{raoult2}, with some details left to
% his unpublished thesis, and for reasons implicit in that he did not 
% control properness aspects except when working with algebraic
% spaces of finite type over a field.

Note also that by taking $X_1 \rightarrow X$ to be $X - Y \rightarrow
X$, it follows that $\pi$ must restrict to an isomorphism over $X -
Y$.  This will also be evident from how $X$ is constructed.
\end{remark}

\begin{proof}[Proof of Theorem $\ref{artincontract}$] 
  By working \'etale-locally on $S$ we may assume that $S$ is an
  excellent noetherian scheme.  Before we address the existence, let
  us grant existence with $X$ quasi-separated and locally of finite
  type over $S$ and settle the finer structural properties at the end
  of part (1) for $X \rightarrow S$.  Since $\pi$ will be a finite
  surjection, the quasi-separated $X$ is necessarily of finite
  presentation over $S$ (equivalently, is quasi-compact) when $X'$ is
  of finite presentation over $S$.  Likewise, granting the existence
  result in general, if $X'$ is $S$-separated then the composite of
  the monomorphism $\Delta_{X/S}:X \rightarrow X \times_S X$ with the
  finite $\pi:X' \rightarrow X$ is proper, so $\Delta_{X/S}$ is proper
  and hence a closed immersion (i.e., $X$ is $S$-separated).  Finally,
  if $X'$ is $S$-proper then $X$ is at least $S$-separated, and so $X$
  is also $S$-proper since $\pi$ is a finite surjection.

  {\bf Step 1}.  We may now turn our attention to the existence
  problem.  As a first step, we establish a pushout property in a
  special case involving affine schemes.  Consider pair of ring maps
  $B_1\to B_1'$, $A_1'\to B_1'$, and form the commutative diagram
  of affine schemes
\begin{equation}\label{aba}
\xymatrix{\Spec(B'_1) \ar[r]^-{j'_1} \ar[d]_-{q_1} & {\Spec(A'_1)} \ar[d]^-{\pi_1} \\
{\Spec(B_1)} \ar[r]_{j_1} & \Spec(A_1)}
\end{equation}
with the fiber product ring $A_1 := B_1 \times_{B'_1} A'_1$.  Assume
that $j_1$ is a closed immersion, $\pi_1$ is a finite surjection, and that
the diagram is {\em cartesian} (i.e., the natural map $B_1
\otimes_{A_1} A'_1 \rightarrow B'_1$ is an isomorphism).  For example,
these hypotheses hold when $j'_1$ is a closed immersion and $q_1$ is a
finite surjection.  Indeed, if $B'_1 = A'_1/J'$ then the fiber product
ring $A_1 = B_1 \times_{A'_1/J'} A'_1$ satisfies $A_1/J = B_1$ for the
ideal $J = \{0\} \times J'$ in $A_1$, and since $J A'_1 = J'$ we have
$B_1 \otimes_{A_1} A'_1 = A'_1/ J A'_1 = A'_1/J' = B'_1$, so the
cartesian property holds.  Finally, the natural map $A_1 \rightarrow
A'_1$ is finite due to finiteness of $B_1 \rightarrow B'_1 = A'_1/J'$,
and in view of the cartesian property the surjectivity of $\pi_1$ is
reduced to the case of points in $\Spec A_1$ where $J$ doesn't vanish.
At such points $P$ we obtain a local isomorphism since the formation
$A_1$ is compatible with localization at an element of $A_1$ (such as
at any $(0,x) \not\in P$ with $x \in J'$).  

We shall prove that any such diagram (\ref{aba}) is a pushout in the
category of algebraic spaces.   Fix an algebraic space $T$.  We aim to
show that the natural map of sets
$$\rho:\Hom(\Spec A_1,T)\to\Hom(\Spec B_1,T)\times_{\Hom(\Spec
  B_1',T)}\Hom(\Spec A_1',T)$$
is a bijection.  Given an algebraic space $Q$ and an $A_1$-algebra
$C$, write $\underline Q(C)$ for the \'etale sheaf on $\Spec A_1$
whose value on an \'etale morphism $W\to\Spec A_1$ is $\Hom(W\otimes_{A_1}
C,Q)$.  Note that $\underline{Q}(A_1)$ is the restriction of $Q$
to the \'etale site of $\Spec A_1$.  We will prove that the map
$$\underline\rho:\underline T(A_1)\to\underline
T(B_1)\times_{\underline T(B_1')}\underline T(A_1')$$
is an isomorphism of sheaves (evaluation on $W = \Spec A_1$ recovers $\rho$). 

Let $(p,q):R\to U\times U$ be an \'etale presentation of $T$.  Suppose
we can show that the natural diagram
\begin{equation}\label{cheese}
\xymatrix{\underline R(B_1)\times_{\underline R(B_1')}\underline R(A_1') \ar@<-2pt>[r] \ar@<2pt>[r]
  & \underline U(B_1)\times_{\underline U(B_1')}\underline U(A_1')\ar[r] &
  \underline T(B_1)\times_{\underline T(B_1')}\underline T(A_1')}
\end{equation}
is a colimit diagram in the category of sheaves of sets.  There is a map of
diagrams
$$
\xymatrix{\underline R(B_1)\times_{\underline R(B_1')}\underline R(A_1')\ar@<-2pt>[r] \ar@<2pt>[r]
  & \underline U(B_1)\times_{\underline U(B_1')}\underline U(A_1')\ar[r] &
  \underline T(B_1)\times_{\underline T(B_1')}\underline T(A_1') \\
\underline R(A_1)\ar[u] \ar@<-2pt>[r] \ar@<2pt>[r]
  & \underline U(A_1)\ar[u]\ar[r] &
  \underline T(A_1)\ar[u]}
$$
in which we wish to show that the right-most vertical arrow is an
isomorphism.  If $U$ is a separated scheme (\'etale cover of $T$) then
$R$ is a separated scheme. Thus, if we can prove the result when $T$
is a separated scheme, we conclude that the left and middle vertical
arrows are isomorphisms, whence the right vertical arrow is an isomorphism 
because both rows are colimit diagrams of sheaves.  But when $T$ is
separated and $U$ is affine we have that $R$ is affine.  Applying the
same argument shows that it suffices to prove the result when $T$ is
an affine scheme.  

Note that the formation of (\ref{aba}) commutes with affine flat (e.g., affine
\'etale) base change on $A_1$ in the following sense.  Let $\Spec A_2
\rightarrow \Spec A_1$ be flat and define $A'_2 = A_2 \otimes_{A_1}
A'_1$ and similarly for $B'_2$ and $B_2$.  The natural exact sequence
of $A_1$-modules
$$0 \rightarrow B_1 \times_{B'_1} A'_1 \rightarrow B_1 \times A'_1 \rightarrow B'_1$$
remains exact after scalar extension by $A_1 \rightarrow A_2$, so the
natural map $A_2 \rightarrow B_2 \times_{B'_2} A'_2$ is an
isomorphism.  Since the morphism $\underline\rho$ is an isomorphism if
and only if it induces bijections when evaluated on all \'etale
$A_1$-algebras $A_1\to A_2$, we are done (when $T$ is affine) because we
know \eqref{aba} is a pushout in the category of affine schemes.

Thus, it suffices to show that \eqref{cheese} is a colimit diagram (of
sheaves).  
Suppose $(t_1,t_2)$ is a section of $\underline
T(B_1)\times_{\underline T(B_1')}\underline T(A_1')$ over an \'etale
$A_1$-scheme $W$.  Since all of the schemes in diagram \eqref{aba} are
finite over $\Spec A_1$, we know that their pullbacks along any strict
henselization of $A_1$ split as finite products of strictly henselian
local rings \cite[IV$_4$, 18.8.10]{ega}.  Thus, since $U\to T$
is an \'etale surjection (hence locally of finite presentation), after
replacing $W$ by an \'etale cover we can choose a lift
$\tilde t_1\in \underline U(B_1)(W)$ of $t_1$.  We
now have a commutative diagram
$$\xymatrix{W\otimes B_1'\ar[r]^-{\tilde t_1|_{B_1'}}\ar[d] & U\ar[d]\\
W\otimes A_1'\ar[r]_-{t_2}\ar[d]\ar@{-->}[ur] & T\\
W}$$
in which the top left vertical arrow is a closed immersion, the bottom
left vertical arrow is a finite surjection, and the
right vertical arrow is an \'etale surjection.  If $W\otimes A_1'$
were strictly henselian, a unique dotted arrow in the commutative diagram would
exist.  Thus, for any strictly henselian local $W$-scheme $H\to
W$, the restricted diagram
$$\xymatrix{H\otimes B_1'\ar[r]\ar[d] & U\ar[d]\\
H\otimes A_1'\ar[r]\ar[d]\ar@{-->}[ur] & T\\
H}$$
admits a unique dotted lift.  Since $U\to T$ is locally of finite
presentation, it follows that there is an \'etale surjection $W'\to W$
such that there is a dotted lift in the resulting diagram
$$\xymatrix{W'\otimes B_1'\ar[r]\ar[d] & U\ar[d]\\
W'\otimes A_1'\ar[r]\ar[d]\ar@{-->}[ur] & T\\
W'}$$
completing the proof that the right arrow of diagram \eqref{cheese} is
a surjection of \'etale sheaves.

Now suppose $(u_1,v_1)$ and $(u_2,v_2)$ are two sections of
$\underline U(B_1)\times_{\underline U(B_1')}\underline U(A_1')$ over
$W\to\Spec A_1$ with the same image in $\underline
T(B_1)\times_{\underline T(B_1')}\underline T(A_1')$.  Since $R\to
U\times_T U$ is an isomorphism, we know that there are unique sections
$r_1\in \underline R(B_1)(W)$ and $r_2\in \underline R(A_1')(W)$ such
that $p(r_i)=u_i$ and $q(r_i)=v_i$.  Moreover, these sections restrict
to sections $r\in \underline R(B_1')(W)$ such that
$p(r)=u_1|_{B_1'}=u_2|_{B_1'}$ and $q(r)=v_1|_{B_1'}=v_2|_{B_1'}$.
Since any such $r$ is unique, we see that $(r_1,r_2)$ defines the
unique section of $\underline R(B_1)\times_{\underline
  R(B_1')}\underline R(A_1')$ over $W$ carrying $(u_1,v_1)$ to
$(u_2,v_2)$.  This completes the proof that \eqref{cheese} is a
colimit diagram and thus the proof that \eqref{aba} is a pushout
diagram in the category of algebraic spaces.

{\bf Step 2}.   The main work
is to handle the case when $X'$ (and hence $Y'$ and $Y$) is of finite
presentation over $S$, with $X$ constructed as also finitely presented
over $S$; the existence result more generally will then be obtained by
simple gluing arguments.  Thus, we now assume (until said otherwise)
that $X'$, $Y'$, and $Y$ are of finite presentation over $S$.  In this
situation, the existence of the pushout in the category of quasi-separated
algebraic spaces is \cite[Thm.~6.1]{artin} when
$S$ is of finite presentation over a field or excellent Dedekind domain, and
the construction in its proof gives that (i) $X$ is of finite
presentation over $S$, (ii) the diagram (\ref{yxyx}) is cartesian, (iii) $Y
\rightarrow X$ is a closed immersion, and (iv) $\pi:X' \rightarrow X$
is proper and restricts to an isomorphism over the open subspace $X -
Y$.  Since $\pi$ is clearly quasi-finite (by the cartesian property),
it must also be finite.  

Artin assumed $S$ is of finite presentation over a
field or excellent Dedekind domain (rather than over an
arbitrary excellent noetherian scheme) because his criterion for
a functor to be an algebraic space was originally proved only for such
$S$.  By \cite[Thm.~1.5]{cd} Artin's proof of that criterion works for any
excellent noetherian $S$, so likewise the above conclusions hold in
such generality.  The pushout constructed in the finitely presented
case has only been shown
to be a pushout in the category of quasi-separated algebraic spaces,
as this is the situation considered by Artin.  To establish the pushout property
relative to maps to any algebraic space, the key step is to check, as we shall now do
(assuming $X'$, $Y'$, and $Y$ are finitely presented over $S$),  
that the pushout property of $X$ is preserved by any quasi-separated \'etale localization on $X$.  
Artin's construction of the quasi-separated pushout $X$ is via an indirect algebraization process, so to be rigorous
will require some care.  If $X_1 \rightarrow X$ is a quasi-separated
 \'etale map, to prove that the $X_1$-pullback diagram is
a pushout it is enough (by \'etale descent for morphisms) to check this property
\'etale-locally on $X_1$.  More specifically, it suffices to treat the case when $X_1 = \Spec A_1$ is affine. 

Consider the resulting pullback diagram (\ref{x1y1}) which consists of
affine schemes, say with $Y_1 = \Spec B_1$, $X'_1 = \Spec A'_1$, and $Y'_1 = 
\Spec B'_1$.  We claim that the natural map $\theta_1:A_1 \rightarrow B_1 \times_{B'_1} A'_1$
is an isomorphism.  Let $J_1 = \ker(A_1 \twoheadrightarrow B_1)$,
so $B'_1 = A'_1/J_1 A'_1$.  Since $A'_1$ is $A_1$-finite, $\theta_1$ is at least
finite.  Also, $\Spec(\theta_1)$ is clearly an isomorphism over the open complement of
$\Spec B_1 = \Spec(A_1/J_1)$ in $\Spec A_1$.  Hence, to prove that $\theta_1$
is an isomorphism it suffices to show that the induced map $\widehat{\theta}_1$
between $J_1$-adic completions is an isomorphism.

Write $\widehat{A}_1$ and $\widehat{A}'_1$ to denote the $J_1$-adic and
$J_1 A'_1$-adic completions respectively, and let the formal algebraic space 
$\mathfrak{X}'$ denote the formal completion of $X'$ along $Y'$.
The \'etale map $\Spec B_1 \rightarrow Y$ has pullback  along $Y' \rightarrow Y$ identified with 
$\Spec B'_1 \rightarrow Y'$, and  (using Proposition \ref{redequiv}) the unique lifting of this latter
\'etale map to a formal noetherian algebraic space formally
\'etale over $\mathfrak{X}'$ is uniquely identified with ${\rm{Spf}}(\widehat{A}'_1) \rightarrow
\mathfrak{X}'$.  Artin's construction of $X$ identifies $\widehat{A}_1$ 
with the ring-theoretic fiber product over $B'_1$
of $B_1$ against the coordinate ring $\widehat{A}'_1$ of this formal lifting, 
which is to say that the natural map
$\widehat{A}_1 \rightarrow B_1 \times_{B'_1} \widehat{A}'_1$ is
an isomorphism.   This isomorphism is the map $\widehat{\theta}_1$, so
$\widehat{\theta}_1$ is an isomorphism and hence $\theta_1$ is an isomorphism. 

With $\theta_1$ now shown to be an isomorphism, the verification of the quasi-separated 
pushout property after
base change to $X_1$ is a special case of a general pushout property for ring-theoretic fiber products that 
we settled in Step 1 (since Artin's pushout diagram is also cartesian). 
Since any quasi-separated \'etale base change on $X$ has now been shown
to yield a pushout diagram in the category 
of quasi-separated algebraic spaces, to prove that 
Artin's quasi-separated pushout is actually a pushout in the category
of all algebraic spaces we can use \'etale descent for morphisms to formally reduce
to the special case when $X'$, $Y$, and $Y'$ are affine.  In this affine case we can form the diagram
as in (\ref{aba}) except with $q_1$ a finite surjection and $j'_1$ a closed immersion
(and $A_1 := B_1 \times_{B'_1} A'_1$).  As we saw early in Step 1, 
in such cases necessarily 
$j_1$ is a closed immersion, $\pi_1$ is a finite surjection, and the diagram is cartesian.
Hence, by Step 1 the affine scheme $\Spec(A_1)$ is a pushout in the category of all algebraic spaces.
This affine scheme is quasi-separated, so it must coincide with the pushout already constructed in the category of
quasi-separated algebraic spaces.  Hence, this latter pushout is also a pushout in the category of all
algebraic spaces.  The proof of compatibility with quasi-separated \'etale base now applies {\em verbatim}
to arbitrary \'etale base change. 

{\bf Step 3}. 
We continue to assume that $X'$, $Y'$, and $Y$ are finitely presented over $S$,
and now we improve upon Step 2 by showing that 
the formation of $X$ is compatible with any flat base change $X_1 \rightarrow X$
(in the sense that the $X_1$-pullback diagram is a pushout in the category of algebraic spaces).
By \'etale descent for morphisms, coupled with the established compatibility with \'etale base change on
$X$, we are reduced to the case when $X$ is an affine scheme
(so $Y$, $Y'$, and $X'$ are also affine) and $X_1$ is affine.  Say $X = \Spec A$, $Y = \Spec B$, 
$X' = \Spec A'$, $Y' = \Spec B'$, and $X_1 = \Spec A_1$.   In the commutative diagram of 
noetherian rings
$$\xymatrix{
B' & A'  \ar[l] \\
B \ar[u] & A \ar[l] \ar[u] }$$
the vertical maps are finite, the horizontal maps are surjective, 
$B' = A' \otimes_A B$, and (as in Step 2)
the natural map of rings $\phi:A \rightarrow B \times_{B'} A'$
is an isomorphism.   
Let $A'_1 := A_1 \otimes_A A'$, and similarly for $B_1$ and $B'_1$. 
By the same calculation with flat scalar extension as in Step 1, the natural map 
$A_1 \rightarrow B_1 \times_{B'_1} A'_1$ is an isomorphism.
Thus, in view of the general pushout result proved in Step 1, we have established 
compatibility with any flat base change on $X$ (when $X'$, $Y'$, and $Y$ are finitely presented over $S$). 

Theorem \ref{artincontract} has now been proved when $X'$, $Y'$, and $Y$ are
finitely presented over $S$.  It remains 
to handle the existence and compatibility with flat base change of a quasi-separated
pushout $X$ when the quasi-separated $X'$, $Y'$, and $Y$ are merely locally
of finite type over 
$S$.  (Note that finite presentation is the same as finite type for algebraic spaces
quasi-separated over $S$.) 

The proof of existence will proceed by a gluing method, and the compatibility with flat
base change on $X$ in general will follow immediately from the gluing construction of $X$
and such compatibility in the finitely presented case.  The key point is that the construction
given above in the
finitely presented case is well-behaved with respect to Zariski localization on the source.  More
precisely, we have:

\begin{lemma}\label{324} In the setup of Theorem $\ref{artincontract}$, 
let $V' \subseteq X'$ and $W \subseteq Y$ be open subspaces such that
$\pi^{-1}(W) = Y' \cap V'$ as open subspaces of $Y'$.
Then the natural map of pushouts 
$$W \coprod_{\pi^{-1}(W)} V' \rightarrow Y \coprod_{Y'} X'$$
is an open immersion.
\end{lemma}

Since we have only proved Theorem \ref{artincontract} 
in the finitely presented case (over $S$), the 
lemma can only be proved at this point in such cases (i.e., 
$X'$, $Y'$, $Y$, $V'$, and $W$ are finitely presented over $S$).  However, 
the proof goes the same way in this case as it will in the general case, so
we write one argument below that is applied first in the finitely presented case
and then (after general existence as in Theorem \ref{artincontract}(1) is proved) 
in the general case.  

\begin{proof}
  By Theorem \ref{artincontract}(2) the formation of the pushout is
  compatible with flat base change on $Y \coprod_{Y'} X'$, so by working
  \'etale-locally on this pushout we may reduce to the case when it is
  affine.  In other words, we have $X' = \Spec(A')$, $Y' = \Spec(B')$,
  $Y = \Spec(B)$, and $Y \coprod_{Y'} X' = \Spec(A)$ with $A = B
  \times_{B'} A'$, $B \rightarrow B'$ finite, $A' \rightarrow B'$
  surjective, and $B' = B \otimes_A A'$ (so $A \rightarrow A'$ is
  finite and $A \rightarrow B$ is surjective).  In particular, $X'
  \rightarrow X$ is an isomorphism over $X - Y$.

The condition $\pi^{-1}(W) = Y' \cap V'$ implies that
$V' = \pi^{-1}(\pi(V'))$, so since $\pi$ is a surjective finite map
we see that $V := \pi(V')$ is an open subset of $X = \Spec(A)$
with complement $Y - W$. 
Giving $V$ the open subscheme structure, we want the commutative diagram 
$$\xymatrix{
\pi^{-1}(W) \ar[d] \ar[r] & V' \ar[d] \\
W \ar[r] & V}$$
to be a pushout.  That is, we want the natural map from the algebraic space
$P := W \coprod_{\pi^{-1}(W)} V'$ to the scheme $V$ to be an isomorphism.
We may work Zariski-locally on $V$ due to the flat base change compatibility of pushouts,
so we may assume $V = \Spec(A_a)$ for some $a \in A = B \times_{B'} A'$.
Writing $a = (b, a')$ where $b$ and $a'$ have the same image $b'$ in $B'$, 
clearly $A_a = B_b \times_{B'_{b'}} A'_{a'}$.
But $\Spec(A'_{a'})$ is the preimage $V'$ of $V$ in $X'$ and $\Spec(B_{b})$
is the preimage $W$ of $V$ in $Y$, so the isomorphism property for $P \rightarrow V$
is reduced to the affine cases for which the pushout has already been shown to be given
by a ring-theoretic fiber product. 
\end{proof}

To complete the proof of existence of $Y \coprod_{Y'} X'$ as a quasi-separated algebraic
space locally of finite type over $S$ in general, let $\{U_i\}$ be a Zariski-open
covering of $Y$ by quasi-compact opens, and let
$\{U'_i\}$ be the pullback cover of $Y'$.  Each $U'_i$ has the form
$U'_i = Y' \cap V'_i$ for a quasi-compact open subspace
$V'_i \subseteq X'$.   Thus, we can form pushouts
$V_i := U_i \coprod_{U'_i} V'_i$ of finite presentation over $S$.
Define $U_{ij} = U_i \cap U_j$, $U'_{ij} = U'_i \cap U'_j$, and
$V'_{ij} = V'_i \cap V'_j$.   We may form the pushout
$V_{ij} := U_{ij} \coprod_{U'_{ij}} V'_{ij}$ and  by Lemma \ref{324}
the natural maps $V_{ij} \rightarrow V_i$ and $V_{ij} \rightarrow V_j$ are quasi-compact 
open immersions.  It is trivial to check the triple overlap compatibility, and so
we may glue the $V_i$'s to obtain a quasi-separated algebraic space $V$ locally of finite type over $S$
equipped with a closed immersion $U \hookrightarrow V$ and a finite surjection
$V' \rightarrow V$ with respect to which $V$ satisfies
the universal property of $Y \coprod_{Y'} V'$ where $V' = \cup V'_i$.  
Either by construction or flat base change compatibility (relative to 
$V - Y \rightarrow V$), the finite surjection $V' \rightarrow V$ restricts to an isomorphism over
$V - Y$.  Hence, we may glue $V$ and $X'$ along the common open subspace
$V - Y \simeq \pi^{-1}(V - Y) = V' - Y' \subset X' - Y'$ inside of $V$ and $X'$.   This gluing is the required $X$
and satisfies all of the desired properties. 
\end{proof}

As an application of Theorem \ref{artincontract}, we can now give a pushout method to
reconstruct certain reduced algebraic spaces from their normalization and their non-normal locus.  

\begin{corollary}\label{lemma9} Let $X$ be a reduced
quasi-separated algebraic space locally of finite type over
an excellent scheme, and let $\pi:\widetilde{X} \rightarrow X$ denote the normalization.
Let $j:Z \hookrightarrow X$ be the closed subspace cut out by the
conductor of $\widetilde{X}/X$.  Let $Y = Z \times_X \widetilde{X}$.  The natural map
$Z \coprod_Y \widetilde{X} \rightarrow X$ is an isomorphism.  
\end{corollary}

\begin{proof} By Theorem \ref{artincontract}, the formation of the
  pushout $Z \coprod_Y \widetilde{X}$ commutes with \'etale
  localization on the pushout, and in particular with \'etale
  localization on $X$.  Since the formation of the conductor is
  \'etale-local on $X$, it suffices to treat the case when $X = \Spec
  A$ is affine, so $\widetilde{X} = \Spec \widetilde{A}$ for the
  normalization $\widetilde{A}$ of $A$, and $Z = \Spec(A/J)$ and $Y =
  \Spec(\widetilde{A}/J)$ for the conductor ideal $J\subset
  A$ (so that $J$ is also an ideal of $\widetilde A$).  The argument
  in Step 1 of the proof of Theorem \ref{artincontract} shows that in
  this case the pushout is identified with $\Spec((A/J)
  \times_{\widetilde{A}/J} \widetilde{A})$, so our
  problem is to prove that the natural map
$$h:A \rightarrow C := (A/J) \times_{\widetilde{A}/J} \widetilde{A}$$
is an isomorphism.   

Since $\widetilde{A}$ is $A$-finite, it is obvious that 
$h$ is finite.  Also, $h$ is injective since $A \rightarrow \widetilde{A}$ is injective
(as $A$ is reduced).   On the other hand, since $J$ is an ideal of
both $A$ and $\widetilde A$, it is a trivial calculation (using the
explicit description of the fiber product of rings) that $h$ is surjective.
\end{proof}

Now we return to the setup with noetherian induction preceding the statement of
Theorem \ref{artincontract}.   Just as Corollary \ref{lemma9} reconstructs
$X$ from $\widetilde{X}$ by contracting along the canonical finite surjective map $\pi:Y =
Z \times_X \widetilde{X} \rightarrow Z$,
we aim to construct an $S$-compactification of $X$ by contracting a suitable choice of 
$\widetilde{X}^{-}$ along a finite surjective map 
$\overline{\pi}:\overline{Y} \rightarrow \overline{Z}$, where $\overline{Y}$ is the closure
of $Y$ in $\widetilde{X}^{-}$.  The first step is to construct $\overline{\pi}$ extending $\pi|_Y$.  

\begin{lemma}\label{lemma11}  For a suitable choice of schematically dense
open immersions
$\widetilde{X} \hookrightarrow \widetilde{X}^{-}$ and $Z \hookrightarrow \overline{Z}$
over $S$ into $S$-proper algebraic spaces, the schematic closure $\overline{Y}$ of
$Y := Z \times_X \widetilde{X}$ in $\widetilde{X}^{-}$ admits
a finite surjective $S$-map $\overline{\pi}:\overline{Y} \rightarrow \overline{Z}$
which restricts to $\pi:Y \rightarrow Z$ over the open subspace $Z \subseteq \overline{Z}$.
\end{lemma}

\begin{proof}
We make an initial choice of $S$-compactifications
$\widetilde{X} \hookrightarrow \widetilde{X}_1^{-}$ and $Z \hookrightarrow
\overline{Z}_1$, which we may and do arrange to be
schematically dense, and we define $\overline{Y}_1$ to be the schematic closure
of $Y$ in $\widetilde{X}_1^{-}$.  Let $\overline{Y}'$ denote the
$S$-proper schematic closure of $Y$ in $\overline{Y}_1 \times_S \overline{Z}_1$.  
The natural map $q':\overline{Y}' \rightarrow \overline{Y}_1$ restricts to 
an isomorphism over the open subspace $Y \subseteq \overline{Y}_1$ because
$Y \rightarrow Y \times_S \overline{Z}_1$ is a closed immersion
(as it is the graph of an $S$-map $Y \rightarrow Z \hookrightarrow \overline{Z}_1$
to an $S$-separated target).   Likewise, due to the definition of
$\overline{Y}'$ as a scheme-theoretic closure, the natural proper $S$-map
$\pi':\overline{Y}' \rightarrow \overline{Z}_1$ restricts to $\pi$ over 
the open subspace $Z \subseteq \overline{Z}_1$ because the
monomorphism $Y \rightarrow \overline{Y}'
\times_S Z$ is a closed immersion (as it is finite, due to finiteness of $\pi:Y \rightarrow Z$).

Now we use some results of Raynaud and Gruson concerning the use of blow-ups
of algebraic spaces to ``improve'' properties of morphisms. 
Since the proper map $q':\overline{Y}' \rightarrow \overline{Y}_1$ 
restricts to an isomorphism over the open subspace
$Y \subseteq \overline{Y}_1$, by \cite[I,~5.7.12]{rg} there is a blow-up 
$q'':\overline{Y}'' \rightarrow \overline{Y}'$ with center disjoint from $Y$
such that $q' \circ q'':\overline{Y}'' \rightarrow \overline{Y}_1$
is a blow-up with center disjoint from $Y$.  
(Blow-ups away from $Y$ are easier to work with than general morphisms that are
isomorphisms over $Y$, since we can focus attention on
the center of the blow-up.)  

Let $\pi'' = \pi' \circ q'':\overline{Y}'' \rightarrow \overline{Z}_1$
denote the natural composite map, so this restricts to 
the finite map $\pi$ over $Z \subseteq \overline{Z}_1$.  Hence, by
\cite[I,~5.7.10]{rg}, there is a blow-up $g:\overline{Z} \rightarrow
\overline{Z}_1$ with center disjoint from $Z$
such that the strict transform $\varphi:\overline{Y} \rightarrow \overline{Y}''$ 
of $\pi''$ with respect to $g$ has $\overline{Y}$ finite over
$\overline{Z}$.  (Note that $\varphi$ is a blow-up of
$\overline{Y}''$ with center disjoint from $Y$; see 
\cite[Lemma~1.1]{conrad} for a proof which adapts immediately to the case
of algebraic spaces. Also, since $Z$ is a schematically dense open in $\overline{Z}_1$, it
is a schematically dense open in $\overline{Z}$.)  
By construction, the finite map $\overline{\pi}:\overline{Y} \rightarrow \overline{Z}$ restricts
to $\pi$ over $Z$, and 
$Y$ and $Z$ are respectively schematically dense open subspaces 
in the $S$-proper $\overline{Y}$ and $\overline{Z}$. 
Since $\pi$ is surjective, it follows from the schematic density of $Z$
in $\overline{Z}$ that $\overline{\pi}$ is surjective. 

The composite map $(q' \circ q'') \circ \varphi:\overline{Y} \rightarrow \overline{Y}_1$
is a composite of blow-ups with center disjoint from $Y$, so by
\cite[I,~5.1.4]{rg} (cf. \cite[Lemma~1.2]{conrad} for a more detailed proof, 
which carries over to the case of algebraic spaces  with the help of
\cite[I,~5.7.8]{rg}) it is
itself a blow-up along a closed subspace
$C \subseteq \overline{Y}_1$ disjoint from $Y$.  
Since $\widetilde{X} \cap \overline{Y}_1 = Y$ as open subspaces 
of $\overline{Y}_1$, when $C$ is viewed as a closed subspace of $\widetilde{X}_1^{-}$
it is disjoint from the open subspace $\widetilde{X}$. 
Thus, since $\widetilde{X}$ is a schematically dense open in
$\widetilde{X}^{-}_1$, the blow-up $\widetilde{X}^{-} := {\rm{Bl}}_C(\widetilde{X}^{-}_1)$ is an $S$-proper
algebraic space naturally containing $\widetilde{X}$ as a schematically
dense open subspace over $S$.   Exactly as in the case of schemes
(see \cite[Lemma~1.1]{conrad}), 
the blow-up $\overline{Y} = {\rm{Bl}}_C(\overline{Y}_1)$ is naturally a closed subspace 
of the blow-up $\widetilde{X}^{-}$.   Hence, $\overline{Y}$ must be the
schematic closure of $Y$ in $\widetilde{X}^{-}$ (since $Y$ is a schematically dense open
in $\overline{Y}_1$). 
Since $\overline{Y}$ was constructed to admit the desired $\overline{\pi}$, we are done. 
\end{proof}

Using $S$-compactifications as in Lemma \ref{lemma11}, define the pushout algebraic space
$$\overline{X} := \overline{Z} \coprod_{\overline{Y}} \widetilde{X}^{-}.$$
This is a pushout of the sort considered in Theorem \ref{artincontract}.
By Theorem \ref{artincontract}, 
$\overline{X}$ is $S$-proper.  By Corollary \ref{lemma9}
and the functoriality of pushouts, there is a natural $S$-map
$$j:X \simeq Z \coprod_Y \widetilde{X} \rightarrow \overline{Z} \coprod_{\overline{Y}}
\widetilde{X}^{-} =: \overline{X}.$$
Thus, to complete the reduction of the proof of Theorem \ref{mainresult} over 
$S$ to the case
when $X$ is normal, it suffices to prove that $j$ is an open immersion.  
Since $\overline{\pi}^{-1}(Z) = Y = \overline{Y} \cap X$ as open subspaces of $\overline{Y}$, 
this is a special case of Lemma \ref{324}. 

This completes the proof of Proposition \ref{P:normalred}.
\end{proof}
 
To prove the general case of Theorem \ref{mainresult} we need a few
results on quotients.

 \subsection{Group quotients}\label{gabberlemma}

The proof of Theorem \ref{mainresult} for normal $X$
(and $S$ a noetherian scheme) in \cite{raoult} rests on a group quotient
result that we shall find useful for other purposes, so we now wish to record
it.   Rather generally, if $X'$ is a quasi-separated algebraic space equipped with
an action by a finite group $G$, we define the {\em quotient} $X'/G$
(if it exists) to be an
initial object $X' \rightarrow X'/G$ in the category of quasi-separated algebraic spaces equipped
with a $G$-invariant map from $X'$ provided that (in addition) the map of sets 
$X'(k)/G \rightarrow (X'/G)(k)$ is bijective for all algebraically closed
fields $k$.  (It makes sense to drop the quasi-separated hypotheses, but we only know an existence
result for such quotients using a universal mapping property within
the category of quasi-separated algebraic spaces.) 
Note that  if $X'$ is reduced and $X'/G$ exists then 
$X'/G$ must be reduced since $X' \rightarrow (X'/G)_{\rm{red}}$ is easily shown
to satisfy the same universal property.   Such quotients are
useful for relating construction problems for normal noetherian algebraic spaces to analogous
problems for normal noetherian schemes, due to the following result. 

\begin{proposition}\label{prop1}
Let $X$ be a $($non-empty$)$ connected normal noetherian algebraic
space. There exists a connected normal noetherian scheme $X'$
equipped with a right action by a finite group $G$ and
a finite $G$-invariant map $\pi:X' \rightarrow X$ such that 
$\pi$ is finite \'etale $G$-torsor over a dense open subspace of $X$ and 
exhibits $X$ as $X'/G$. 
$($In particular, $X'/G$ exists.$)$
\end{proposition}

This result is \cite[Cor.~16.6.2]{lmb}, and it is also proved in \cite[Cor.]{raoult};
in both of these references, the existence of 
$X'/G$ is part of the assertion.  For the convenience of the reader, we give a proof 
after we first 
discuss the general existence problem for 
$X'/G$ when one does not have a candidate for this quotient already in hand.  
Such an existence result is required for applications to compactification.
Deligne proved the existence of $X'/G$ when $X'$ is a separated 
algebraic space; this is recorded in Theorem \ref{delthm} below and is not used in this paper. 
% and this seems to have been used in
% \cite{raoult}.
We wish to avoid such separatedness hypotheses
on total spaces, only ever assuming separatedness for morphisms (if at all).

The best approach we know for existence results
for quotients $X'/G$ is to use the work of Keel and Mori \cite{km} (or its generalizations; see 
\cite{rydh-qt})
on coarse moduli spaces for Artin stacks.  This allows one to relate 
the quotient $X'/G$ to the Deligne--Mumford stack $[X'/G]$.
It is therefore convenient to now recall the definition and main existence
theorem for coarse moduli spaces of quasi-separated Artin stacks.  

 If $\mathscr{X}$ is a quasi-separated Artin stack then
a {\em coarse moduli space} is a morphism
$\pi:\mathscr{X} \rightarrow X$ to a quasi-separated algebraic space $X$ such that it is initial in the category of
maps from $\mathscr{X}$ to quasi-separated algebraic spaces and the map of sets 
$(\mathscr{X}(k)/\simeq) \rightarrow X(k)$ is bijective
for every algebraically closed field $k$.  
It was proved by Keel and Mori \cite{km}
that there exists a coarse moduli space $X$ 
whenever $\mathscr{X}$ is of finite presentation (and hence quasi-separated) 
over a locally noetherian scheme $S$
and the inertia stack $I_S(\mathscr{X}) = \mathscr{X} \times_{\mathscr{X} \times_S \mathscr{X}}
\mathscr{X}$ is $\mathscr{X}$-finite (under either projection map).  Moreover, it is
proved there that the following additional properties hold in such cases: $\pi$ is proper and quasi-finite,
$X$ is of finite presentation over $S$, $X$ is $S$-separated if
$\mathscr{X}$ is $S$-separated, and the formation of $\pi$ commutes
with any quasi-separated flat base change morphism $X' \rightarrow X$ that is locally of finite type (i.e., 
$X'$ is the coarse moduli space of the quasi-separated Artin stack $\mathscr{X} \times_X X'$).   These results were generalized in \cite{conrad3} and \cite{rydh-qt}, where the noetherian hypotheses were eliminated, but the above results for $\mathscr{X}$ finitely presented over a locally noetherian scheme 
are enough for what
we need.   Note that a special case of
the compatibility with quasi-separated flat base change locally of finite type 
is that the formation of the coarse moduli space $X$ is
compatible with quasi-separated \'etale base change on $X$.

By using the universal properties of coarse moduli spaces and quotient stacks, one easily proves:

\begin{lemma}\label{maxlemma} Let $Y'$ be a quasi-separated
algebraic space equipped with an action by a finite group
$H$.  The quotient $Y'/H$ exists if and only if the Deligne--Mumford stack $[Y'/H]$
admits a coarse moduli space $Q$, in which case the 
natural map $Y'/H \rightarrow Q$ is an isomorphism. 
\end{lemma} 

We shall be interested in the special case of Lemma \ref{maxlemma} when $Y'$ is separated  of finite
presentation (even finite) over a noetherian algebraic space $S$ and $H$ acts on
$Y'$ over $S$.  In this case the quotient stack $[Y'/H]$ is of finite
presentation over $S$ with diagonal $\Delta_{[Y'/H]/S}$ that is
separated (as it is a subfunctor of the separated Isom-functor between
pairs of $H$-torsors over $S$-schemes), and the projections
$I_S([Y'/H]) \rightrightarrows [Y'/H]$ are finite (using the fact that $Y'$ is separated).
  Hence, if $S$ is a scheme then by \cite{km} the
quotient $Y'/H$ does exist as an algebraic space of finite
presentation over $S$, and the map $Y' \rightarrow Y'/H$ is a finite
surjection because $Y' \rightarrow [Y'/H]$ is an $H$-torsor and
$[Y'/H] \rightarrow Y'/H$ is a proper quasi-finite surjection. 

 In
particular, since  $Y'$ is $S$-separated the quotient $Y'/H$ is also 
$S$-separated (as could also be deduced from $S$-separatedness of
$[Y'/H]$), so if $Y'$ is $S$-proper then $Y'/H$ is also $S$-proper.
The same conclusions hold if $S$ is merely a noetherian algebraic
space rather than a noetherian scheme.  Indeed, since quotients by
\'etale equivalence relations always exist in the category of
algebraic spaces (Corollary \ref{qtspace}), the \'etale-localization
compatibility of the formation of coarse spaces in the setup of \cite{km} allows us to work \'etale-locally
over $S$ (and to thereby reduce to the case when $S$ is a noetherian scheme) for
the existence result as well as for the finer asserted properties of
the quotient over $S$.  The following is a special case.

\begin{example}\label{exqt} If $X'$ is a noetherian algebraic space equipped with 
an action by a finite group $G$ and there is a $G$-invariant finite map $X' \rightarrow S$ 
to a noetherian algebraic space $S$ then $X'/G$ exists
and the map $X'/G \rightarrow S$ is proper and quasi-finite, hence finite. 
\end{example}

Our proof of Proposition \ref{prop1} will use the irreducible component decomposition
for locally noetherian algebraic spaces, and we refer the reader to
Proposition \ref{irredcomp} for a general discussion of this decomposition 
(avoiding the local separatedness hypotheses imposed in \cite[II,~\S8.5]{knutson}).

One final issue we address before taking up the proof of Proposition \ref{prop1}
is normalization in function field extensions for quasi-separated algebraic spaces.  
Let $X$ be a
reduced and irreducible 
locally noetherian algebraic space (so $X$ is quasi-separated;
see Definition \ref{noethdef}).  Let $\eta$ be the unique generic point of $X$,
so $X$ contains a dense open subspace around $\eta$ that is a scheme. 
The {\em function field} $k(X)$ is the henselian local ring of $X$ at $\eta$,
or more concretely it is the common function field of any (necessarily reduced
and irreducible) open scheme neighborhood of $\eta$ in $X$, so
there is a canonical map $\Spec k(X) \rightarrow X$.

By using an open scheme neighborhood of $\eta$ in $X$ we see that 
for any \'etale map $V \rightarrow X$ from a scheme, 
the pullback $V_{\eta}$ over $\Spec k(X)$ is the scheme of
generic points of $V$.  Hence, for any finite reduced $k(X)$-algebra $K$, 
the pullback $V_K$ along $\Spec K \rightarrow X$ 
is an \'etale $K$-scheme
that is a finite flat cover of the scheme of generic points of $V$.
The normalization $V'$ of $V$ in $V_K$ therefore makes sense as an affine 
surjective map $V' \rightarrow V$ that is finite when either $X$ is 
locally of finite type over an excellent scheme or $X$ is normal and
$K/k(X)$ is separable, 
and the scheme of generic points of $V'$ is identified with $V_K$.  

The
formation of the normalization $V'$ is \'etale-local on $V$, so by \'etale descent the affine surjective maps 
$V' \rightarrow V$ (for varying schemes $V$ \'etale over $X$) uniquely descend to a common 
affine surjective map of algebraic spaces $\pi:X' \rightarrow X$.  In particular,
$X'$ is normal and ${\ms O}_X \rightarrow \pi_{\ast}({\ms O}_{X'})$ is injective.  We call
$X' \rightarrow X$ the {\em normalization of $X$ in $K/k(X)$}. 
In the special case $K = k(X)$ we call $X'$ the {\em normalization of $X$}. 
More generally, if $X$ is a reduced locally noetherian algebraic space
that has finite many irreducible components $\{X_i\}$ then we can define
the affine surjective normalization $X' \rightarrow X$ 
of $X$ in any finite reduced faithfully flat algebra over
$\prod k(X_i)$. 

Now assume that the reduced
and irreducible locally noetherian $X$ is locally of finite type over
an excellent scheme or that $X$ is normal and $K/k(X)$ is separable, so
the normalization $\pi:X' \rightarrow X$ is finite.  
By construction, the fiber $X'_{\eta} = X' \times_X \Spec k(X)$
is finite \'etale over $k(X)$ and it is identified with $\coprod \Spec k(X'_i)$,
where $\{X'_i\}$ is the finite set of irreducible components of $X'$.
(This is called the {\em scheme of generic points} of $X'$.) 
The following lemma is a straightforward generalization (via \'etale descent)
of its well-known analogue for schemes. 

\begin{lemma}\label{normg}
Let $X$ be an irreducible and reduced locally noetherian algebraic space
that is locally of finite type over an excellent scheme $($resp.\ is normal$)$, 
and let $\eta$ denote its unique generic point.
Let $\mathscr{N}_X$ denote the category of finite  $($resp.\ finite
generically \'etale$)$ maps 
$f:X' \rightarrow X$ from normal algebraic spaces $X'$
such that ${\ms O}_X \rightarrow f_{\ast}{\ms O}_{X'}$ is injective.

The functor $X' \rightsquigarrow X'_{\eta}$ is an equivalence from
the category $\mathscr{N}_X$ to the category of non-empty finite
$($resp.\ non-empty finite generically \'etale$)$ reduced $k(X)$-schemes,
and normalization of $X$ in nonzero finite reduced $($resp.\ nonzero 
finite \'etale$)$ $k(X)$-algebras is a quasi-inverse.
\end{lemma}

Now we can give the proof of Proposition \ref{prop1}.   By Proposition \ref{irredcomp},
$X$ is irreducible.  Let $\eta$ denote its unique generic point.
Choose an \'etale covering $U \rightarrow X$ by an affine scheme, and 
let $L/k(X)$ be a finite Galois extension which splits the finite \'etale
$k(X)$-scheme $U_{\eta}$.  Let $\pi:X' \rightarrow X$ 
denote the normalization of $X$ in $L$.  Let $G = {\rm{Gal}}(L/k(X))$, so
by the equivalence in
Lemma \ref{normg} there is a natural right action by $G$ on $X'$ over $X$.
In particular, $G$ acts on the coherent ${\ms O}_X$-algebra
$\pi_{\ast}({\ms O}_{X'})$, so there is a natural injective map
${\ms O}_X \rightarrow \pi_{\ast}({\ms O}_{X'})^G$ of coherent ${\ms O}_X$-algebras.  
We claim that this is an isomorphism.  By normality it suffices to work over
a Zariski-dense open subspace of $X$, so taking such a subspace that is
an affine scheme does the job.  Since $L/k(X)$ is Galois, we likewise
see by working over such a dense open subscheme that $\pi$ is an
\'etale $G$-torsor over a dense open subspace of $X$. 

Since $X' \rightarrow X$ is finite, by Example \ref{exqt} the quotient
$X'/G$ exists and  
the natural map $X'/G \rightarrow X$ is finite.  We can say more:

\begin{lemma} The natural map $X'/G \rightarrow X$ is an isomorphism.
\end{lemma}

\begin{proof}  The finite map $X' \rightarrow X$ between irreducible noetherian algebraic spaces
is dominant, so the same holds for 
$X'/G \rightarrow X$.  The algebraic space $X'/G$ is also reduced
since $X'$ is reduced.  
The function field of $X'/G$ contains $k(X)$ and is contained in
$k(X')^G = L^G = k(X)$, so 
the finite map $X'/G \rightarrow X$ is birational.
It remains to use the fact that a finite birational
map between reduced noetherian algebraic spaces
is an isomorphism when the target is normal (as we may check by working \'etale-locally
to reduce to the known case of schemes).
\end{proof}

We have not yet used the precise way in which $L/k(X)$ was defined.   This is essential
to prove the next lemma, which will complete the proof of Proposition \ref{prop1}.

\begin{lemma} The algebraic space $X'$ is a scheme.
\end{lemma}

This assertion is \cite[Prop.~1]{raoult}, where a proof resting on
Zariski's Main Theorem is given.  
% (The statement of
% \cite[Prop.~1]{raoult} requires a noetherian hypothesis to ensure
% finiteness for integral closures in the proof.)  
We now give an
alternative proof below for the convenience of the reader.

\begin{proof}
Recall that $L/k(X)$ is a finite Galois splitting field of
the scheme of generic points $\Spec A_U$ of an affine \'etale scheme
cover $U \rightarrow X$.   Consider the algebraic space $P = X' \times_X U$.
This is finite over $U$, so it is an affine scheme, and it is clearly a quasi-compact \'etale cover
of $X'$.   In particular, $P$ is normal.  Each connected component $P_i$ 
of $P$ maps birationally to $X'$ since the scheme of generic points of $P$ is
$\Spec(L \otimes_{k(X)} A_U) \simeq \coprod \Spec L$
due to $L/k(X)$ being a Galois splitting field for each factor field of $A_U$. 
We shall prove that each $P_i$ maps to $X'$ via a monomorphism.
Any \'etale monomorphism of algebraic spaces is an open immersion
(as we deduce from the scheme case via descent), so it would follow
that the \'etale covering map $P \rightarrow X'$ realizes the $P_i$'s as a collection
of open subspaces that cover $X'$ and are schemes, whence $X'$ is a scheme as desired.

Now we show that each map $P_i \rightarrow X'$ is a monomorphism, or in other words
that the diagonal map $P_i \rightarrow P_i \times_{X'} P_i$ is an isomorphism.  
This diagonal is a closed immersion since $P_i \rightarrow X'$ is separated
(as the affine scheme $P_i$ is separated over $\Spec \Z$) and it is also \'etale,
so it is an open immersion too.   In other words, this diagonal realizes
$P_i$ as a connected component of $P_i \times_{X'} P_i$.  But this
fiber product has scheme of generic points $\Spec(L \otimes_L L) = \Spec(L)$
since $P_i \rightarrow X'$ is \'etale, so $P_i \times_{X'} P_i$ is irreducible.
Therefore $\Delta_{P_i/X'}$ is an isomorphism, as desired.
\end{proof}

\subsection{Proof of Theorem \ref{mainresult} when $X$ is normal}\label{normalproof}
 
 The aim of this section is to use 
 the known Nagata compactification theorem for schemes 
 (together with Proposition \ref{prop1}) to prove
 the following special case of Theorem \ref{mainresult}, from which we will deduce the general case. 
  
 \begin{theorem}\label{thm4}
 Let $f:X \rightarrow S$ be a separated map of finite type between algebraic
 spaces, with $S$ of finite presentation over an excellent noetherian
 scheme and $X$ normal.
 Then $X$ admits an $S$-compactification;
 i.e., $f$ factors through an open immersion $j:X \rightarrow \overline{X}$ into
 an $S$-proper algebraic space.
 \end{theorem}
 
 \begin{proof}
 {\bf Step 1}. We first reduce to the case when 
 $S$ is normal and both $X$ and $S$ are irreducible.  The main subtlety is
 that irreducibility is not \'etale-local.  We shall
 use the irreducible component decomposition 
 of noetherian algebraic spaces; see Proposition \ref{irredcomp}.  
 We may replace $S$ with the schematic image of the separated finite type
map $f:X \rightarrow S$, so ${\ms O}_S
\rightarrow f_{\ast}{\ms O}_X$ is injective. Thus, $S$ is reduced and
for each irreducible component $S_j$ of $S$, there is an irreducible
component $X_{i(j)}$ of $X$ that is carried onto a dense subset of 
$S_j$ by $f$.  In particular, the generic point of $X_{i(j)}$ is carried to the generic point
of $S_j$.  Writing $R_X$ and $R_S$ to denote the coordinate rings
of the schemes of generic points, the preceding says exactly that
$R_S \rightarrow R_X$ is a faithfully flat
ring extension.  This latter formulation has the advantage
that (unlike irreducible components) it is compatible with passing to quasi-compact \'etale covers
of $X$ and $S$.  

Let 
the finite map $\widetilde{S} \rightarrow S$ denote the normalization of
$S$ in its scheme of generic points (see Lemma \ref{normg} 
and the discussion preceding it).  We claim that $f$ uniquely factors through a 
(necessarily separated, finite type, and schematically dominant) map
$\widetilde{f}:X \rightarrow \widetilde{S}$.  This is well-known in the scheme case,
and to handle the general case we use \'etale descent for
morphisms:  by the claimed uniqueness we may work \'etale-locally
on $S$ to reduce to the case when it is a scheme, and we can then work
over an \'etale scheme cover of $X$ to reduce to the case when $X$ is also a scheme.  

By using $\widetilde{f}$, we may replace $S$
with $\widetilde{S}$ to reduce to the case
when $S$ is normal.  We may pass to connected components so
that $X$ and $S$ are both connected and hence are irreducible.  

{\bf Step 2}.  Now we make a digression, and 
prove Theorem \ref{mainresult} whenever $S$ is an excellent
noetherian scheme and $X$ is arbitrary (e.g., not necessarily normal). 
This generality will be useful in Step 3. 
By Proposition \ref{P:normalred} it suffices to treat the case when $X$ is normal. 
By Step 1 we may assume that both $X$ and $S$ are normal and
connected.  
% This case is asserted in \cite[Prop.~2]{raoult}, granting
% general facts about quotients of algebraic spaces by finite groups.
For the convenience of the reader, we explain the argument in terms of
the theory of quotients that we reviewed in \S\ref{gabberlemma}
(resting on the use of quotient stacks).

By Proposition \ref{prop1}, there is a normal noetherian scheme $X'$ equipped
with a right action by a finite group $G$ and a $G$-invariant finite map
$X' \rightarrow X$ inducing an isomorphism
$X'/G \simeq X$.   By Nagata's compactification theorem for schemes,
there is an $S$-compactification $j:X' \hookrightarrow \overline{X}'$
with $\overline{X}'$ a proper $S$-scheme.  For each $g \in G$, let
$j^{(g)} = j \circ [g]$ where $[g]:X' \simeq X'$ is the $S$-automorphism given by the right 
action by $g \in G$ on $X'$.    Thus, the fiber product
$P = \prod_{g \in G} \overline{X}'$ over $S$ is a proper $S$-scheme
admitting a right $G$-action via $[g_0]:(\overline{x}'_g)_g \mapsto
(\overline{x}'_{g_0g})_g$ for varying $g_0 \in G$.
The map $X' \rightarrow P$ defined by $x' \mapsto (j^{(g)}(x'))_g$
is an immersion, and it is $G$-equivariant since
$j^{(g_0g)}(x') = j([g_0g](x')) = j^{(g)}([g_0](x'))$. 
Hence, the schematic closure $\overline{X}'_1$ of $X'$ in $P$
is a $G$-equivariant compactification of $X'$ over $S$.

Passing to quotient stacks, $[X'/G] \rightarrow [\overline{X}'_1/G]$
is an open immersion over $S$.  Passing to the coarse moduli spaces,
we get an $S$-map $X \simeq X'/G \rightarrow \overline{X}'_1/G$ 
with $\overline{X}'_1/G$ proper over $S$.  This is also an open immersion because
the formation of coarse moduli spaces is compatible with passage to open substacks
(as a special case of compatibility with quasi-separated flat base change on
the coarse space).
Hence, we have constructed an $S$-compactification of $X$
as an algebraic space.  
This proves Theorem \ref{mainresult} whenever $S$ is an excellent noetherian scheme.

{\bf Step 3}.  Finally, we return to the situation at the end of Step 1,
so $S$ is a connected normal algebraic space
of finite presentation over an excellent noetherian scheme and $X$ is
normal and connected, hence generically flat over $S$.  
We seek to construct an $S$-compactification of $X$. 
The quotient technique that was used in Step 2 will be applied in order to reduce to
the settled case when $S$ is a scheme (and $X$ is arbitrary).  

Since $S$ is normal, by
Proposition \ref{prop1} we have $S \simeq S'/G$ for some normal
noetherian scheme $S'$ equipped with a right action by a finite group
$G$ and a finite surjective $G$-invariant map $S' \rightarrow S$.
Thus, $S'$ is an excellent noetherian scheme.  Let $X' = (X \times_S
S')_{\text{red}}$, so $X'$ has a natural $G$-action over $X$ and $X'
\rightarrow X$ is a finite surjective $G$-invariant map.  Since $X$ is
generically flat over $S$, the induced finite map $X'/G\rightarrow X$ is an isomorphism between dense opens,
 and thus is an isomorphism because $X$ is normal.
Beware that $X'$ may not be normal.

Since $S'$ is an excellent noetherian scheme and
in Step 2 we proved Theorem \ref{mainresult} 
whenever the base
is an excellent noetherian scheme, there is an $S'$-compactification
$j:X' \hookrightarrow \overline{X}'$.  
For each $g \in G$, let $[g]_{S'}:S' \simeq S'$ and
$[g]_{X'}:X' \simeq X'$ denote the action maps for $g$ on $S'$ and $X'$
respectively (so $[g]_{X'}$ is a map over $[g]_{S'}$).   Let 
${\overline{X}'}^{(g)} = S' \times_{[g]_{S'}, S'} \overline{X}'$. 
Since $j \circ [g]_{X'}:X' \rightarrow \overline{X}'$ is an open immersion over
the automorphism $[g]_{S'}$ of $S'$, it induces an open immersion
$j^{(g)}:X' \hookrightarrow {\overline{X}'}^{(g)}$ over $S'$.

For the fiber product $P = \prod_{g \in G} {\overline{X}'}^{(g)}$ over $S'$,  
any $g_0 \in G$ induces an isomorphism ${\overline{X}'}^{(g)} \rightarrow
{\overline{X}'}^{(g_0^{-1}g)}$ over $[g_0]_{S'}$.  For a fixed $g_0 \in G$
these isomorphisms between
the factors, all over a common automorphism of $S'$, combine to define
an automorphism of $P$ over the automorphism $[g_0]_{S'}$ of $S'$,
and this is a right $G$-action on $P$ over the right $G$-action on $S'$.  Moreover,
the immersion $X' \rightarrow P$ defined by $x' \mapsto (j^{(g)}(x))$
is $G$-equivariant, exactly as in Step 2.   Hence, the schematic
closure $Z$ of $X'$ in $P$ is an $S'$-proper algebraic
space equipped with a right $G$-action over the one on $S'$,
and this action is compatible with the given $G$-action on the open subscheme $X'$.
The induced $S$-map $[X'/G] \rightarrow [Z/G]$ of Deligne--Mumford stacks
is therefore an open immersion, 
so (exactly as in Step 2) it induces an open immersion of algebraic spaces $X \simeq X'/G \hookrightarrow
Z/G$ over $S$.  
Since $Z$ is $S'$-proper, so is $Z/G \rightarrow S$. 
\end{proof}

\subsection{Conclusion of the proof of Theorem \ref{mainresult}}
\label{sec:concl}

By Proposition \ref{P:normalred} it is enough to prove Theorem
\ref{mainresult} when $X$ is normal, and this is covered by Theorem \ref{thm4}.

\section{Approximation results}\label{approxsec}

This section is devoted to establishing several general technical results which allow
us to reduce problems to the noetherian (quasi-separated) 
case, and even to the case of algebraic spaces
of finite presentation over $\Z$.   In particular, we will prove
Theorem \ref{mainresult} by using the settled case from
\S\ref{exccase} with $S$ of finite presentation over $\Z$. 

\subsection{Absolute noetherian approximation}\label{absnoetherian}

The key to avoiding noetherian hypotheses in Theorem \ref{mainresult} is the absolute noetherian 
approximation result in Theorem \ref{absapprox}.  We will prove Theorem \ref{absapprox} by
reducing it to the known case when $S$ is a scheme \cite[Thms.~C.7,~C.9]{tt}, in which case 
all $S_{\lambda}$ can be taken to be schemes. 
The reduction to the scheme case rests on the
fact that any qcqs algebraic space admits a special kind of finite stratification by schemes: 

\begin{theorem}[Raynaud--Gruson]\label{rgstrat}
Let $S$ be a qcqs algebraic space.  There is a finite rising chain
$$\emptyset = U_0 \subseteq U_1 \subseteq \dots \subseteq U_r = S$$
of quasi-compact open subspaces such that
for each $i > 0$ the open subspace $U_i$ admits an 
\'etale cover $\varphi_i:Y_i \rightarrow U_i$ by a quasi-compact separated scheme $Y_i$ 
with $\varphi_i$ restricting to an isomorphism over the closed subspace
$Z_i = U_i - U_{i-1}$ in $U_i$ endowed with its reduced structure.  Moreover, 
each $\varphi_i$ is separated, and 
each $Z_i$ is a separated and quasi-compact scheme. 
\end{theorem}

\begin{proof}  The statement of \cite[I,~5.7.6]{rg} gives the existence result, except with 
\'etale covers $\varphi'_i:Y'_i \rightarrow (U_i)_{\rm{red}}$ by quasi-compact
separated schemes $Y'_i$ such that each $\varphi'_i$ restricts to an isomorphism over $Z_i$.  
The construction of
$\varphi'_i$ is as a pullback of an \'etale cover
$\varphi_i:Y_i \rightarrow U_i$ by a quasi-compact
separated scheme $Y_i$.   The map 
 $\varphi_i$ is necessarily separated since the composition of 
$\Delta_{\varphi_i}$ with the monomorphism $Y_i \times_{U_i} Y_i
\rightarrow Y_i \times_{\Spec \Z} Y_i$ is a closed immersion.
%The map $\varphi_1:Y_1 \rightarrow U_1$ is a quasi-compact and separated
%\'etale map of algebraic spaces which restricts to an isomorphism over $Z_1 = (U_1)_{\rm{red}}$.
%A quasi-compact and separated \'etale map
%of algebraic spaces is representable in schemes \cite[II,~Cor.~6.17]{knutson}  and 
%so is an isomorphism if it has constant fiber-rank equal to 1
%\cite[Lemma~A.1.4]{conrad2}.   Hence, $\varphi_1$ is an isomorphism, so 
%$U_1$ is a scheme. 
\end{proof}

\begin{definition}\label{defapprox}  Let $\Lambda$ be a noetherian ring and 
$S$ an algebraic space over $\Lambda$.  We say that $S$ is {\em $\Lambda$-approximable}
if there is a $\Lambda$-isomorphism
$S \simeq \invlim S_{\alpha}$ where $\{S_{\alpha}\}$ is an inverse system
of algebraic spaces of finite presentation over $\Lambda$ having affine transition maps
$S_{\beta} \rightarrow S_{\alpha}$ for all $\alpha$ and all $\beta \ge \alpha$.
In case $\Lambda = \Z$, we say that $S$ is {\em approximable}. 
\end{definition}

Observe that we use ``finite presentation''
rather than just ``finite type'' in Definition \ref{defapprox}.  This is essential,
as we indicated in \S\ref{notation}.  
Any inverse limit as
in Definition \ref{defapprox} is necessarily qcqs
(over $\Lambda$, or equivalently over $\Z$), and our aim is to prove that
conversely every qcqs algebraic space over $\Lambda$
is $\Lambda$-approximable.  The most interesting case is $\Lambda = \Z$, 
and in fact this is enough to settle the general case:

\begin{lemma}\label{zapprox}
 Let $\Lambda$ be a noetherian ring, and $S$ a $\Lambda$-approximable
algebraic space.  The inverse system $\{S_{\alpha}\}$ as in Definition $\ref{defapprox}$
can be taken to have schematically dominant affine transition maps.  Moreover,
if $\Lambda \rightarrow \Lambda'$ is a map of noetherian rings
and $S$ admits a compatible structure of algebraic space over $\Lambda'$ then
$S$ is also $\Lambda'$-approximable.
\end{lemma}

\begin{proof}
Choose an inverse system $\{S_{\alpha}\}$ of algebraic
spaces of finite presentation over $\Lambda$
with affine transition maps such that $S \simeq \invlim S_{\alpha}$
over $\Lambda$.  Each map $q_{\alpha}:S \rightarrow S_{\alpha}$
is affine, so it admits a scheme-theoretic image $S'_{\alpha} \subseteq S_{\alpha}$
that is the closed subspace corresponding to the quasi-coherent kernel of
${\ms O}_{S_{\alpha}} \rightarrow q_{\alpha\ast}({\ms O}_S)$.  By working \'etale-locally over 
a fixed $S_{\alpha_0}$ we see that the map $q':S \rightarrow \invlim S'_{\alpha}$ is
an isomorphism and 
$q'_{\alpha}:S \rightarrow S'_{\alpha}$ and $q'_{\alpha \beta}: S'_{\beta} \rightarrow S'_{\alpha}$
are schematically dominant and affine for all $\alpha$ and all $\beta \ge \alpha$.  

Now assume there is given a 
(necessarily quasi-separated) $\Lambda$-morphism $S \rightarrow \Spec \Lambda'$ for 
a noetherian $\Lambda$-algebra
$\Lambda'$.  Fix $\alpha_0$ and define the quasi-coherent
sheaf
$$\calA_{{\alpha}} := \Lambda' \cdot q'_{\alpha_0,\alpha \ast}({\ms O}_{S'_{\alpha}})
\subseteq q'_{\alpha_0 \ast}({\ms O}_S)$$
of ${\ms O}_{S_{\alpha_0}}$-algebras 
for $\alpha \ge \alpha_0$.
The algebraic spaces $S''_{\alpha} = \Spec_{S'_{\alpha_0}}(\calA_{\alpha})$ 
of finite presentation over $\Lambda'$ form an inverse system with schematically
dominant and affine transition maps such that $\invlim S''_{\alpha} \simeq S$
over $\Lambda'$. 
\end{proof}

To make effective use of Theorem \ref{rgstrat} in the proof that $S$ is approximable, 
the idea is to use induction on $r$ to get to the situation of extending an inverse system approximation
for $U_{r-1}$ across the complementary {\em scheme} $Z_r$ to get an inverse system approximation 
for $S$.  Making such an extension requires
studying the algebraic spaces $U_{r-1}$ and $S$ in terms of 
\'etale scheme presentations, which must themselves be replaced with compatible inverse limit approximations.
The strategy is to reconstruct $S$ as a kind of pushout of $Z_r$ against an \'etale
scheme chart of $U_{r-1}$, but everything must be done at the level of inverse systems.
Moreover, the constructions we give have to be sufficiently ``cartesian'' so that we can compute
what happens in the limit.  This leads us to the following initial lemma in which a cartesian property
at the level of \'etale scheme charts implies a cartesian property for the quotients.
It will be more convenient to work with sheaf-functors instead of schemes:  

\begin{lemma}\label{lemma3} Let $R \rightrightarrows U$ and $R' \rightrightarrows U'$
be equivalence relations in sheaves of sets
on the \'etale site of the category of schemes.  Assume that there is
given a map $f:U' \rightarrow U$ such that $f \times f:U' \times U' \rightarrow U \times U$ carries
$R'$ into $R$ and the co-commutative diagram 
$$\xymatrix{ {R'} \ar[r] \ar@<-2pt>[d]_-{p'_1} \ar@<2pt>[d]^-{p'_2} & {R} \ar@<-2pt>[d]_-{p_1} 
\ar@<2pt>[d]^-{p_2} \\ {U'} \ar[r]_-{f} & {U}}$$
is cartesian for each pair $(p_i, p'_i)$.  Then the induced commutative square
$$\xymatrix{ {U'} \ar[r] ^-{f} \ar[d] & {U} \ar[d] \\ {U'/R'} \ar[r] & {U/R}}$$
is cartesian.
\end{lemma}

\begin{proof}
We have to show that the natural map $U' \rightarrow (U'/R') \times_{U/R} U$
is an isomorphism as \'etale sheaves on any scheme $T$.  It suffices to check on 
stalks, which is to say on $T$-valued points at geometric points
$\overline{t}:\Spec k \rightarrow T$, with $k$ a separable closure of the residue field
at an arbitrary point $t \in T$.  
Hence, $(U/R)_{\overline{t}} = 
U_{\overline{t}}/R_{\overline{t}}$ and $(U'/R')_{\overline{t}} = U'_{\overline{t}}/R'_{\overline{t}}$, so we need to prove that the natural map
\begin{equation}\label{xxr}
U'_{\overline{t}} \rightarrow (U'_{\overline{t}}/R'_{\overline{t}}) \times_{U_{\overline{t}}/R_{\overline{t}}} U_{\overline{t}}
\end{equation}
is bijective.  We will do this using the commutative diagram
\begin{equation}\label{cartr}
\xymatrix{
{R'_{\overline{t}}} \ar[d]_-{g=f\times f} \ar[r]^-{p'_i} & {U'_{\overline{t}}} \ar[r]^-{\pi'} \ar[d]^-{f} & {U'_{\overline{t}}/R'_{\overline{t}}} \ar[d] \\
{R_{\overline{t}}} \ar[r]_-{p_i} & {U_{\overline{t}}} \ar[r] & {U_{\overline{t}}/R_{\overline{t}}}}
\end{equation}
with cartesian left square for $i \in \{1, 2\}$.

To prove surjectivity of (\ref{xxr}), choose $\alpha \in (U'_{\overline{t}}/R'_{\overline{t}}) \times_{U_{\overline{t}}/R_{\overline{t}}} U_{\overline{t}}$, so
$\alpha = (u' \bmod R'_{\overline{t}}, u)$ for some $u' \in U'_{\overline{t}}$ and $u \in U_{\overline{t}}$.  
The fiber product condition on $\alpha$ 
says $(f(u'), u) \in R_{\overline{t}}$ inside $U_{\overline{t}} \times U_{\overline{t}}$.  The cartesian property of the left square
in (\ref{cartr})
with $i = 1$ therefore gives a unique point $r' \in R'_{\overline{t}}$ with $g(r') = (f(u'),u)$ and $p'_1(r') = u'$.
The commutativity of the left square with $i = 2$ says that if $y' = p'_2(r')$ then
$f(y') = u$.  Hence, (\ref{xxr}) carries $y' \in U'_{\overline{t}}$ over to 
$(y' \bmod R'_{\overline{t}}, u) = (u' \bmod R'_{\overline{t}}, u) = \alpha$.

Now pick $u'_1, u'_2 \in U'_{\overline{t}}$ that are carried to the same point under (\ref{xxr}), which is to say that 
$(u'_1, u'_2) \in R'_{\overline{t}}$ and $f(u'_1) = f(u'_2)$ in $U_{\overline{t}}$.  Letting $r' = (u'_1, u'_2)$, clearly 
$g(r') = (f(u'_1), f(u'_2)) = (f(u'_1), f(u'_1))$.  That is, $g(r') = \Delta(f(u'_1))$
where $\Delta:U \rightarrow R$ is the diagonal section.   But the point $\Delta(u'_1) \in
R'_{\overline{t}}$ satisfies
$$g(\Delta(u'_1)) = (f(u'_1), f(u'_1)) = g(r'),\,\,\,
p'_1(\Delta(u'_1)) = u'_1 = p'_1(r'),$$
so the cartesian property of the left square in (\ref{cartr}) for $i = 1$ implies that
$r' = \Delta(u'_1) = (u'_1, u'_1)$.  Since $r' = (u'_1, u'_2)$ by definition, we get $u'_2 = u'_1$ as required
for injectivity of (\ref{xxr}).
\end{proof}

\begin{corollary}\label{cor4}
In the setup of Lemma $\ref{lemma3}$, if $U$, $U'$, $R$, and $R'$ are algebraic spaces
with maps $p_i$ and $p'_i$ \'etale, and if 
$U' \rightarrow U$ satisfies a property
$\P$ of morphisms of algebraic spaces that is \'etale-local on the base, then 
$U'/R' \rightarrow U/R$ satisfies property $\P$. 
\end{corollary}

By Corollary \ref{qtspace}, $U/R$ and $U'/R'$ are algebraic spaces. 

\begin{proof}
To analyze the asserted property of $U'/R' \rightarrow U/R$ it suffices to check after
pullback to the \'etale covering $U$ of $U/R$.  By Lemma \ref{lemma3},
this pullback is identified with the map $U' \rightarrow U$.
\end{proof}

To apply Corollary \ref{cor4}, we wish to 
describe a situation in which the setup of Lemma \ref{lemma3}
naturally arises.  We first require one further lemma, concerning
the existence and properties of certain pushouts.  

\begin{lemma}\label{pushoutlemma} Consider a diagram of algebraic spaces
$$\xymatrix{
U' \ar[r]^-{j'} \ar[d]_-{p'} & X' \\ U}$$
in which $j'$ is an open immersion and $p'$ is an \'etale surjection.
\begin{enumerate}
\item There exists a pushout $X = U \coprod_{U'} X'$ in the category
  of algebraic spaces, and the associated diagram
\begin{equation}\label{uuxx}
\xymatrix{
  U' \ar[r]^-{j'} \ar[d]_-{p'} & X' \ar[d]^-{p} \\
  U \ar[r]_-{j} & X}
  \end{equation} is cartesian, with $j$ an open immersion and
$p$ an \'etale surjection.   The formation of this pushout commutes with any base change on $X$
in the sense that if $X_1 \rightarrow X$ is any map of algebraic spaces then the pullback diagram
\begin{equation}\label{uuxx1} 
\xymatrix{
  U'_1 \ar[r]^-{j'_1} \ar[d]_-{p'_1} & X'_1 \ar[d]^-{p_1} \\
  U_1 \ar[r]_-{j_1} & X_1}
  \end{equation}
  is also a pushout. 
\item If $j'$ is quasi-compact then $j$ is quasi-compact, if $p'$ is
  separated then $p$ is separated, and if $j'$ is quasi-compact and
  $p'$ is finitely presented then $p$ is also finitely presented.
\item If $U$, $U'$, and $X'$ are qcqs then so is $X$.
\end{enumerate}
\end{lemma}

Note that in (1), the asserted compatibility with any base change is only being made for the kinds of pushouts
considered there, namely pushout of an \'etale surjection along an open immersion.  This is a mild (but useful)
generalization of the familiar base change compatibility of a Zariski gluing, which can be viewed as a pushout of
one open immersion along another. 

\begin{proof}
Since $\Delta_{p'}:U' \rightarrow U' \times_U U'$ is a section to an \'etale map
of algebraic spaces, it is an \'etale map. 
Thus,  $\Delta_{p'}$ is an \'etale monomorphism.
An \'etale monomorphism of algebraic spaces is always an open immersion.
Indeed, in the special case of quasi-separated algebraic spaces we may work
Zariski-locally to reduce to the finitely presented case, which is \cite[II,~Lemma~6.15b]{knutson}.
This handles the general case when the target is affine, as then the source is separated
(due to separatedness of monomorphisms).  In 
 general, we may work \'etale-locally on the base  to reduce to the settled case
when the target is affine.

Since $\Delta_{p_1}$ is an open immersion, 
it makes sense to form the gluing $$R = X' \coprod_{U'} (U' \times_U U')$$
of $X'$ and $U' \times_U U'$ along the common open subspace $U'$.
Using the maps $U' \times_U U' \rightrightarrows U' \stackrel{j'}{\hookrightarrow} X'$, we arrive at 
natural maps $R \rightrightarrows X'$ that clearly constitute an \'etale equivalence
relation in algebraic spaces.  The quotient $X := X'/R$ is an
algebraic space, and by construction $p:X' \rightarrow X$ is an \'etale surjection.
By the definition of $R$, the map $j':U' \hookrightarrow X'$ induces a canonical map
$$j:U = U'/(U' \times_U U') \rightarrow X'/R  = X,$$
so we obtain the commutative diagram (\ref{uuxx}).  By using the definition of $R$, this diagram 
is easily checked to be 
cartesian and a pushout.   Due to the cartesian property and the fact that
$p$ is an \'etale surjection, $j$ is an open immersion because $j'$ is an open immersion.

To complete the proof of (1), we have to verify the compatibility with base
change on $X$.   Letting $X_1 \rightarrow X$ be a map from an algebraic space, we have
to prove that (\ref{uuxx1}) is a pushout.  That is, we want the natural $X_1$-map
$U_1 \coprod_{U'_1} X'_1 \rightarrow X_1$ to be an isomorphism.  
In view of the construction of this latter pushout as a quotient, it is equivalent to say that 
the map $p_1:X'_1 \rightarrow X_1$ is an \'etale cover and 
$$h_1:R_1 = X'_1 \coprod_{U'_1} (U'_1 \times_{U_1} U'_1) \rightarrow X'_1 \times_{X_1} X'_1$$
is an isomorphism, where $R_1$ is a gluing for the Zariski topology.   

The map $p_1$ is a base change of $p$, so it is an \'etale surjection.   To prove that $h_1$ is an isomorphism,
we observe that $X'_1 \times_{X_1} X'_1 = (X' \times_X X') \times_X X_1$ and likewise
$R_1 = R \times_X {X_1}$ because the formation of Zariski gluings of objects over a common base
(such as $X'$ and $U' \times_U U'$ glued along the open $U'$, all over the base $X$) commutes
with any base change.   Hence, $h_1$ is identified with the base change along $X_1 \rightarrow X$ of
the analogous map $h:R \rightarrow X' \times_X X'$ that is an isomorphism (because $X := X'/R$).  This
completes the proof of (1). 

Consider the claims in (2).   By descent through $p$, the map $j$ is
 quasi-compact when $j'$ is quasi-compact.
Assuming that $p'$ is separated, the open subspace
$U'$ in $U' \times_U U'$ via the diagonal is also closed
and hence splits off as a disjoint union:  $U' \times_U U' = \Delta(U') \coprod V$
for an algebraic space $V$ that is separated over $U'$ (and hence
over $X'$) via either projection.  Thus, in such cases $R = X' \coprod V$ is separated
over $X'$ via either projection, so $p:X' \rightarrow X$ is separated.  
In case $j'$ is quasi-compact (so $U' \rightarrow X'$ is finitely presented)
and $p'$ is finitely presented, the map
$p:X' \rightarrow X$ is finitely presented because it is a descent of either of the projection maps 
$R \rightrightarrows X'$ which express $R$ as a gluing of two finitely presented
$X'$-spaces along a common finitely presented open subspace. 

Finally, to prove (3), observe that (by construction) if
$U$, $U'$, and $X'$ are qcqs then $R$ is qcqs, so the maps $R \rightrightarrows X'$
are qcqs and hence the quotient map $X' \rightarrow X$ is qcqs.  Thus, in such cases
$X$ is qcqs.
\end{proof}

\begin{notation}
  We will sometimes refer to the \'etale equivalence relation $X'
  \times_X X' \rightrightarrows X'$ constructed in Lemma
  \ref{pushoutlemma} as being obtained from the \'etale equivalence
  relation $U' \times_{U} U' \rightrightarrows U$ via \emph{extension along
  the diagonal}.
\end{notation}

To prove Theorem \ref{mainresult}, we wish to inductively construct limit presentations of
qcqs algebraic spaces by means of stratifications as in Theorem \ref{rgstrat}.  This will
be achieved by using the following result.  

\begin{proposition}\label{prop8}
Let $X$ be a qcqs algebraic space, and suppose there is given a diagram
\begin{equation}\label{uxz}
\xymatrix{{U} \ar[r]^-{j} \ar[d]_-{\pi} & {Z} \\ {X} & }
\end{equation} 
in which $\pi$ is a finitely presented \'etale scheme covering and $j$ is 
an open immersion into a qcqs scheme $Z$.
Form the cartesian pushout diagram in algebraic spaces
\begin{equation}\label{ypush}
\xymatrix{{U} \ar[r]^-{j} \ar[d]_-{\pi} & {Z} \ar[d] \\ {X} \ar[r] & {Y}}
\end{equation}
as in Lemma $\ref{pushoutlemma}$, so
the bottom side is an open immersion, the right side is
an \'etale surjection, and $Y$ is qcqs.

If $X$ is approximable then so is $Y$.
\end{proposition}

To prove Proposition \ref{prop8}, we first need to study pairs of diagrams of the type in (\ref{uxz})
that are connected to each other via affine and schematically dominant maps.  Thus, we now briefly
digress to consider such diagrams and their corresponding pushouts as in Lemma \ref{pushoutlemma}.

Let $X' \rightarrow X$ be an affine and schematically dominant map of
algebraic spaces, and let $U \rightarrow X$ be an \'etale covering by a scheme, so
$h:U' := U \times_X X' \rightarrow U$ is affine (hence $U'$ is a scheme)
and $U' \rightarrow X'$ is an \'etale covering.  Note that the affine map $U' \rightarrow U$ is 
schematically dominant.  Suppose that there is a cartesian square of schemes
$$\xymatrix{{U'} \ar[r]^-{j'} \ar[d]_-{h} & {Z'} \ar[d] \\ {U} \ar[r]_{j} & {Z}}$$
in which the horizontal maps are open immersions and 
the right vertical map is affine (like the left side).  

The respective algebraic space quotients $X$ and $X'$ of $U$ and $U'$
give rise to \'etale equivalence relations in schemes
$$U \times_X U \rightrightarrows U,\,\,\,
U' \times_{X'} U' \rightrightarrows U',$$
and we extend these to \'etale equivalence relations in schemes
$$R \rightrightarrows Z,\,\,\,R' \rightrightarrows Z'$$
via extension along the diagonal, exactly as in the proof
of Lemma \ref{pushoutlemma}:  define the subfunctors $R \subseteq Z \times Z$ and
$R' \subseteq Z' \times Z'$ to respectively be the gluings along common open subspaces
\begin{equation}\label{rdz}
R = \Delta(Z) \coprod_{\Delta(U)} (U \times_X U),\,\,\,
R' = \Delta(Z') \coprod_{\Delta'(U')} (U' \times_{X'} U').
\end{equation} 
In particular, the diagrams
$$\xymatrix{{U} \ar[r] \ar[d]_-{h} & {Z} \ar[d] \\ {X} \ar[r] & {Z/R}}\,\,\,\,\,\,\,\,
\xymatrix{{U'} \ar[r] \ar[d]_-{h} & {Z'} \ar[d] \\ {X'} \ar[r] & {Z'/R'}}$$
are cartesian and are pushouts with open immersions along the bottom, and if (as 
in applications below) $X$, $U$, and $Z$ are qcqs (so likewise for $X'$, $U'$, and $Z'$) then 
$R$, $R'$, $Z/R$, and $Z'/R'$ are qcqs.

\begin{corollary}\label{cor5}
In the above situation, the 
co-commutative
diagram 
\begin{equation}\label{rrz}
\xymatrix{ {R'} \ar[r] \ar@<2pt>[d]^-{p'_2} \ar@<-2pt>[d]_-{p'_1} & {R} \ar@<2pt>[d]^-{p_2} 
\ar@<-2pt>[d]_-{p_1} \\ {Z'} \ar[r] & {Z}}
\end{equation}
is cartesian for each pair $(p_i, p'_i)$, and the 
map of pushouts $Z'/R' \rightarrow Z/R$ is affine and schematically dominant.
\end{corollary}

\begin{proof}  It is straightforward to check that the cartesian
  property holds, using
that $U' = U \times_X X'$ by definition. 
Thus, the hypotheses of Lemma \ref{lemma3} are satisfied.  By Corollary \ref{cor4} we are done.
\end{proof}

The reason for our interest in Corollary \ref{cor5} is that it arises in the proof of 
Proposition \ref{prop8}, which we now give:

\begin{proof}[Proof of Prop.~$\ref{prop8}$] 
Since $X$ is approximable, we may choose an isomorphism
$X \simeq \invlim X_{\alpha}$ with $\{X_{\alpha}\}$ an inverse system of
algebraic spaces of finite presentation over $\Z$ with affine transition maps.
We may and do arrange that the transition maps are also schematically dominant.  
By Proposition \ref{lemma7}, we may also assume (by requiring $\alpha$ to be
sufficiently large) that this isomorphism is covered by an isomorphism
$U \simeq \invlim U_{\alpha}$ where $\{U_{\alpha}\}$ is an inverse system
of finitely presented {\em schemes} over $\{X_{\alpha}\}$ such that the maps 
$U_{\beta} \rightarrow X_{\beta} \times_{X_{\alpha}} U_{\alpha}$ are isomorphisms 
whenever $\beta \ge \alpha$.  (In particular, $\{U_{\alpha}\}$ has affine
transition maps, so $\invlim U_{\alpha}$ makes sense.)  

By Corollary \ref{mapdescend2} we may
and do require $\alpha$ to be sufficiently large so that 
the finitely presented maps $h_{\alpha}:U_{\alpha} \rightarrow X_{\alpha}$ are  \'etale coverings. 
Thus, by flatness of $h_{\alpha}$, the inverse system 
$\{U_{\alpha}\}$ has schematically dominant transition maps 
since the same holds for $\{X_{\alpha}\}$. 
Moreover, each scheme $U_{\alpha}$ is of finite type over $\Z$
since $X_{\alpha}$ is of finite presentation over $\Z$ (and $h_{\alpha}$ is finitely presented).  
Hence, by applying Lemma \ref{lemma2} to the quasi-compact open immersion
$U \hookrightarrow Z$, at the expense of possibly modifying the indexing system
we can arrange that there is a {\em cartesian} inverse system of quasi-compact 
open immersions $j_{\alpha}:U_{\alpha} \hookrightarrow Z_{\alpha}$ of
finite type $\Z$-schemes 
such that $\{Z_{\alpha}\}$ has affine and schematically
dominant transition maps and $\invlim j_{\alpha}$
is the given open immersion $j:U \hookrightarrow Z$.  
We emphasize that it is the application of Lemma \ref{lemma2} that is
the entire reason we had to make the affine transition maps in our initial inverse systems 
be schematically dominant. 

Consider the system of diagrams
$$\xymatrix{
{U_{\alpha}} \ar[r]^-{j_{\alpha}} \ar[d]_-{h_{\alpha}} & {Z_{\alpha}} \\
{X_{\alpha}} & }$$
in which the maps on the left are 
 \'etale scheme coverings. 
This is a ``cartesian'' system of diagrams in the sense that 
the diagrams 
\begin{equation}\label{lastpush}
\xymatrix{
{U_{\alpha'}} \ar[r]^-{j_{\alpha'}} \ar[d] & {Z_{\alpha'}} \ar[d] \\
{U_{\alpha}} \ar[r]_-{j_{\alpha}} & {Z_{\alpha}}},\,\,\,\,\,\,
\xymatrix{{U_{\alpha'}} \ar[d]_-{h_{\alpha'}} \ar[r] & {U_{\alpha}} \ar[d]^-{h_{\alpha}} \\ {X_{\alpha'}}
 \ar[r] & {X_{\alpha}}}
 \end{equation}
for $\alpha' \ge \alpha$ are cartesian.  Thus, 
the setup preceding Corollary \ref{cor5} is applicable
to the system of \'etale scheme coverings $h_{\alpha}:U_{\alpha} \rightarrow X_{\alpha}$
with affine and schematically dominant transition maps, equipped
with the compatible open immersions $j_{\alpha}:U_{\alpha} \hookrightarrow
Z_{\alpha}$.  

By Corollary \ref{cor5}, we thereby obtain a cartesian system of
\'etale equivalence relations $R_{\alpha} \rightrightarrows Z_{\alpha}$
in qcqs schemes akin to (\ref{rdz}) and 
the resulting qcqs algebraic space quotients $Y_{\alpha} = Z_{\alpha}/R_{\alpha}$
naturally fit into an inverse system with affine and schematically dominant transition maps. 
These quotients are exactly the pushouts $X_{\alpha} \coprod_{U_{\alpha}} Z_{\alpha}$ as
constructed in Lemma \ref{pushoutlemma}.  
Each $Y_{\alpha}$ is of finite type over $\Z$ since the same holds
for its \'etale scheme covering $Z_{\alpha}$, so each $Y_{\alpha}$ is finitely
presented over $\Z$ (as each $Y_{\alpha}$ is quasi-separated). 

Consider the pushout diagram
\begin{equation}\label{uzx}
\xymatrix{{U_{\alpha}} \ar[r]^-{j_{\alpha}} \ar[d] & {Z_{\alpha}} \ar[d] \\ 
{X_{\alpha}} \ar[r] & {Y_{\alpha}}}
\end{equation}
of algebraic spaces as in (\ref{uuxx}), 
so these are cartesian and have $X_{\alpha} \rightarrow Y_{\alpha}$
an open immersion.   For any $\alpha' \ge \alpha$,
the diagram (\ref{uzx}) for $\alpha'$ maps to the one for $\alpha$, and we claim that
this resulting inverse system of diagrams is a cartesian system in the sense that 
the pullback of (\ref{uzx})  along $Y_{\alpha'} \rightarrow Y_{\alpha}$
is identified (via the natural maps) with the $\alpha'$-version of (\ref{uzx}). 
This cartesian system claim along the right side of the diagrams (\ref{uzx}) is 
Lemma \ref{lemma3}, and along the top and left sides it is the cartesian property observed in (\ref{lastpush}). 
Since $X_{\alpha}$ is the image of the open subspace $U_{\alpha} \hookrightarrow Z_{\alpha}$
along the \'etale quotient map $Z_{\alpha} \rightarrow Y_{\alpha}$, and likewise
for $\alpha' \ge \alpha$, the cartesian property along the bottom sides of the diagrams (\ref{uzx}) follows.

The definition of $Y$ as a pushout
provides maps $Y \rightarrow Y_{\alpha}$
respecting change in $\alpha$, and we shall prove that
the induced map $Y \rightarrow \invlim Y_{\alpha}$ is an isomorphism.
This would show that $Y$ is approximable, as desired. 

Define $R = \invlim R_{\alpha} \subseteq 
\invlim (Z_{\alpha} \times_{\Spec \Z} Z_{\alpha}) = Z \times_{\Spec \Z} Z$, 
so $R$ is a qcqs scheme
and the pair of maps $p_1, p_2:R \rightrightarrows Z$ obtained from passage
to the limit on the cartesian system 
$p_{1,\alpha}, p_{2,\alpha}:R_{\alpha} \rightrightarrows Z_{\alpha}$
is an \'etale equivalence relation.  
Lemma \ref{lemma3} ensures that the natural maps
$Z_{\beta} \rightarrow Z_{\alpha} \times_{Y_{\alpha}} Y_{\beta}$
are isomorphisms for all $\beta \ge \alpha$, so passing to the limit on $\beta$ 
with a fixed $\alpha$ gives that the natural map
$Z \rightarrow Z_{\alpha} \times_{Y_{\alpha}} Y$ is an isomorphism
(since inverse limits of algebraic spaces under affine transition maps 
commute with fiber products).  Similarly, for each fixed $i \in \{1, 2\}$ 
and $\beta \ge \alpha$, the natural map 
$R_{\beta} \rightarrow R_{\alpha} \times_{p_{i,\alpha}, Z_{\alpha}} Z_{\beta}$
over $p_{i,\beta}:R_{\beta} \rightarrow Z_{\beta}$ is an isomorphism
due to the cartesian observation preceding Corollary \ref{cor5}.
Hence, passing to the limit on $\beta$ with a fixed $\alpha$
gives that the natural map $R \rightarrow R_{\alpha} \times_{p_{i,\alpha},Z_{\alpha}} Z$
over $p_i:R \rightarrow Z$ is an isomorphism for all $\alpha$.  
But $R_{\alpha} = Z_{\alpha} \times_{Y_{\alpha}} Z_{\alpha}$, 
so taking $i=2$ gives $$R \simeq Z_{\alpha} \times_{Y_{\alpha}} Z$$
for all $\alpha$.  Hence, passing to the limit on $\alpha$ gives
$R = Z \times_Y Z$.   In other words, $R \rightrightarrows Z$ is
an \'etale chart in qcqs schemes for the algebraic space $Y$.  

Our problem is now reduced to showing that the natural map of algebraic spaces 
$\phi:Z/R \rightarrow \invlim (Z_{\alpha}/R_{\alpha})$ 
is an isomorphism, where the inverse system
of algebraic spaces $\{Z_{\alpha}/R_{\alpha}\} = \{Y_{\alpha}\}$
has affine and schematically dominant transition maps.  The map $\phi$ is affine
and schematically dominant since Corollary \ref{cor5} implies
that each map $Z/R \rightarrow Z_{\alpha}/R_{\alpha}$ is affine and schematically
dominant.   But the qcqs \'etale coverings $Z_{\alpha} \rightarrow Y_{\alpha} =
Z_{\alpha}/R_{\alpha}$ are cartesian with respect to change in $\alpha$, so passing to the limit
gives that $Z = \invlim Z_{\alpha}$ is a qcqs \'etale scheme cover of $\invlim (Z_{\alpha}/R_{\alpha})$.
This covering by $Z$ is compatible with $\phi$, so the affine map $\phi$ is an \'etale surjection. 
Since $R = \invlim R_{\alpha}$ inside of $Z \times Z = \invlim (Z_{\alpha} \times Z_{\alpha})$,
it follows that $\phi$ is a monomorphism.  Being affine and \'etale, it is therefore
also a quasi-compact open immersion.  But $\phi$ is an \'etale cover, so it is an isomorphism.
\end{proof}

To apply Proposition \ref{prop8} repeatedly in the context of
Theorem \ref{rgstrat}, we require one more lemma. 

\begin{lemma}\label{lemma10} Let $S$ be a qcqs algebraic space,
and choose a finite rising chain $\{U_i\}$ of quasi-compact open subspaces in $S$
and quasi-compact \'etale scheme covers $\varphi_i:Y_i \rightarrow U_i$
with separated $Y_i$ as in Theorem $\ref{rgstrat}$.
For each $i > 0$, consider the diagram of cartesian squares
$$\xymatrix{
{U'_{i-1}} \ar[d] \ar[r] & {Y_i} \ar[d]^-{\varphi_i} & {Z_i} \ar[l] \ar@{=}[d] \\
{U_{i-1}} \ar[r] & {U_i} & {Z_i} \ar[l] }$$
where $Z_i := U_i - U_{i-1}$ endowed with the reduced structure. 

The \'etale equivalence relation in schemes $R_i = Y_i \times_{U_i} Y_i 
\rightrightarrows Y_i$ is the extension along the diagonal 
$($in the sense of Lemma $\ref{pushoutlemma}$$)$ of the
\'etale equivalence relation $U'_{i-1} \times_{U_{i-1}} U'_{i-1}
\rightrightarrows U'_{i-1}$. 
\end{lemma}

\begin{proof}
The subfunctor $R_i \subseteq Y_i \times Y_i$
contains the subfunctors $\Delta(Y_i)$ and $U'_{i-1} \times_{U_{i-1}} U'_{i-1}$
which overlap along the common open subfunctor
$\Delta(U'_{i-1})$ (openness in $U'_{i-1} \times_{U_{i-1}} U'_{i-1}$
due to $U'_{i-1} \rightarrow U_{i-1}$ being \'etale).  Our aim is to prove that the inclusion
$$\eta_i:\Delta(Y_i) \coprod_{\Delta(U'_{i-1})} (U'_{i-1} \times_{U_{i-1}} U'_{i-1})
\subseteq R_i$$
between subfunctors of $Y_i \times Y_i$ is an isomorphism.   Restricting over
the open subscheme $U'_{i-1} \times U'_{i-1}$ clearly gives an isomorphism,
and since $\varphi_i$ is \'etale and
separated we see that $\Delta(Y_i)$ is an open and closed subscheme of $R_i$.
Thus, $\eta_i$ is an open immersion of schemes, so it suffices to check equality
on geometric fibers over $U_i$.  Over $U_{i-1}$ the situation is clear, and over
$U_i - U_{i-1} = Z_i$ it is also clear since $\varphi_i$ restricts to an isomorphism over $Z_i$
(so the part of $R_i$ lying over $Z_i
\subseteq U_i$ is contained in $\Delta(Y_i)$ on geometric points). 
\end{proof}

Now we are finally in position to prove Theorem \ref{absapprox}. 

\begin{proof}[Proof of Theorem $\ref{absapprox}$]
  Fix a stratification and associated \'etale coverings as in Theorem
  \ref{rgstrat}.  We shall prove Theorem \ref{absapprox} by induction
  on $r$, the case $r = 0$ being the trivial case of empty $S$.
  In general, 
  by  induction we may assume $r \ge 1$ and that $U_{r-1}$ is approximable.  
  (Note that if $r = 1$ then $S_{\rm{red}}$ is a scheme, but we do not yet know
  that this forces $S$ to be a scheme in general.  Hence, approximation
  for schemes \cite[Thm.~C.9]{tt} does not suffices to settle the case $r=1$.)  
  By Lemma
  \ref{lemma10}, the open immersion $U_{r-1} \hookrightarrow U_r = S$
  arises along the bottom side of a pushout diagram as in (\ref{ypush}). 
   Thus, by Proposition \ref{prop8} the approximability
  of $S$ follows from that of $U_{r-1}$.  This completes the proof.
\end{proof}

\begin{remark} 
In case the reader is wondering where the scheme case of Theorem \ref{absapprox}
(i.e., \cite[Thm.~C.9]{tt}) is lurking in our proof for the case of algebraic spaces, it is
used in the proof of Lemma \ref{scheme}, which in turn is an essential ingredient in
proofs of subsequent results in \S\ref{limprop} that were used in our treatment of the case
of algebraic spaces. 
\end{remark}

\begin{corollary}\label{corred}
Let $S$ be an algebraic space.  If $S_{\rm{red}}$ is a scheme then $S$ is a scheme.
\end{corollary}

\begin{proof}
Working Zariski-locally on $S_{\rm{red}}$ is the same as working Zariski-locally on $S$,
so we may arrange that $S_{\rm{red}}$ is an affine scheme.  Hence,
$S_{\rm{red}}$ is quasi-compact and separated, so $S$ is quasi-compact and separated.
By Theorem \ref{absapprox}, we may therefore write $S \simeq \invlim S_i$ where
$\{S_i\}$ is an inverse system of algebraic spaces of finite presentation over $\Z$.   Since 
$\invlim (S_i)_{\rm{red}} \simeq S_{\rm{red}}$, this limit is an affine scheme.  
Thus, by Lemma \ref{scheme}, there is an $i_0$ such that $(S_i)_{\rm{red}}$ is a scheme for
all $i \ge i_0$.  But each $S_i$ is a noetherian algebraic space, so by 
\cite[III,~Thm.~3.3]{knutson} it follows that $S_i$ is a scheme for all $i \ge i_0$.  Hence, $S$ is a scheme
since each map $S \rightarrow S_i$ is affine.
\end{proof}

As another application of Theorem \ref{absapprox}, we can prove the following result via reduction
to the finitely presented case over $\Z$:

\begin{theorem}[Deligne]\label{delthm}
 Let $X \rightarrow S$ be a separated map of quasi-separated algebraic spaces
equipped with a left $S$-action on $X$ by a finite group $G$. The 
quotient stack $[X/G]$ admits a coarse moduli space $\pi :[X/G]\rightarrow X/G$ such that the following hold:
\begin{enumerate}
\item [(i)] $X/G$ is separated over $S$.
\item [(ii)] The projection morphism $\pi _X:X\rightarrow X/G$ is affine and integral.
\item [(iii)] If $U\rightarrow X/G$ is a separated \'etale morphism, then the induced $G$-equivariant map
$$
U\times _{X/G}X\rightarrow U
$$
identifies $U$ with the quotient $(U\times _{X/G}X)/G$.
\item [(iv)] The pullback map
$$
\pi _X^{-1}:(\text{open subspaces of $X/G$})\rightarrow (\text{$G$-invariant open subspaces of $X$})
$$
is a bijection.
\end{enumerate}
\end{theorem}

A stronger version, allowing non-constant $G$, is given in \cite[Thm.~5.4]{rydh-qt}. 
In our proof below, we rely on some results proved in the unpublished note \cite{conrad3}, but
this seems permissible since \cite{rydh-qt} will be published and establishes a more general result. 

\begin{proof}
By \cite[1.1]{conrad3} the result holds when $X\rightarrow S$ is locally of finite presentation, and by 
the proof of \cite[3.1]{conrad3} (using an argument from \cite[V,~4.1]{sga3}) it also holds when
$X$ is an affine scheme (for which it suffices to treat the case $S = \Spec \Z$). 
In fact, the proof shows that if $X = \Spec (A)$, then $X/G = \Spec (A^G)$; the verification of the
universal property uses that we only consider maps to quasi-separated algebraic spaces
in the definition of $X/G$ in \S\ref{gabberlemma}.

The general case can be deduced from these two cases as follows.  By property (iii) and the existence of quotients
by \'etale equivalence relations in the category of 
algebraic spaces, to prove the existence we may work \'etale locally  on $S$.  Hence, 
we may assume that $S$ is an affine scheme, so in particular $X$ is separated.  Furthermore,
by (iii) and (iv) we can use an elementary gluing argument for the Zariski topology
(as at the start of \cite[\S2]{conrad}) to reduce to the case 
that $X$ is also quasi-compact. We apply absolute noetherian approximation
in Theorem \ref{absapprox} to the quasi-compact separated algebraic space $X$ (over $\Spec \Z$) to get
an isomorphism 
$$
X \simeq \varprojlim _{\lambda }X_\lambda, 
$$
where $\{X_\lambda \}$ is an inverse system of separated algebraic spaces 
of finite presentation (equivalently, finite type) over $\Z$ with affine transition maps.   
In particular, we can find an affine morphism
$$
X\rightarrow Y'
$$
with $Y'$ separated and finitely presented over $\Z$.    Fix such a morphism.  

The algebraic space $Y:= S\times _{\Spec \Z}Y'$ 
is affine over $Y'$, so we obtain an affine morphism
$$
\gamma :X\rightarrow Y
$$
over $S$ to an algebraic space that is separated
and of finite presentation over $S$.  Let $Z$ denote the $|G|$-fold fiber product of $Y$ with itself over $S$: 
$$
Z:= \prod _{g\in G}Y.
$$
There is a left action of $G$ on $Z$ given by permuting the factors via $g_0.(y_g)_{g \in G} =
(y_{gg_0})_{g \in G}$.  Moreover, we have a $G$-equivariant map
$$
f :X\rightarrow Z, \ \ x\mapsto (\gamma (g.x))_{g\in G}.
$$
This map $f $ is affine since it is equal to the composite
$$
\xymatrix{
X\ar[r]^-{\Delta }& \prod _{g\in G}X\ar[r]^-{\prod {\gamma \circ g}}& \prod _{g\in G}Y}
$$
in which $\Delta $ is affine since $X$ is $S$-separated. 

By the settled case when the structure map to $S$ is of finite presentation,
the quotient $Z/G$ exists as an $S$-separated algebraic space
and the projection map $Z\rightarrow Z/G$ is
affine (by (ii)). 
We conclude that we can find an $S$-morphism
$$
\tau :[X/G]\rightarrow W
$$
from the stack-quotient $[X/G]$ to an $S$-separated algebraic space $W$ such that the composite map
$$
X\rightarrow [X/G]\rightarrow W
$$
is affine.  We may replace $S$ with $W$, and then by working \'etale-locally on $S$ as before we are 
reduced to the settled case when $X$ is an affine scheme. 
\end{proof}

%The following result will be useful when reducing the proof of Theorem \ref{mainresult}
%to the case of objects of finite type over $\Z$. 

%\begin{corollary}\label{afflimit}
%Let $\Lambda$ be a ring
%and $\{X_{\alpha}\}$ and $\{S_{\alpha}\}$ be inverse systems of algebraic
%spaces of finite presentation over $\Lambda$
%with affine transition maps, and suppose there are given compatible maps
%$f_{\alpha}:X_{\alpha} \rightarrow S_{\alpha}$ over $\Lambda$ satisfying the cartesian condition
%that $X_{\beta} \rightarrow S_{\beta} \times_{S_{\alpha}} X_{\alpha}$
%is an isomorphism whenever $\beta \ge \alpha$.  Let
%$X = \invlim X_{\alpha}$ and $S = \invlim S_{\alpha}$,
%and let $f = \invlim f_{\alpha}:X \rightarrow S$.

%The map $f$ is affine $($resp. a closed immersion, resp. separated$)$ if and only if
%$f_{\alpha}$ is affine $($resp. a closed immersion, resp. separated$)$ for sufficiently large $\alpha$.
%\end{corollary}

%\begin{proof}
%By passing to diagonal maps it suffices to treat the
%``$\Rightarrow$'' direction for the property of being affine or a closed immersion.
%Fix an $\alpha_0$.  We may work \'etale-locally over $S_{\alpha_0}$,
%so we may assume that all $S_{\alpha}$ and $S$ are affine schemes, so $X$ is affine.
%By Lemma \ref{scheme}, since $\invlim X_{\alpha} = X$ is an affine scheme,
%so is $X_{\alpha}$ for large $\alpha$.  This settles
%the affine case, and the case of closed immersions is now clear. 
%\end{proof}

\subsection{Finite type and finite presentation}\label{fpresec}

In \cite[Thm.~4.3]{conrad} it is proved that if
$X \rightarrow S$ is a map of finite type between
qcqs schemes then there is a closed immersion $i:X \hookrightarrow \overline{X}$ over $S$ into
a finitely presented $S$-scheme $\overline{X}$, and that $\overline{X}$ can be taken to be separated
over $S$ if $X$ is separated over $S$.  This is the trick that, together with absolute
noetherian approximation for qcqs schemes, allows one to reduce the proof of Nagata's theorem
in the general scheme case to the case of schemes of finite type over $\Z$. 
We require an analogue for algebraic spaces, so we now aim to prove:

\begin{theorem}\label{fpresred}
Let $f:X \rightarrow S$ be a map of finite type between qcqs algebraic spaces.
There exists a closed immersion $i:X \hookrightarrow \mathscr{X}$
over $S$ into an algebraic space $\mathscr{X}$ of finite presentation over $S$.
If $X$ is $S$-separated then $\mathscr{X}$ may be taken to be $S$-separated.
\end{theorem}

To prove Theorem \ref{fpresred}, we first need a gluing result for closed subspaces
of algebraic spaces of finite presentation over a qcqs algebraic space $S$.
Consider a commutative diagram of $S$-maps
\begin{equation}\label{pushdiagram}
  \xymatrix{
    U' \ar@{->>}[d]_-{q} \ar[r]^-{j'} \ar[dr]^-{i'_1} & X'  \ar[dr]^-{i} & \\
    U \ar[dr]_-{i_1} & \mathscr{U}' \ar[r]_-{j} \ar@{-->}[d]_-{?}^-{\pi} &  \mathscr{X}' \\
    & \mathscr{U} &}
\end{equation}
in which $q$ is a quasi-compact separated \'etale cover,  $j$ is an open immersion, the maps
$i_1$, $i'_1$, and $i$ are closed immersions into algebraic spaces
that are finitely presented and separated over $S$, and the top part
is cartesian $($so $j'$ is an open immersion$)$.  We wish to study the possibility that there
exists a suitable map $\pi$ as indicated in (\ref{pushdiagram}). 

\begin{lemma}\label{fpresredpre1}
With notation and hypotheses as in
$(\ref{pushdiagram})$, 
let $X = U \coprod_{U'} X'$ be the pushout of the upper left triangle formed by $j'$ and $q$.
If there is a quasi-compact separated \'etale map $\pi$
as shown in \eqref{pushdiagram} that makes the left part cartesian
then the natural $S$-map from $X$ to the algebraic space 
pushout $\mathscr{X} = \mathscr{U} \coprod_{\mathscr{U}'} \mathscr{X}'$ of $j$ along $\pi$
is a closed immersion, and
$\mathscr{X}$ is finitely
presented over $S$.
\end{lemma}
\begin{proof}
%   By Lemma \ref{pushoutlemma}, the algebraic space $X$ over $S$ is
%   implicit in this diagram in the sense that the natural cartesian
%   diagram
% \begin{equation}\label{pushx}
% \xymatrix{
% U' \ar[r]^-{j'} \ar@{->>}[d]_-{q} & X' \ar[d] \\ U \ar[r] & X}
% \end{equation}
% arising from the definitions (with bottom side equal to the quasi-compact open immersion
% $X_{r-1} \hookrightarrow X$) expresses $X$ as a pushout.  

Given $\pi$, form the cartesian pushout
diagram in qcqs algebraic spaces over $S$
\begin{equation}\label{firstdiag}
\xymatrix{
\mathscr{U}' \ar[r]^-{j} \ar[d]_-{\pi} & \mathscr{X}' \ar[d] \\
\mathscr{U} \ar[r] & \mathscr{X}}
\end{equation}
as in Lemma \ref{pushoutlemma}, so the bottom side is a quasi-compact open immersion
(as $j$ is) and the right side is a quasi-compact separated \'etale surjection. 
In particular, $\mathscr{X}$ is finitely presented over $S$ since $\mathscr{X}'$ is
finitely presented over $S$.
Consider the $S$-map of pushouts
$$X = U \coprod_{U'} X' \rightarrow \mathscr{U} \coprod_{\mathscr{U}'} \mathscr{X}' = \mathscr{X}.$$
It remains to  prove that this is a closed immersion.

It is elementary to check that for $j \in \{1,2\}$, the diagram
$$\xymatrix{
U' \times_U U' \ar[d]_-{p_j} \ar[r] & \mathscr{U}' \times_{\mathscr{U}} \mathscr{U}' \ar[d]^-{p_j} \\
X' \ar[r]_-{i} & \mathscr{X}'}$$
is cartesian, due to the cartesian property of both parallelogram sub-diagrams in
(\ref{pushdiagram}).   Since $X = X'/R'$ and $\mathscr{X} = \mathscr{X}'/\mathscr{R}'$ for
$$R' := (U' \times_U U') \coprod X',\,\,\,\,\,
\mathscr{R}' := (\mathscr{U}' \times_{\mathscr{U}} \mathscr{U}') \coprod \mathscr{X}',$$
we may apply Lemma \ref{lemma3} to infer that the commutative diagram
$$\xymatrix{
X' \ar[r]^-{i} \ar[d] & \mathscr{X}' \ar[d] \\
X \ar[r] & \mathscr{X}}$$
with \'etale surjective vertical maps is cartesian.   Hence, the bottom side is a closed immersion since
the top side is a closed immersion.
\end{proof}

Given an arbitrary diagram of type \eqref{pushdiagram}
(with hypotheses on the maps as indicated there), 
the existence of $\pi$ as in Lemma \ref{fpresredpre1}
is quite subtle (and likely false).  However, we can always modify $\mathscr U$, $\mathscr U'$, and $\mathscr X'$ so that the resulting diagram admits a $\pi$.  More precisely, we have the following.

\begin{proposition}\label{gluingyourmama}
  Given a diagram \eqref{pushdiagram}, there is another one with the same $j'$ and $q$ 
  and for which an arrow $\pi$ as in Lemma $\ref{fpresredpre1}$ exists. 
\end{proposition}
This says that, given a pushout diagram in the category of algebraic spaces of the type encountered in the upper left of \eqref{pushdiagram} with $j'$ and $q$, 
if the objects $U$, $U'$, and $X'$ individually admit closed immersions 
over $S$ into finitely presented algebraic spaces separated over $S$ 
(and satisfy a compatibility as expressed by the auxiliary map $j$) 
then we can choose such immersions so that an arrow $\pi$ as in Lemma \ref{fpresredpre1} exists.
\begin{proof}[Proof of Proposition $\ref{gluingyourmama}$] The proof proceeds in several steps.

{\bf Step 1 (replacing $\mathscr U$)}: Let us first show that there is a finitely presented closed
subspace $\mathscr{V}$ in  $\mathscr{U}$ which is better 
than $\mathscr{U}$ in the sense that it not only contains $U$ as a closed subspace
but also admits a quasi-compact, \'etale, separated covering $h:\mathscr{V}' \rightarrow \mathscr{V}$
whose pullback over $U$ is $q:U' \rightarrow U$.  The defect, which will have to be confronted
in subsequent steps, is that $\mathscr{V}'$ will have no evident connection with $\mathscr{U}'$
(or even with $\mathscr{X}'$) apart from containing $U'$ as a closed subspace. 

\begin{lemma}\label{qetdescend}
There is a cartesian diagram of algebraic spaces over $S$
$$\xymatrix{
U' \ar[d]_-{q} \ar[r] & \mathscr{V}' \ar[d]^-{h} \\
U \ar[r] & \mathscr{V}}$$
such that the bottom arrow is an inclusion of closed subspaces of $\mathscr{U}$
with $\mathscr{V}$ finitely presented over $S$, and 
$h$ is quasi-compact, \'etale, and separated.
\end{lemma}

If $U$ were finitely presented over $S$ then we could take $\mathscr{V} = U$ and $h = q$.
The point of the lemma is to ``spread out'' $q$ over a {\em finitely presented}
closed subspace of $\mathscr{U}$, whereas $U$ is merely finite type (and quasi-separated) over $S$
 so its quasi-coherent ideal in $\mathscr{O}_{\mathscr{U}}$ may not be locally finitely generated. 

\begin{proof} 
Consider the given closed immersion $i_1:U \hookrightarrow \mathscr{U}$ over $S$, 
with $\mathscr{U}$ a finitely
presented algebraic space over $S$. Since $\mathscr{U}$ is a qcqs algebraic space, 
the quasi-coherent ideal $\mathscr{I} \subseteq {\ms O}_{\mathscr{U}}$
cutting out $U$ can be expressed as the direct limit $\varinjlim \mathscr{I}_{\lambda}$
of its quasi-coherent subsheaves of finite type \cite[I,~5.7.8]{rg}.
Hence, $U = \invlim \mathscr{U}_{\lambda}$ where
$\mathscr{U}_{\lambda} \hookrightarrow \mathscr{U}$ is cut out by 
$\mathscr{I}_{\lambda}$.  

The inverse system $\{\mathscr{U}_{\lambda}\}$
of finitely presented algebraic spaces over $S$
has affine transition maps and limit $U$, 
and every qcqs algebraic space (such as any $\mathscr{U}_{\lambda}$)
is affine over an algebraic space of
finite presentation over $\Z$ (by Theorem \ref{absapprox}, which was proved
in \S\ref{absnoetherian}).   Thus, by Proposition \ref{lemma7}
and Corollary \ref{mapdescend2} we deduce that 
the quasi-compact \'etale separated cover $U' \rightarrow U$
descends to a quasi-compact \'etale separated
cover $\mathscr{V}' \rightarrow \mathscr{U}_{\lambda_0}$
for some sufficiently large $\lambda_0$.  Rename this $\mathscr{U}_{\lambda_0}$
as $\mathscr{V}$.
\end{proof}

We may and do replace $\mathscr{U}$ in (\ref{pushdiagram}) with $\mathscr{V}$, but we emphasize again that 
$\mathscr{V}'$ is unrelated to $\mathscr{U}'$ (except that both contain $U'$ as a closed subspace).

Returning to the task of constructing $\pi$ as in \eqref{pushdiagram} after a suitable
change in $\mathscr{U}'$ (and leaving $j'$ and $q$ unchanged), 
the strategy is to make an initial change in 
$\mathscr{U}'$ (as a quasi-compact open subspace of $\mathscr{X}'$
meeting $X'$ in $U'$) and to then show that by replacing
$\mathscr{X}'$ and $\mathscr{U} = \mathscr{V}$ with suitable finitely presented
closed subspaces (respectively containing $X'$ and $U$) and 
replacing the $\mathscr{U}$-spaces $\mathscr{U}'$ and $\mathscr{V}'$ with the respective pullback closed
subspaces  containing $U'$ we eventually get to a situation in which
we {\em can} identify $\mathscr{V}'$ and $\mathscr{U}'$ over $S$ in a manner
that respects how $U'$ is a closed subspace of each.  In such a favorable situation the map
$h:\mathscr{V}' \rightarrow \mathscr{U}$ then serves as the desired map $\pi$.  

To carry this out, we have to overcome the essential difficulty in comparison
with the Zariski-gluing problem faced in the scheme case as considered
in \cite[Thm.~4.3]{conrad}:  whereas $U'$ is open in $X'$, it is only 
\'etale (rather than Zariski-open) over $U$, and so rather than trying
to spread $U'$ to a common open subspace of
$\mathscr{X}'$ and $\mathscr{U}$ (after suitable shrinking on these two spaces)
we are instead trying to spread $U'$ to an open subspace
of $\mathscr{X}'$ that is quasi-compact, \'etale, and separated over $\mathscr{U}$.

Let $\{\mathscr{X}'_{\alpha}\}$ denote the inverse system 
 consisting of all finitely presented closed subspaces
 in $\mathscr{X}'$ containing $X'$ and let
 $\{\mathscr{U}_{\beta}\}$ denote the inverse system of finitely 
presented closed subspaces in $\mathscr{U}$ containing $U$, so 
$$X' = \invlim \mathscr{X}'_{\alpha},\,\,\,U = \invlim \mathscr{U}_{\beta}.$$
Consider the diagram of cartesian squares
$$\xymatrix{
U' \ar[d] \ar[r] & \mathscr{U}'_{\alpha} \ar[d] \ar[r] & \mathscr{U}' \ar[d] \\
X' \ar[r] & \mathscr{X}'_{\alpha} \ar[r] & \mathscr{X}'}$$
over $S$ in which the horizontal maps are closed immersions and the vertical
maps are quasi-compact open immersions. 
Passing to the limit gives $U' \simeq \invlim \mathscr{U}'_{\alpha}$  over $S$, 
with $\{\mathscr{U}'_{\alpha}\}$ an inverse system of finitely presented closed
subspaces of $\mathscr{U}'$ containing $U'$. 

 Likewise, consider the diagram of cartesian squares
$$\xymatrix{
U' \ar[d] \ar[r] & \mathscr{V}'_{\beta} \ar[d] \ar[r]^-{j_{\beta}} & \mathscr{V}' \ar[d]^-{h} \\
U \ar[r] & \mathscr{U}_{\beta} \ar[r] & \mathscr{U}}$$
over $S$ in which the horizontal maps are closed immersions and the vertical
maps are quasi-compact, \'etale, and separated.   Passing to the limit 
gives another limit presentation of $U'$ over $S$, namely $U' = \invlim
\mathscr{V}'_{\beta}$.  The situation we would like to reach is
that for some large $\alpha_0$ and $\beta_0$
there is an $S$-isomorphism $\mathscr{V}'_{\beta_0} \simeq \mathscr{U}'_{\alpha_0}$
respecting the closed immersions of $U'$ into each side, as we could
then rename $\mathscr{V}'_{\beta_0}$ as $\mathscr{V}'$ and $\mathscr{U}'_{\alpha_0}$
as $\mathscr{U}'$ to get to the case
$\mathscr{V}' = \mathscr{U}'$ that provides $\pi$ as desired.  We will not find such a pair
$(\alpha_0, \beta_0)$ after making changes in various
auxiliary spaces, but we will do something close to this which is sufficient for the purpose
of proving Proposition \ref{gluingyourmama}. 

{\bf Step 2 (shrinking $\mathscr{U}'$ and $\mathscr{V}'$ so that they become related)}: 
A key observation is that we may replace $\mathscr{U}'$ with any finitely presented closed subspace containing
$U'$.    To justify this, we have to prove that any finitely presented
closed subspace $\mathscr{U}'_1$ in $\mathscr{U}'$ containing $U'$ extends 
to a finitely presented closed subspace $\mathscr{X}'_1$ in $\mathscr{X}'$
containing $X'$.   (Recall that $\mathscr{U}'$ is open in $\mathscr{X}'$
and $X'$ is closed in $\mathscr{X}'$, with
$\mathscr{U}' \cap X' = \mathscr{U}' \times_{\mathscr{X}'} X' = U'$.)
This is a standard extension argument with quasi-coherent sheaves, except
that we are working with algebraic spaces rather than with schemes, and it goes as follows. 
Let $\mathscr{I}_{X'} \subseteq {\ms O}_{\mathscr{X}'}$ denote the quasi-coherent ideal
that cuts out $X'$ inside of $\mathscr{X}'$.  
The finitely presented closed subspace $\mathscr{U}'_1$ in the quasi-compact open subspace
$\mathscr{U}' \subseteq \mathscr{X}'$ is cut out by a finite type 
quasi-coherent ideal sheaf $\mathscr{I} \subseteq {\ms O}_{\mathscr{U}'}$, 
and since $\mathscr{U}' \cap X' = U' \subseteq \mathscr{U}'_1$ 
we have $\mathscr{I} \subseteq \mathscr{I}_{X'}|_{\mathscr{U}'}$.
Let $\mathscr{K} \subseteq {\ms O}_{\mathscr{X}'}$ be the quasi-coherent ideal
corresponding to the schematic closure of $\mathscr{U}'_1$ inside of $\mathscr{X}'$.
(Explicitly, if $j:\mathscr{U}' \hookrightarrow \mathscr{X}'$ denotes the quasi-compact
open immersion then $\mathscr{K}$ is  the preimage of $j_{\ast}(\mathscr{I})$
under ${\ms O}_{\mathscr{X}'} \rightarrow j_{\ast}({\ms O}_{\mathscr{U}'})$.)
This is a quasi-coherent ideal sheaf whose intersection with
$\mathscr{I}_{X'}$ inside of ${\ms O}_{\mathscr{X}'}$
restricts to the finite type $\mathscr{I}$ on the quasi-compact open $\mathscr{U}'$. 
On the qcqs algebraic space $\mathscr{X}'$ we can express
$\mathscr{K}$ as the direct limit of its finite type quasi-coherent subsheaves
\cite[I,~5.7.8]{rg}, so the finite type quasi-coherent ideal $\mathscr{I}$ on $\mathscr{U}'$ extends to 
to a finite type quasi-coherent ideal $\mathscr{J} \subseteq \mathscr{I}_{X'}$. 
The finitely presented closed subspace $\mathscr{X}'_1 \hookrightarrow
\mathscr{X}'$ cut out by $\mathscr{J}$ fits into a commutative diagram
\begin{equation}\label{xumaps}
\xymatrix{
X' \ar[r] & \mathscr{X}'_1 \ar[r] & \mathscr{X}' \\
U' \ar[u] \ar[r] & \mathscr{U}'_1 \ar[u] \ar[r] & \mathscr{U}' \ar[u]}
\end{equation} 
in which the vertical maps are open immersions, the horizontal maps are
closed immersions, the right square is cartesian by definition of
$\mathscr{X}'_1$, and the outside edge is cartesian since
$X' \cap \mathscr{U}' = U'$ as quasi-compact subspaces of $\mathscr{X}'$. 
Hence, the {\em left} square in (\ref{xumaps}) is cartesian, so
$\mathscr{U}'_1$ and $\mathscr{X}_1$ are as required.  

To make a useful 
initial shrinking of $\mathscr{U}'$ as just considered, we shall use $\mathscr{V}'$ as a guide. 
Since $\mathscr{V}'$ is finitely presented over $S$, 
by  Proposition \ref{fpres} 
we have 
$\Hom_S(U', \mathscr{V}') = \varinjlim \Hom_S(\mathscr{U}'_{\alpha}, \mathscr{V}')$.
Applying this to the closed immersion $\invlim \mathscr{U}'_{\alpha} = U' \hookrightarrow \mathscr{V}'$ over $S$ 
as in Lemma \ref{qetdescend}
gives an $\alpha_0$ for which there is a commutative diagram of $S$-maps
\begin{equation}\label{raoultsmother1}
\xymatrix{
{U'} \ar[rr] \ar[dr] && \mathscr{V}' \\
 & \mathscr{U}'_{\alpha_0} \ar[ur]_-{\phi} & }
 \end{equation}
in which the left diagonal and top arrows are the canonical closed immersions.
Thus, upon renaming such a $\mathscr{U}'_{\alpha_0}$ as $\mathscr{U}'$
(as we may), we get a commutative diagram of $S$-maps 
\begin{equation}\label{raoultsmother2}
\xymatrix{
{U'} \ar[rr] \ar[dr] && \mathscr{V}' \\
 & \mathscr{U}' \ar[ur]_-{\phi} & }
 \end{equation}
in which the unlabeled maps are the canonical closed immersions.  
Likewise, since the new $\mathscr{U}'$ is finitely presented over $S$ we have
 $$\Hom_S(U', \mathscr{U}') = \varinjlim \Hom_S(\mathscr{V}'_{\beta}, \mathscr{U}'),$$
 so there is a $\beta_0$ for which we may construct a commutative diagram
 \begin{equation}\label{raoultsmother3}
 \xymatrix{
{U'} \ar[rr] \ar[dr] && \mathscr{U}' \\
 & \mathscr{V}' _{\beta_0} \ar[ur]_-{\psi_{\beta_0}} & }
 \end{equation}
  over $S$ in which the unlabeled maps are the canonical closed immersions. 

For $\beta \ge \beta_0$, define $\psi_{\beta}:\mathscr{V}'_{\beta} \rightarrow \mathscr{U}'$
to be the composition of $\psi_{\beta_0}$ with the closed immersion
$\mathscr{V}'_{\beta} \hookrightarrow \mathscr{V}'_{\beta_0}$.
We concatenate (\ref{raoultsmother2}) and the variant of (\ref{raoultsmother3})
for $\beta \ge \beta_0$ to get  the commutative diagram
$$\xymatrix{
\mathscr{V}'_{\beta} \ar[r]^-{\psi_{\beta}} & \mathscr{U}' \ar[r]^-{\phi} & \mathscr{V}' \\
& U' \ar[ul] \ar[u] \ar[ur] & }$$
over $S$ in which the maps from $U'$ are the canonical closed immersions. 
Since the family of closed immersions $U' \rightarrow \mathscr{V}'_{\beta}$ over $S$ 
becomes the $S$-isomorphism $U' \simeq \invlim \mathscr{V}'_{\beta}$ upon passing to the limit, 
the family of $S$-maps $\phi \circ \psi_{\beta}:\mathscr{V}'_{\beta} \rightarrow \mathscr{V}'$
and the family of $S$-maps $j_{\beta}:\mathscr{V}'_{\beta} \hookrightarrow \mathscr{V}'$
given by the canonical closed immersions have limit maps $U' \rightrightarrows
\mathscr{V}'$ which coincide (with the canonical closed immersion).
But $U' = \invlim V'_{\beta}$, so Proposition \ref{fpres} gives
$$\Hom_S(U', \mathscr{V}') = \varinjlim \Hom_S(\mathscr{V}'_{\beta}, \mathscr{V}').$$
Thus, for some $\beta_1 \ge \beta_0$ the $S$-map 
$\phi \circ \psi_{\beta_1}:\mathscr{V}'_{\beta_1} \rightarrow \mathscr{V}'$
is equal to the canonical closed immersion $j_{\beta_1}$.

Now consider the diagram of $S$-maps (with $\psi := \psi_{\beta_1}$) 
\begin{equation}\label{maxdiag}
\xymatrix{
\mathscr{V}'_{\beta_1} \ar[r]^-{\psi} \ar@{-->}[dr]_-{\psi_1} & \mathscr{U}' \ar[r]^-{\phi} &
\mathscr{V}' \\
& \mathscr{U}'_1 \ar[u]_-{i_{\beta_1}} \ar[r]^-{\phi_1} & \mathscr{V}'_{\beta_1} \ar[u]_-{j_{\beta_1}} \\
U' \ar[uu] \ar@{-->}[ur]^-{?} \ar[urr] & &}
\end{equation}
in which the new algebraic space
$\mathscr{U}'_1$ and the maps $i_{\beta_1}$ and $\phi_1$ from it to
$\mathscr{U}'$ and $\mathscr{V}'_{\beta_1}$ are 
{\em defined} to make the right square be cartesian
(so $\mathscr{U}'_1$ is a finitely presented closed subspace of
$\mathscr{U}'$) and the composite map $\phi \circ \psi$ across the top 
equals the map $j_{\beta_1}$ along the right side.  Thus, there is a unique 
 $\psi_1$ making the upper triangle commute (uniqueness because
 $i_{\beta_1}$ is monic), and the composite
 $\phi_1 \circ \psi_1$ is the identity map either by construction or because 
 $$j_{\beta_1} \circ (\phi_1 \circ \psi_1) = \phi \circ i_{\beta_1} \circ \psi_1 = 
 \phi \circ \psi = j_{\beta_1}$$
 with $j_{\beta_1}$ monic. 
The canonical closed immersion from $U'$ to $\mathscr{V}'_{\beta_1}$
is used along the left and bottom sides of the diagram, and we 
define the dotted arrow from $U'$ to $\mathscr{U}'_1$ by composing
the left side with $\psi_1$.  
In particular, $\mathscr{U}'_1$ as a (finitely presented) closed subspace of $\mathscr{U}'$
contains $U'$. Since $\phi_1 \circ \psi_1 = {\rm{id}}_{\mathscr{V}'_{\beta_1}}$,
the bottom triangle also commutes.  Thus, we have filled in all of the dotted arrows to
make (\ref{maxdiag}) a commutative diagram over $S$. 

Using the argument at the beginning of Step 2, we can extend 
$\mathscr{U}'_1$ 
to a finitely presented closed subspace $\mathscr{X}'_1$ in $\mathscr{X}'$
containing $X'$.    We may replace 
$\mathscr{U}' \hookrightarrow  \mathscr{X}'$ with 
$\mathscr{U}'_1 \hookrightarrow \mathscr{X}'_1$,
$\mathscr{V} = \mathscr{U}$ with $\mathscr{U}_{\beta_1}$, 
$\mathscr{V}'$ with $\mathscr{V}'_{\beta_1} = \mathscr{V}' \times_{\mathscr{U}} \mathscr{U}_{\beta_1}$,
and the maps $\phi$ and $\psi$ with $\phi_1$ (renamed as $s$) and $\psi_1$
(renamed as $t$) respectively to arrive at the case when
the conclusion of Lemma \ref{qetdescend} still holds with
$\mathscr{V} = \mathscr{U}$ but there is also a commutative diagram
\begin{equation}\label{stdiagram}
\xymatrix{
\mathscr{V}' \ar[r]^-{t} & \mathscr{U}' \ar[r]^-{s} & \mathscr{V}' \\
& U' \ar[ul] \ar[u] \ar[ur] & }
\end{equation}
whose composite across the top is the identity map (and whose non-horizontal arrows
are the canonical closed immersions).  
If $t \circ s = {\rm{id}}_{\mathscr{U}'}$ then $s$ and $t$ would be inverse
isomorphisms respecting the closed immersions from $U'$, so we would be done
(taking $\pi$ to be $\mathscr{U}' \simeq \mathscr{V}' \stackrel{h}{\rightarrow} \mathscr{V} =
\mathscr{U}$). 
It remains to handle the possibility $t \circ s \ne {\rm{id}}_{\mathscr{U}'}$.  

{\bf Step 3 (making inverse maps)}:
We shall construct suitable further shrinkings of
$\mathscr{U}'$ and $\mathscr{V}'$ compatibly with $s$ and $t$ simultaneously so
that we can get to the case where $t \circ s = {\rm{id}}_{\mathscr{U}'}$. 
The commutativity of (\ref{stdiagram}) implies that the composite $S$-map
$$U' \hookrightarrow \mathscr{U}' \stackrel{s}{\rightarrow} \mathscr{V}'
\stackrel{t}{\rightarrow} \mathscr{U}'$$
is the canonical closed immersion.  But $U' = \invlim \mathscr{U}'_{\alpha}$
and the target $\mathscr{U}'$ is finitely presented over $S$, so
Proposition \ref{fpres}  provides
an $\alpha_2$ such that  the composite $S$-map
$$\mathscr{U}'_{\alpha_2} \hookrightarrow \mathscr{U}' \stackrel{s}{\rightarrow} \mathscr{V}'
\stackrel{t}{\rightarrow} \mathscr{U}'$$
is the canonical closed immersion. 
Let $\mathscr{U}'_2 = \mathscr{U}'_{\alpha_2}$
and $\mathscr{V}'_2 = t^{-1}(\mathscr{U}'_{\alpha_2})$,
so it makes sense to restrict $s$ and $t$ to get $S$-maps
$$\mathscr{U}'_2 \stackrel{s_2}{\rightarrow} \mathscr{V}'_2 \stackrel{t_2}{\rightarrow}
\mathscr{U}'_2$$
whose composition is the identity map (as may be checked upon composing on the right
with $\mathscr{U}'_2 \hookrightarrow \mathscr{U}'$).   Moreover,
since $s \circ t = {\rm{id}}_{\mathscr{U}'}$ we have that $s_2 \circ t_2 =
{\rm{id}}_{\mathscr{V}'_2}$.    Hence, $s_2$ and $t_2$ are inverse $S$-isomorphisms!  But beware
that we are not done, since the commutative square
$$\xymatrix{
U' \ar[r] \ar[d]_-{q} & \mathscr{V}'_2 \ar[d]^-{h|_{\mathscr{V}'_2}} \\
U \ar[r] & \mathscr{U}}$$
may fail to be cartesian, and more seriously
the map $h|_{\mathscr{V}'_2}$ may not be
\'etale (though it is certainly quasi-compact and separated).
Hence, we cannot simply rename
$\mathscr{V}'_2$ as $\mathscr{V}'$ without destroying the  properties that make 
$\mathscr{V}'$ useful in the first place, as in Lemma \ref{qetdescend}. 

We can fix this difficulty as follows.  Since
$U = \invlim \mathscr{U}_{\alpha}$ for the inverse system
$\{\mathscr{U}_{\alpha}\}$ of finitely presented
closed subspaces of $\mathscr{U}$, the cartesian property in Lemma \ref{qetdescend}
implies that $U' = \invlim h^{-1}(\mathscr{U}_{\alpha})$.  Hence, there exists 
$\alpha_3 \ge \alpha_2$ such that the finitely
presented closed subspace $\mathscr{V}'_3 := h^{-1}(\mathscr{U}_{\alpha_3})$ in
$\mathscr{V}'$ is contained in 
$\mathscr{V}'_2$ because $\mathscr{V}'_2$ is finitely presented inside of $\mathscr{V}'$. 
Obviously the commutative diagram
$$\xymatrix{
U' \ar[r] \ar[d] & \mathscr{V}'_3 \ar[d] \\
U \ar[r] & \mathscr{U}_{\alpha_3}}$$
is cartesian and it is \'etale along the right side.  
Thus, if we let $\mathscr{U}'_3 \subseteq \mathscr{U}'_2$ correspond
to $\mathscr{V}'_3 \subseteq \mathscr{V}'_2$ via the inverse $S$-isomorphisms
$s_2$ and $t_2$ then we can replace $\mathscr{U}$ with
$\mathscr{U}_{\alpha_3}$ and use $\mathscr{U}'_3$ in the role of
$\mathscr{U}'_1$ back in the construction of
$\mathscr{X}'_1$ in (\ref{xumaps})
to get a finitely presented closed subspace $\mathscr{X}'_3
\subseteq \mathscr{X}'$ for which the preceding cartesian difficulty does not arise. 
Hence, the $S$-isomorphism $\mathscr{U}'_3 \simeq \mathscr{V}'_3$ 
via $s_2$ and $t_2$ has the required properties to provide the equality
$\mathscr{U}' = \mathscr{V}'$ that we sought to construct. 
\end{proof}

We are now in a position to prove Theorem \ref{fpresred}.

\begin{proof}[Proof of Theorem $\ref{fpresred}$] We again proceed in three steps.

{\bf Step 1 (separatedness)}: Let us first check that once the existence of $i$ is proved then
we can find such an $\mathscr{X}$ that is $S$-separated if
$X$ is $S$-separated.   Make an initial choice of $\mathscr{X}$, and let
$\mathscr{I}$ denote the quasi-coherent ideal of $X$ in $\mathscr{X}$.
By \cite[I,~5.7.8]{rg}, $\mathscr{I} = \varinjlim \mathscr{I}_{\alpha}$
for the directed system of finite type quasi-coherent ideal sheaves
$\mathscr{I}_{\alpha}$ in $\mathscr{I}$.  Let $\mathscr{X}_{\alpha} \subseteq \mathscr{X}$
be the closed subspace cut out by $\mathscr{I}_{\alpha}$, so each
$\mathscr{X}_{\alpha}$ is finitely presented over $S$ and
$\invlim \mathscr{X}_{\alpha} = X$ is $S$-separated.  It suffices to prove that 
$\mathscr{X}_{\alpha}$ is $S$-separated for sufficiently large $\alpha$.
We may work \'etale-locally on $S$, so we may assume $S$ is affine. 

Since $\mathscr{X} \rightarrow S$ is finitely presented, by absolute
noetherian approximation for affine schemes (such as $S$) and Proposition \ref{lemma7} 
there is a cartesian diagram of algebraic spaces 
$$\xymatrix{
\mathscr{X} \ar[r] \ar[d] & S \ar[d] \\
\mathscr{X}_0 \ar[r] & \mathscr{S}}$$
with $\mathscr{S}$ affine of finite type over $\Z$ and the bottom side of finite presentation.
The normalization of $\mathscr{X}_{0,{\rm{red}}}$ is therefore
a finite surjective map of algebraic spaces $\mathscr{X}_1 \twoheadrightarrow \mathscr{X}_0$
(so $\mathscr{X}_1$ is finitely presented over $\Z$), 
and by Proposition \ref{prop1} (which rests on \S\ref{repgeneral} and \S\ref{topsec})
there is a finite surjective map
$\mathscr{X}'_1 \rightarrow \mathscr{X}_1$ with $\mathscr{X}'_1$ a scheme.
This finite map is finitely presented, since its target
is finitely presented over $\Z$.  
Applying the affine base change by $S \rightarrow \mathscr{S}$, we get a finite and finitely presented
surjective map $\pi:\mathscr{X}' \twoheadrightarrow \mathscr{X}$ from a scheme. 

Define $\mathscr{X}'_{\alpha} = \pi^{-1}(\mathscr{X}_{\alpha})$,
so $\invlim \mathscr{X}'_{\alpha} = \pi^{-1}(X)$ is separated.  
Each $\mathscr{X}'_{\alpha}$ is finitely presented over $S$
since $\mathscr{X}_{\alpha}$ is finitely presented over $S$ and $\pi$ is finitely presented.  
Since $\{\mathscr{X}'_{\alpha}\}$ is an inverse system of finitely presented $S$-schemes
with affine transition maps, by applying \cite[Thm.~C.7]{tt} over the affine $S$ 
it follows that $\mathscr{X}'_{\alpha}$ is separated for sufficiently large $\alpha$.
Finiteness and surjectivity of $\pi$ then gives that the 
algebraic space $\mathscr{X}_{\alpha}$ of finite presentation over the affine $S$ is separated for
sufficiently large $\alpha$, as desired. 

We may now ignore separatedness considerations and focus on merely constructing
a closed immersion $i:X \hookrightarrow \mathscr{X}$ over $S$ into an algebraic
space $\mathscr{X}$ of finite presentation over $S$.

{\bf Step 2 (proof when $S$ is a scheme, by induction on $X$-stratification)}:

Assume $S$ is a scheme.  Applying Theorem \ref{rgstrat} to the qcqs algebraic space $X$ yields 
a rising union of quasi-compact open subspaces
$$\emptyset = X_0 \subseteq X_1 \subseteq \dots \subseteq X_{s-1} \subseteq X_s = X$$
whose (locally closed) strata are schemes (with their reduced structure).  
Let us now proceed by induction on $s$, the case $s = 0$ being trivial.  
Thus, we may assume $s \ge 1$ and that the result is known for
$U = X_{s-1}$ (a quasi-compact open subspace of $X$).

Let $X' \rightarrow X$ be an \'etale covering
by an affine scheme (so $X' \rightarrow X$ is quasi-compact and separated as well, hence
finitely presented).  
The pullback $U' = U \times_X X'$ is a quasi-compact open subspace of the scheme $X'$ 
(so it is a scheme) and 
the projection $q:U' \rightarrow U$ is finitely presented, separated, and \'etale.
The inductive hypothesis provides a closed 
immersion $U \hookrightarrow \mathscr{U}$ into a finitely presented algebraic
space over the scheme $S$.  Moreover, since $X' \rightarrow S$ is a finite type  
map of schemes, by the known scheme case \cite[Thm.~4.3]{conrad} there is a closed
immersion $i:X' \hookrightarrow \mathscr{X}'$ into a scheme $\mathscr{X}'$ of finite
presentation over $S$.   The quasi-compact open subscheme $U' \subseteq X'$
has the form $X' \cap \mathscr{U}'$ for a quasi-compact open subscheme
$\mathscr{U}' \subseteq \mathscr{X}'$.  Applying Proposition \ref{gluingyourmama} to the resulting diagram of type \eqref{pushdiagram}, we  can use 
Lemma \ref{fpresredpre1} to find a closed immersion $X\hookrightarrow\mathscr X$
over $S$ with $\mathscr X$ of finite presentation over $S$, as desired.

{\bf Step 3 (proof for general $S$, by induction on $S$-stratification)}: 
Consider a rising union 
$$\emptyset  = U_0 \subseteq U_1 \subseteq \dots \subseteq U_r = S$$ 
of quasi-compact open subspaces as in Theorem \ref{rgstrat}.
Let $X_i = X \times_S {U_i}$.  We shall proceed by induction on $r$.  
The case $r = 0$ is trivial, so we can assume $r \ge 1$
and that $U := X_{r-1}$ has a closed
immersion over $U_{r-1}$ into a finitely presented algebraic space $\mathscr{U}$ over $U_{r-1}$
(which in turn is finitely presented over $S$ since $U_{r-1}$ is a quasi-compact open subspace
of the qcqs $S$).  

Let $\varphi:Y \rightarrow U_r = S$ be an \'etale covering by a quasi-compact separated (e.g., affine)
scheme, so $\varphi$ is a quasi-compact
separated \'etale map (hence of finite presentation)
and $X' := X \times_S Y$ is an algebraic space that is separated
and of finite type over $Y$.  Since $U_{r-1}$ is a quasi-compact open subspace
of $S$, the algebraic space
 $X'$ contains $U' = X' \times_S U_{r-1}$ as a quasi-compact open subspace  
such that the natural map $q:U' \rightarrow X \times_S U_{r-1} = X_{r-1} =: U$
is quasi-compact, \'etale, and separated.  But we noted that by induction on $r$ 
there is a closed immersion $U \hookrightarrow \mathscr{U}$ over 
$U_{r-1}$ into a finitely presented algebraic space over $U_{r-1}$, and that 
$\mathscr{U}$ is also finitely presented over $S$.  

By the settled case when the base is a scheme (Step 2), the separated finite type
morphism $X' \rightarrow Y$ to the qcqs scheme $Y$
admits a closed immersion $i:X' \hookrightarrow \mathscr{X}'$
over $Y$ into an algebraic space of finite presentation over $Y$.  By composing 
$\mathscr{X}' \rightarrow Y$ with
the finitely presented 
$\varphi:Y \rightarrow S$, $i$ may also be viewed
as a closed immersion over $S$ into an algebraic space
that is separated and of finite presentation over $S$. 
Since $X'$ is a closed subspace of $\mathscr{X}'$ (via $i$),
the quasi-compact open subspace $U' \subseteq X'$ has the form
$X' \cap \mathscr{U}'$ for a quasi-compact open subspace $\mathscr{U}' \subseteq \mathscr{X}'$.  

These spaces fit together into a diagram of the form \eqref{pushdiagram}, but possibly without
a map $\pi: \mathscr{U}' \rightarrow \mathscr{U}$ as in 
Lemma \ref{fpresredpre1}.  But applying Proposition \ref{gluingyourmama} puts us in the case when
the map $\pi$ does exist, so we may then apply Lemma \ref{fpresredpre1} to complete the proof
of Theorem \ref{fpresred}.  
\end{proof}

To conclude this section, let us use Theorem \ref{fpresred} 
to prove Theorem \ref{mainresult} in general; this includes the general noetherian case,
which has not yet been proved without excellence hypotheses.  Note that in 
\S\ref{exccase} we already settled the case when 
$S$ is of finite presentation over $\Z$.   By Theorem \ref{fpresred}, there is 
a closed immersion $i: X \hookrightarrow X'$ (over $S$)
into an algebraic space $X'$ that is separated and of finite presentation over
$S$.  If $X'$ admits an $S$-compactification $j:X' \hookrightarrow
\overline{X}'$ then the composite immersion
$j \circ i:X \rightarrow \overline{X}'$ over $S$ is quasi-compact, so its
scheme-theoretic image $\overline{X} \subseteq \overline{X}'$ is
an $S$-compactification of $X$.  Hence, it suffices
to treat $X'$ instead of $X$, so we may and do now assume that $X$ is finitely
presented (and separated) over $S$.  

By Theorem \ref{absapprox} (which was proved in
\S\ref{absnoetherian}) there is an isomorphism
$S \simeq \invlim S_{\alpha}$
with $\{S_{\alpha}\}$ an inverse system with affine transition maps between
algebraic spaces of finite presentation over $\Z$.   Choose an \'etale cover $U \rightarrow X$
with $U$ an affine scheme, so $U$ and $R = U \times_X U$ are schemes
of finite presentation over $S$.  By Proposition \ref{lemma7}, 
for sufficiently large $\alpha_0$ there are
 schemes $R_{\alpha_0}$ and $U_{\alpha_0}$
of finite presentation over $S_{\alpha_0}$ such that for
$U_{\alpha} := S_{\alpha} \times_{S_{\alpha_0}} U_{\alpha_0}$
and $R_{\alpha} := S_{\alpha} \times_{S_{\alpha_0}} R_{\alpha_0}$
we have $U = \invlim U_{\alpha}$ and $R = \invlim R_{\alpha}$.   

By Corollary \ref{mapdescend1} and Corollary \ref{mapdescend2}, for sufficiently large $\alpha$
the closed immersion  $R \rightarrow U \times_S U$
over $S$ arises by base change from a closed immersion 
$R_{\alpha} \rightarrow U_{\alpha} \times_{S_{\alpha}} U_{\alpha}$
over $S_{\alpha}$ such that the projections $R_{\alpha} \rightrightarrows
U_{\alpha}$ are \'etale.  A further application of Corollary \ref{mapdescend1} 
gives that $R_{\alpha}$ is thereby an equivalence relation on $U_{\alpha}$
for sufficiently large $\alpha$.  Thus, $X_{\alpha} := U_{\alpha}/R_{\alpha}$
is a separated algebraic space of finite presentation over $S_{\alpha}$
and it induces $X \rightarrow S$ after base change along $S \rightarrow S_{\alpha}$.
Theorem \ref{mainresult} is already known for 
$X_{\alpha} \rightarrow S_{\alpha}$ since $S_{\alpha}$ is
of finite presentation over $\Z$, so the theorem is now fully proved in general.

\begin{appendix}

\section{Some foundational facts for algebraic spaces}

In \S\ref{notation} we explained why the avoidance of quasi-separatedness hypotheses
in the foundations of the theory of algebraic spaces is convenient when carrying out
some of the gluing constructions that we need in the proof of Nagata's theorem for qcqs algebraic spaces.
We also need to use the formalism of inverse limits (with affine transition maps) when the base object
$S = \invlim S_{\alpha}$ is an inverse limit of
algebraic spaces $S_{\alpha}$ rather than of schemes, and so we have to prove
certain results in this generality that are analogues of known results in the case of limits of schemes.
It is convenient to not have to keep track of quasi-separatedness {\em a priori}
in such limit constructions as well. 

This appendix collects the results that we need from the theory of algebraic spaces and for which
we do not know of a reference in the literature.  We also provide some instructive examples 
of the surprises that can occur in the absence of quasi-separatedness
(but the reader need not get worried about such things; for our purposes, the reason
to avoid quasi-separatedness in the foundations is so that we can make sense of certain 
constructions {\em prior} to checking that the construction is quasi-separated in cases that we require).

\subsection{Algebraic spaces without quasi-separatedness}\label{repgeneral}

In \cite[I,~5.7.2]{rg} it is proved that if $R \rightrightarrows U$ is an \'etale equivalence relation
in schemes and $X = U/R$ is  the quotient sheaf on the \'etale site of the category of schemes
then the quotient map of sheaves $U \rightarrow X$ is representable by \'etale surjections
of schemes.  That is, for any map $W \rightarrow X$ from a scheme $W$,
$U \times_X W$ is a scheme.  Note that in this fiber product, the map $U \rightarrow X$ is rather special:
it arises from an \'etale presentation of $X$.

To make the theory of algebraic spaces workable without quasi-separatedness hypotheses
(with algebraic spaces defined to be \'etale sheaves admitting such a quotient presentation
$U/R$), the essential step is to prove that if $V$ is any scheme then 
any map $V \rightarrow X$ is representable; that is, $V \times_X W$ is a scheme for
any pair of maps $V \rightarrow X$ and $W \rightarrow X$ from schemes $V$ and $W$. 
Since $V \times_X W = (V \times W) \times_{X \times X} X$ (with absolute products
over $\Spec \Z$), this representability is a consequence of the following result.

\begin{proposition}\label{schemprod}
 For $X$ as above, the diagonal map $\Delta_X:X \rightarrow X \times_{\Spec \Z} X$
is representable.  More generally, if $X \rightarrow S$ is any map of algebraic spaces then 
$\Delta_{X/S}:X \rightarrow X \times_S X$ is representable.
\end{proposition}

Recall that a map $X \rightarrow Y$ of contravariant functors on the category of schemes is {\em representable} if
$X \times_Y T$ is a scheme for any scheme $T$. 

\begin{proof}
If $V$ is a scheme and $h:V \rightarrow X \times_S X$ is a map corresponding to a pair
of maps $h_1,h_2:V \rightrightarrows X$ then $h$ is the composition 
$V \rightarrow V \times_S V \rightarrow X \times_S X$, where the second step is $h_1 \times h_2$.
Provided that $\Delta_S$ is representable, so $V \times_S V$ is a scheme, the representability
of $\Delta_{X/S}^{-1}(V)$  reduces to the representability of
$\Delta_{X/S}^{-1}(V \times_S V) = \Delta_X^{-1}(V \times_{\Spec \Z} V)$. 
Hence, the representability of $\Delta_{X/S}$ in general reduces to the special
case $S = \Spec \Z$.   Moreover, as we just saw, to handle this case
it suffices to prove representability of pullbacks $\Delta_X^{-1}(V \times W) = V \times_X W$
for any pair of maps $v:V \rightarrow X$ and $w:W \rightarrow X$
from schemes $V$ and $W$.  By working
Zariski-locally on $W$ and then on $V$ we may assume each is affine. 
Now we shall adapt the proof of \cite[I,~5.7.2]{rg} (which is the case when one of $V$ or $W$ is
the cover of $X$ in an \'etale quotient presentation with a representable sheaf of relations). 

Fix an \'etale quotient presentation in schemes $R \rightrightarrows U$ for $X$, as in the
definition of $X$ as an algebraic space.  Since $U \rightarrow X$ is a surjection
for the \'etale topology 
and $V$ is an affine scheme, there is an \'etale affine cover $V' \rightarrow V$ 
fitting into a commutative diagram 
$$\xymatrix{
{V'} \ar[d] \ar[r] & {U} \ar[d] \\
{V} \ar[r]_-{v} & {X}}$$
Since $U \times_X W$ is a scheme (by \cite[I,~5.7.2]{rg}), $V' \times_X W = V' \times_U (U \times_X W)$
is also a scheme.    If we let 
$V'' = V' \times_V V'$ then there is a canonical map
$V'' \rightarrow X$ factoring through either of the two evident maps
$V'' \rightrightarrows U$, and using either of these shows similarly that
$V'' \times_X W$ is a scheme.  The fiber square of the sheaf quotient
map $V' \times_X W \rightarrow V \times_X W$
is identified with the \'etale equivalence relation in schemes
\begin{equation}\label{vvw}
V'' \times_X W \rightrightarrows V' \times_X W
\end{equation}
which can be interpreted as a descent datum on the $V'$-scheme $V' \times_X W$
relative to the \'etale covering map of affine schemes $V' \rightarrow V$.  Hence,
our problem is to prove effectivity of this descent problem in
the category of $V$-schemes.

The map of schemes $V' \times_X W \rightarrow V' \times_{\Spec \Z} W$ with affine target
is separated (even a monomorphism), so $V' \times_X W$ is a separated scheme.  Since the
maps in (\ref{vvw}) are quasi-compact \'etale surjections, it follows from separatedness of
$V' \times_X W$ that we can cover $V' \times_X W$ by quasi-compact open
subschemes which are stable under the descent datum over $V' \rightarrow V$.  
Now we can argue as in the proof of \cite[I,~5.7.2]{rg}:  since \'etale descent is
effective for schemes that are quasi-affine over the base \cite[VIII,~7.9]{sga1}, 
it suffices to show that every quasi-compact open subscheme of $V' \times_X W$ is
quasi-affine.   Any quasi-finite separated map of schemes is quasi-affine
\cite[IV$_4$,~18.12.12]{ega}, so it suffices to show that 
the separated scheme $V' \times_X W$ admits a locally quasi-finite morphism to an affine scheme.
 The map of schemes 
$V' \times_X W = V' \times_U (U \times_X W) \rightarrow V' \times_U W$
is \'etale since $U \times_X W \rightarrow W$ is \'etale (by \cite[I,~5.7.2]{rg}),
and the map $V' \times_U W \rightarrow V' \times_{\Spec \Z} W$ with affine target 
is an immersion since
it is a base change of $\Delta_U$.  
\end{proof}

\begin{corollary}\label{qtspace}
 If $R \rightrightarrows U$ is an \'etale equivalence relation in algebraic spaces
then the quotient $X = U/R$ is an algebraic space.
\end{corollary}

\begin{proof}
Let $U' \rightarrow U$ be an \'etale cover by a separated scheme, so $X = U'/R'$ where
$R' = U' \times_X U' = U' \times_{U,p_1} R \times_{p_2,U} U'$ is clearly an algebraic space.
The natural maps of algebraic spaces $R' \rightrightarrows U'$ are easily checked to be \'etale, so 
we just need to check that $R'$ is a scheme.  Since 
$R' \rightarrow U' \times_{\Spec \Z} U'$ is a monomorphism, hence separated, 
$R'$ is a separated algebraic space (since we chose $U'$ to be separated).
Thus, $p'_1:R' \rightarrow U'$ is a separated \'etale map from an algebraic space to a 
(separated) scheme. 
It then follows from \cite[II,~Cor.~6.16]{knutson} that $R'$ is a scheme. 
\end{proof}

The ``topological invariance'' of the \'etale site of a scheme holds more generally for algebraic
spaces:

\begin{proposition}\label{redequiv} Let $X$ be an algebraic space.  The functor
$E \rightsquigarrow E_{\rm{red}} = E \times_X X_{\rm{red}}$
from the category of algebraic spaces \'etale over $X$ to the category of 
algebraic spaces \'etale over $X_{\rm{red}}$ is an equivalence. 

In particular, the pullback and pushforward functors
between \'etale topoi ${\rm{\acute{E}t}}(X)$ and ${\rm{\acute{E}t}}(X_{\rm{red}})$ are naturally inverse to each other.
\end{proposition}

\begin{proof}  Let 
$U \rightarrow X$ be an \'etale cover by a scheme and let $R = U \times_X U$, so
$R$ is also a scheme.  Since $U \rightarrow X$ is \'etale, 
$U_{\rm{red}} = U \times_X X_{\rm{red}}$ and 
$R_{\rm{red}} = U_{\rm{red}} \times_{X_{\rm{red}}} U_{\rm{red}}$.
Hence, $R_{\rm{red}} \rightrightarrows U_{\rm{red}}$ is an \'etale chart for $X_{\rm{red}}$.  

For any algebraic space $E$ equipped with a map $E \rightarrow X$, we may identify
$E$ with the quotient of the \'etale equivalence relation in algebraic spaces
$R \times_X E \rightrightarrows U \times_X E$, and in case $E \rightarrow X$ is
\'etale we see that $E_{\rm{red}}$ is the quotient of the analogous equivalence relation
$$R_{\rm{red}} \times_{X_{\rm{red}}} E_{\rm{red}} \rightrightarrows U_{\rm{red}}
\times_{X_{\rm{red}}} E_{\rm{red}}.$$
By Corollary \ref{qtspace}, the category of algebraic spaces is stable under the formation of
quotients by arbitrary \'etale equivalence relations.  
Thus, since the maps $R \rightrightarrows U$ are \'etale
and pull back to the maps $R_{\rm{red}} \rightrightarrows U_{\rm{red}}$, we easily
reduce to proving the lemma for $U$ and $R$ in place of $X$.  That is, we
can assume that $X$ is a scheme.

The equivalence result for categories of {\em schemes} \'etale over $X$ and $X_{\rm{red}}$
is \cite[IV$_4$,~18.1.2]{ega}.  Thus, if $E_0$ is an algebraic space
\'etale over $X_{\rm{red}}$ and $U_0 \rightarrow E_0$ is an \'etale scheme cover
then $R_0 = U_0 \times_{E_0} U_0$ is a scheme and $R_0 \rightrightarrows U_0$
is an equivalence relation in \'etale $X_{\rm{red}}$-schemes with
quotient $E_0$.  This uniquely lifts to an equivalence relation
$R \rightrightarrows U$ in \'etale $X$-schemes, and it is clearly an
\'etale equivalence relation.  Hence, the quotient $E = U/R$ makes sense as an algebraic
space \'etale over $X$, and $E \times_X X_{\rm{red}} \simeq U_0/R_0 = E_0$.
Thus, we have proved essential surjectivity.  

For full faithfulness, 
let $E$ and $E'$ be algebraic spaces \'etale over $X$.  To show that the map 
$\Hom_X(E,E') \rightarrow \Hom_{X_{\rm{red}}}(E_{\rm{red}}, E'_{\rm{red}})$ is
bijective, by using an \'etale scheme chart for $E$ we may easily reduce to the case when
$E$ is a scheme.  That is, it suffices to prove that 
$E'(V) \rightarrow E'_{\rm{red}}(V_{\rm{red}})$ is bijective for any
\'etale $X$-scheme $V$ and any algebraic space $E'$ over $X$.  
Let $j:X_{\rm{red}} \rightarrow X$ denote the canonical map.  The functors
$j_{\ast}$ and $j^{\ast}$ between small \'etale sites are inverse equivalences of categories
(essentially due to the equivalence result for schemes \cite[IV$_4$,~18.1.2]{ega} that we already used), 
and if $V \rightarrow X$ is an \'etale map of schemes
we naturally have 
\begin{equation}\label{jv}
j^{\ast}(V) \simeq V \times_X X_{\rm{red}} = V_{\rm{red}}
\end{equation}
as sheaves on the small \'etale site of $X_{\rm{red}}$.  Thus, 
it suffices to show that the natural map $E'_{\rm{red}} \rightarrow j^{\ast}(E')$
of sheaves on the small \'etale site of $X_{\rm{red}}$ is an isomorphism.
But $j^{\ast}$ is an exact functor, so by using an \'etale scheme chart
for $E'$ we reduce to the case when $E'$ is a scheme, in which case  the map of 
interest is identified with the map in (\ref{jv}) that is an isomorphism. 
\end{proof}

\subsection{Topology of algebraic spaces}\label{topsec}

We briefly review a few facts about the topology of quasi-separated
algebraic spaces.  Rydh has recently developed a theory
of points and topological spaces associated to general algebraic
spaces.

\begin{definition}\label{noethdef}
 Let $S$ be an algebraic space.  It is {\em noetherian} if it is quasi-separated
and admits an \'etale cover by a noetherian scheme.  It is {\em locally noetherian}
if it admits a  Zariski-open covering by noetherian open subspaces.
\end{definition}

It is clear that any locally noetherian algebraic space is quasi-separated. 
In \cite[II,~8.5]{knutson} the 
irreducible component decomposition of noetherian 
algebraic spaces $S$ is developed, assuming $S$ is (quasi-separated and) locally separated
in the sense that  the map $\Delta_S$ is a quasi-compact immersion. 
The theory of irreducible components is developed in \cite{lmb}
for quasi-separated Artin stacks.  (Quasi-separatedness is 
a running assumption throughout \cite{lmb}, along with
separatedness of the diagonal.)  More specifically, specializing
 \cite[5.7.2]{lmb} to the case of algebraic spaces, one gets:
 
\begin{lemma}\label{genpt}  Let $S$ be a quasi-separated algebraic space.  Every irreducible closed
set in $|S|$ has a unique generic point,
and if $f:S' \rightarrow S$ is an \'etale surjective map of quasi-separated algebraic spaces then every
irreducible component of $S$ has generic point of the form $f(s')$ for some generic
point $s'$ of $S'$. 
\end{lemma}

We do not know if Lemma \ref{genpt} holds without quasi-separatedness.
% ; Example \ref{galex}
% and Example \ref{nongenex} illustrate the subtleties of $|S|$ when $S$ is not quasi-separated.

The following proposition is known more generally for quasi-separarted Artin stacks, but
we include a proof for the convenience of the reader. 

\begin{proposition}\label{irredcomp} Let $S$ be a quasi-separated algebraic space.  Its irreducible
components constitute 
the unique collection $\{S_i\}$ of pairwise distinct irreducible and reduced
closed subspaces of $S$ such that 
they cover $S$ and no $|S_i|$ is contained in  $\cup_{j \ne i} |S_j|$.
If $V$ is open in $S$ then the corresponding collection for $V$ is
the collection of those $V \cap S_i$ that are non-empty.  

If $S$ is locally noetherian then $\{S_i\}$ is Zariski-locally finite
and the normalization of $S_{\rm{red}}$ is the disjoint union
$\coprod \widetilde{S}_i$ where $\widetilde{S}_i$ is
the normalization of $S_i$ and each $\widetilde{S}_i$ is irreducible.
In particular, if $S$ is normal and locally noetherian then $\{S_i\}$ is its set of connected components.
\end{proposition}

\begin{proof} 
To show that an irreducible component of $S$ is not covered (as a topological space)
by the union of the other irreducible components of $S$, it suffices to show that if $S$ is irreducible
then $|S|$ is not covered by proper closed subsets.  The existence of a generic point 
as in Lemma \ref{genpt} makes this clear.
It is likewise clear via generic points that if $V$ is open in $S$ and $\{S_i\}$ is the set of
irreducible components of $S$ then the non-empty $V \cap S_i$'s are the irreducible components
of $V$ (as topological closure preserves irreducibility).  

For the uniqueness claim concerning
the collection $\{S_i\}$, consider any collection $\{S'_i\}$ of irreducible and reduced closed subspaces
of $S$ such that the $|S'_i|$'s cover $|S|$ and each $|S'_i|$ is not contained in
the union of the $|S'_j|$'s for $j \ne i$.  Each $S'_i$ must be an irreducible component, since
if $S'_i$ is strictly contained in another irreducible and reduced closed subspace
$Z$ then the generic point of $Z$ cannot lie in $S'_i$ and so lies in $S'_j$ for some $j \ne i$,
a contradiction since then $S'_i \subseteq Z \subseteq S'_j$ by closedness of each $S'_j$ in $S'$. 
Moreover, the generic point of each irreducible component $Y$ in $S$ lies in some
$S'_i$ and hence $Y = S'_i$ by maximality for $Y$.  This completes the proof of uniqueness. 

In the locally noetherian case the topological space
$|S|$ is locally noetherian, so $\{S_i\}$ is Zariski-locally finite in $S$.  
Consider the normalization $\widetilde{S}$ of $S_{\rm{red}}$.
The description of its irreducible components is known in the scheme case, and so
follows in general by working over a dense open subspace $U \subseteq S$
that is a scheme \cite[II,~6.8]{knutson} (with the normalization $\widetilde{U}$ likewise
a dense open subspace of $\widetilde{S}$).    The irreducible components of
$\widetilde{U}$ are its connected components, so to check that the 
irreducible components of $\widetilde{S}$ are pairwise
disjoint it suffices to show that idemponents on $\widetilde{U}$ uniquely
extend to idempotents on $\widetilde{S}$.  More generally, if $N$ is a normal algebraic
space and if $V$ is a dense open subspace, then idempotents of $V$ uniquely extend
to idempotents of $N$.  Indeed, to prove this it suffices to work \'etale-locally on $N$, so
we reduce to the known case when $N$ is a normal scheme.
\end{proof}

\subsection{Properties of inverse limits of algebraic spaces}\label{limprop}

The results in this section extend to the case of algebraic spaces a variety of known facts concerning
inverse limits of schemes.

\begin{proposition}\label{fpres}
Let $X \rightarrow Y$ be a map locally 
of finite presentation between algebraic spaces such that
$X$ and $Y$ are covered by Zariski-open subspaces that are quasi-separated.
If $\{S_{\alpha}\}$ is an inverse system of algebraic spaces  over $Y$ with affine
transition maps and each $S_{\alpha}$ is qcqs then 
$$\varinjlim \Hom_Y(S_{\alpha},X) \rightarrow \Hom_Y(\invlim S_{\alpha}, X)$$
is bijective.
\end{proposition}

In the case of schemes the Zariski covering hypotheses on $X$ and $Y$ are always
satisfied, so when all objects (including
the $S_{\alpha}$ are schemes this proposition is part of \cite[IV$_3$,~8.14.2(b)]{ega}. 

\begin{proof}
The first part of the proof consists of some tedious but mechanical arguments to reduce
to the case when $X$ is qcqs and $Y$ and all $S_{\alpha}$ are affine.  Then the real argument can begin. 

Since any particular $S_{\alpha_0}$ is quasi-compact, so a Zariski-open covering of $Y$
has pullback to $S_{\alpha_0}$ admitting a finite quasi-compact open subcover,
it is straightforward to reduce to general problem over $Y$
to the general problem over each member of a Zariski-open covering of $Y$.
Hence, we can assume that $Y$ is quasi-separated, and then repeating the same
argument reduces us to the case when $Y$ is qcqs.  Hence, 
$Y$ admits an \'etale cover $Y' \rightarrow Y$ by an affine scheme $Y'$,
and this covering map is quasi-compact (since $Y$ is qcqs) and separated.
An elementary \'etale descent argument thereby reduces our problem to the case when $Y$ is affine.

Let $U_{\alpha_0} \rightarrow S_{\alpha_0}$ be an \'etale cover by an
affine scheme for some $\alpha_0$ 
and let $R_{\alpha_0} = U_{\alpha_0} \times_{S_{\alpha_0}} U_{\alpha_0}$, 
so $R_{\alpha_0}$ is a quasi-compact and separated scheme. 
We may and do assume $\alpha \ge \alpha_0$ for all $\alpha$.  
Define $R_{\alpha} \rightrightarrows U_{\alpha}$ and $R \rightrightarrows U$
by base change along $S_{\alpha} \rightarrow S_{\alpha_0}$ and $S := \invlim
S_{\alpha} \rightarrow S_{\alpha_0}$,
so these are respective \'etale charts for $S_{\alpha}$ and $S$.  It is easy to reduce
to the analogous bijectivity claims for the inverse systems $\{U_{\alpha}\}$ and $\{R_{\alpha}\}$.
In particular, we may assume that all $S_{\alpha}$ are schemes, and then even quasi-compact
and separated (and even affine).

With all $S_{\alpha}$ now quasi-compact
separated schemes, any quasi-compact open subscheme
of $\invlim S_{\alpha}$ descends to a quasi-compact open subscheme of
some $S_{\alpha_0}$.   Thus, if $\{X_i\}$ is any Zariski-open covering of $X$ and 
we can solve the problem for each $X_i$ separately (for any inverse system
$\{S_{\alpha}\}$ of quasi-compact and separated schemes 
equipped with a compatible map to $Y$) then it is straightforward to use 
quasi-compact open refinements
of pullbacks of the Zariski-covering $\{X_i\}$ get the result for
$X$.   But we assumed that $X$ has a Zariski-open covering by quasi-separated
algebraic spaces, so we may arrange that $X$ is qcqs.  Hence, $X \rightarrow Y$ is finitely presented.
We may also now run through the localization argument on the $S_{\alpha}$'s to reduce to the case
when all $S_{\alpha}$ are affine.  Keep in mind that we also reduced to the case when $Y$ is affine. 

Consider a $Y$-map $f:S \rightarrow X$.  Let $V \rightarrow X$ be an \'etale cover by an affine scheme, 
so this covering map is quasi-compact (since $X$ is quasi-separated) and separated
(since $V$ is separated).  Hence, $V$ is finitely presented over $Y$
and  the \'etale scheme cover $W := f^{-1}(V) \rightarrow S$ is finitely presented
(i.e., qcqs). It follows that
$W \rightarrow S$ 
descends to a qcqs \'etale scheme cover $W_{\alpha_0} \rightarrow S_{\alpha_0}$
for some $\alpha_0$.  The pullbacks $W_{\alpha} = W_{\alpha_0} \times_{S_{\alpha_0}} S_{\alpha}$
form an inverse system of 
$Y$-schemes, with affine transition maps, so since $V$ is a scheme that is finitely presented over
$Y$ it follows that the natural map 
$$\varinjlim \Hom_Y(W_{\alpha}, V) \rightarrow \Hom_Y(W,V)$$
is bijective by \cite[IV$_3$,~8.14.2(b)]{ega}.
Thus, the natural $Y$-map $h:W = f^{-1}(V) \rightarrow V$
factors through a $Y$-map $h_{\alpha_1}:W_{\alpha_1} \rightarrow V$ for some large $\alpha_1$.

The pair of maps $h_{\alpha_1} \circ p_1, h_{\alpha_1} \circ p_2$
in
$\Hom_Y(W_{\alpha_1} \times_{S_{\alpha_1}} W_{\alpha_1}, V)$
may not coincide, but they do yield the same map
$h \circ p_1 = h \circ p_2$ 
after composing back to $W \times_S W = f^{-1}(V \times_X V)$.
Thus, since $\{W_{\alpha} \times_{S_{\alpha}} W_{\alpha}\}$
is an inverse system of qcqs schemes with affine
transition maps (over $Y$) and limit $W \times_S W$, by the injectivity of the map 
$$\varinjlim \Hom_Y(W_{\alpha} \times_{S_{\alpha}} W_{\alpha}, V)
\simeq \Hom_Y(W \times_S W, V)$$
it follows that if we enlarge $\alpha_1$ enough then
$h_{\alpha_1} \circ p_1 = h_{\alpha_1} \circ p_2$.  
That is, $h_{\alpha_1}$ factors through a $Y$-map
$f_{\alpha_1}:S_{\alpha_1} \rightarrow V$.  It suffices to check that this induces the original map
$f:S \rightarrow V$ after composing with $S \rightarrow S_{\alpha_1}$. 
Such an equality may be checked after composing with
the \'etale covering map $W \rightarrow S$, so we are done because the diagram
$$\xymatrix{
{W} \ar[r] \ar[d] & {S} \ar[d]\\
{W_{\alpha}} \ar[r] & {S_{\alpha}}}$$
commutes for all $\alpha$ (due to how $W_{\alpha}$ was defined). 
\end{proof}

\begin{corollary}\label{mapdescend1}
Let $\{S_{\alpha}\}$ be an inverse system of qcqs 
algebraic spaces with affine transition maps and limit $S$.  
Let $\{X_{\alpha}\}$ and $\{Y_{\alpha}\}$ be
compatible cartesian systems of finitely presented algebraic spaces over $\{S_{\alpha}\}$,
and let $X = \invlim X_{\alpha}$ and $Y = \invlim Y_{\alpha}$. 

The natural map
$$\varinjlim \Hom_{S_{\alpha}}(X_{\alpha}, Y_{\alpha}) \rightarrow
\Hom_S(X,Y)$$
is bijective. 
\end{corollary}

In Corollary \ref{mapdescend2} 
we will address the descent of properties of morphisms through such limits.
We postpone it because the verification of the properties 
requires an additional result given in Proposition \ref{lemma7} below.

\begin{proof}
Fix some $\alpha_0$ and consider only $\alpha \ge \alpha_0$ without loss of generality.  
Since $Y_{\alpha} = S_{\alpha} \times_{S_{\alpha_0}} Y_{\alpha_0}$ and 
$Y = S \times_{S_{\alpha_0}} Y_{\alpha_0}$, the map of interest may be expressed as
$$\varinjlim \Hom_{S_{\alpha_0}}(X_{\alpha}, Y_{\alpha_0}) \rightarrow
\Hom_{S_{\alpha_0}}(\invlim X_{\alpha}, Y_{\alpha_0}).$$
Thus, the desired bijectivity is a consequence of Proposition \ref{fpres}. 
\end{proof}

\begin{lemma}\label{scheme}
Let $\{S_{\alpha}\}$ be an inverse system of qcqs algebraic spaces with affine 
transition maps.  Let $S = \invlim S_{\alpha}$, and assume
that some $S_{\alpha'}$ is affine over an algebraic space
of finite presentation over $\Z$.  

If $S$ is a scheme $($resp. affine scheme$)$
then so is $S_{\alpha}$ for sufficiently large $\alpha$.
\end{lemma}

The hypothesis on some $S_{\alpha'}$ is satisfied by all $S_{\alpha}$
when they are schemes, by absolute noetherian approximation for qcqs schemes 
\cite[Thm.~C.9]{tt}.  The proof of absolute noetherian approximation for qcqs algebraic
spaces in \S\ref{absnoetherian} requires Lemma \ref{scheme} in the special case
that all $S_{\alpha}$ are of finite presentation over $\Z$, but once that
result is proved in general then the hypothesis on some $S_{\alpha'}$
in Lemma \ref{scheme} may be dropped (as it is then satisfied for any $\alpha'$). 

  \begin{proof}
Fix an $\alpha'$ as in the statement of the lemma,
 and consider only $\alpha \ge \alpha'$ without loss of generality.
 In particular, all $S_{\alpha}$ are affine over an algebraic space $\mathscr{S}$ of finite
 presentation over $\Z$.  When the $S_{\alpha}$'s are all schemes, by absolute
 noetherian approximation for qcqs schemes \cite[Thm.~C.9]{tt}
 we may choose $\mathscr{S}$ to be a scheme; note that $\mathscr{S}$ is an auxiliary
 device, so changing it at the outset is not a problem.
 (Even when the $S_{\alpha}$ are all schemes and $S$ is affine, our desired
 conclusion that $S_{\alpha}$ is affine for all large $\alpha$ is slightly more general than
 \cite[Prop.~C.6]{tt} in the sense that the argument there assumes the $S_{\alpha}$ are
 all finite type over $\Z$, or over some noetherian ring.)
  
  Fix some $\alpha_0$ and let 
$U_{\alpha_0} \rightarrow S_{\alpha_0}$ be an \'etale scheme cover with 
$U_{\alpha_0}$ a qcqs scheme (e.g., affine), 
and define $R_{\alpha_0} = U_{\alpha_0} \times_{S_{\alpha_0}} U_{\alpha_0}$,
so $R_{\alpha_0}$ is a qcqs scheme since $S_{\alpha_0}$ is quasi-separated and
$U_{\alpha_0}$ is quasi-compact.  Define $R_{\alpha} \rightrightarrows U_{\alpha}$
and $R \rightrightarrows U$ by base change to $S_{\alpha}$ and $S$ respectively,
so these are \'etale scheme charts for $S_{\alpha}$ and $S$ respectively. 
Let $p:U \rightarrow S$ be the covering map, and $p_1, p_2:R \rightrightarrows U$
be the two projections.

We claim (exactly as is known for inverse systems of qcqs {\em schemes}
with affine transition maps) that if $V \subseteq S$
is a quasi-compact open subspace then $V$ descends to a quasi-compact
open subspace $V_{\alpha_0} \subseteq S_{\alpha_0}$ for some $\alpha_0$, and that
if $V'_{\alpha_0} \subseteq S_{\alpha_0}$ is another such descent then the pullbacks of
$V_{\alpha_0}$ and $V'_{\alpha_0}$ to $S_{\alpha}$ coincide for sufficiently large $\alpha \ge \alpha_0$.

To prove these claims, consider the quasi-compact open subscheme $W := p^{-1}(V) \subseteq U$.
Since $\{U_{\alpha}\}$ is an inverse system of qcqs schemes with affine
transition maps and limit $U$, certainly $W$ descends to a quasi-compact open subscheme
$W_{\alpha_0} \subseteq U_{\alpha_0}$ for some $\alpha_0$.  The two open subschemes
$p_{1,\alpha_0}^{-1}(W_{\alpha_0})$ and $p_{2,\alpha_0}^{-1}(W_{\alpha_0})$ in
$R_{\alpha_0}$ have pullbacks $p_1^{-1}(W)$ and $p_2^{-1}(W)$ in $R = \invlim R_{\alpha}$
which coincide.  Since $\{R_{\alpha}\}$ is an inverse system of qcqs schemes
with affine transition maps, by increasing $\alpha_0$ we can ensure that
the two pullbacks of $W_{\alpha_0}$ into $R_{\alpha_0}$ coincide, so
$W_{\alpha_0}$ descends to a quasi-compact open subspace $V_{\alpha_0}
\subseteq S_{\alpha_0}$, and this clearly descends $V$.   If
$V'_{\alpha_0}$ is another such descent then the desired equality $V_{\alpha} = V'_{\alpha}$
of pullbacks into $S_{\alpha}$ for sufficiently large $\alpha$
can be checked by working in the inverse system $\{U_{\alpha}\}$ of
qcqs \'etale scheme covers of $\{S_{\alpha}\}$, so we are reduced to the known scheme case.

Choose a finite covering of the scheme $S$ by affine open subschemes
$\{V_1, \dots, V_n\}$.   By what we just proved, for sufficiently large $\alpha_0$ this descends
to a finite collection of quasi-compact open subspaces
$\{V_{1,\alpha_0}, \dots, V_{n, \alpha_0}\}$ in $S_{\alpha_0}$, and by taking $\alpha_0$ big enough
it constitutes a covering.   For each $1 \le i \le n$, the resulting inverse system
$\{V_{i,\alpha}\}_{\alpha \ge \alpha_0}$ of qcqs algebraic spaces
with affine transition maps has limit $V_i$ that is an affine scheme,
and $\{V_{1,\alpha}, \dots, V_{n,\alpha}\}$ is an open cover of
$S_{\alpha}$.  Thus, $S_{\alpha}$ is a scheme if each
$V_{i,\alpha}$ is a scheme, so for each fixed $i$ we 
may rename $V_i$ as $S$ and $V_{i,\alpha}$ as $S_{\alpha}$
to reduce to the case when $S$ is an affine scheme, provided of course that 
we prove the stronger conclusion that $S_{\alpha}$ is an affine scheme for sufficiently large $\alpha$.

Recall that $\{S_{\alpha}\}$ is an inverse system of objects that are affine
over an algebraic space $\mathscr{S}$ of finite presentation over $\Z$, and that when all
$S_{\alpha}$ are schemes we may take $\mathscr{S}$ to be a scheme.  
Since the map $S \rightarrow S_{\alpha}$ is affine, we can form its 
schematic image $S'_{\alpha} \subseteq S_{\alpha}$ for all $\alpha$, and
clearly $\{S'_{\alpha}\}$ is an inverse system with schematically dominant
and affine transition maps such that the natural map $S \rightarrow \invlim S'_{\alpha}$
is an isomorphism.  Each map $S'_{\alpha} \rightarrow \mathscr{S}$
is affine and so corresponds to a quasi-coherent sheaf of ${\ms O}_{\mathscr{S}}$-algebras
$\mathscr{A}'_{\alpha}$ with injective transition maps such
that $S$ corresponds to $\mathscr{A} = \varinjlim \mathscr{A}'_{\alpha}$.
By \cite[III,~Thm.~1.1,~Cor.~1.2]{knutson}, each $\mathscr{A}'_{\alpha}$ is the direct
limit of its coherent ${\ms O}_{\mathscr{S}}$-submodules.  Passing to the ${\ms O}_{\mathscr{S}}$-subalgebras
generated by these, we get $\mathscr{A}'_{\alpha} = \varinjlim_{\beta} \mathscr{B}_{\alpha,\beta}$
with $\mathscr{B}_{\alpha,\beta}$ ranging through the set of finite type
quasi-coherent ${\ms O}_{\mathscr{S}}$-subalgebras of $\mathscr{A}'_{\alpha}$. 
Let $\{\gamma\}$ denote the set of such pairs $(\alpha, \beta)$
and define $\gamma' \ge \gamma$ if $\mathscr{B}_{\gamma'}$ contains
$\mathscr{B}_{\gamma}$ as subsheaves of $\mathscr{A}$.  
Hence
$S = \invlim T_{\gamma}$ over $\mathscr{S}$ with $T_{\gamma} = \Spec_{\mathscr{S}}(\mathscr{B}_{\gamma})$. 

By Proposition \ref{fpres} applied to the finitely presented map 
$T_{\gamma} \rightarrow \mathscr{S}$ and the inverse system 
$\{S_{\alpha}\}$ over $\mathscr{S}$ with limit $S$, 
for each $\gamma_0$ the canonical $\mathscr{S}$-map $S \rightarrow T_{\gamma_0}$ factors
through an $\mathscr{S}$-map $S_{\alpha} \rightarrow T_{\gamma_0}$ for
some sufficiently large $\alpha$ (depending on $\gamma_0$).  
But $S_{\alpha}$ and $T_{\gamma_0}$ are each affine over $\mathscr{S}$, so 
any $\mathscr{S}$-map between them is automatically affine.  Thus, 
if $T_{\gamma_0}$ is affine for some $\gamma_0$  then $S_{\alpha}$ is affine
for all large $\alpha$, as desired.  It is therefore permissible to replace
$\{S_{\alpha}\}$ with $\{T_{\gamma}\}$, so we thereby reduce to the case
when all $S_{\alpha}$ are of finite presentation over $\Z$.  This reduction step preserves
the property that all $S_{\alpha}$ are schemes if we began in that case
and took $\mathscr{S}$ to be a scheme (as we may when all
$S_{\alpha}$ are schemes).  
Now that all $S_{\alpha}$ are of finite presentation over $\Z$, if 
they are all schemes then $S$ being affine forces
$S_{\alpha}$ to be affine for all large $\alpha$ by \cite[Prop.~C.6]{tt}.
Thus, in the original problem, if the $S_{\alpha}$ are schemes
and $S$ is affine then we have deduced the desired
stronger conclusion that the $S_{\alpha}$ are
affine for all large $\alpha$.  (In other words,  we have slightly generalized \cite[Prop.~C.6]{tt} by
eliminating the finite presentation hypotheses there.)

Consider the original general problem, with
$S$ a scheme and each $S_{\alpha}$ just an algebraic space.
To conclude that the $S_{\alpha}$ are schemes for large
$\alpha$, we have already reduced ourselves to the case when $S$ is an affine
scheme and the $S_{\alpha}$ are algebraic spaces 
of finite presentation over $\Z$.  In this case
we want to conclude that once again the $S_{\alpha}$ are affine
for all large $\alpha$.  Since Proposition \ref{fpres} is available in the category of algebraic
spaces, we may use it (with $Y = \Spec \Z$) in place of \cite[C.5]{tt} to make the proof of
the scheme case in \cite[Prop.~C.6]{tt} work verbatim with each $S_{\alpha}$ an 
algebraic space of finite presentation over $\Z$. 
\end{proof}

The following result is useful for both descent through limits of algebraic spaces
and  relating such descents to limit procedures with schemes.  

\begin{proposition}\label{lemma7}
Let $\{S_{\alpha}\}$ be an inverse system of qcqs algebraic spaces
with affine transition maps and limit $S$.  Assume that some
$S_{\alpha'}$ is affine over an algebraic space of finite presentation over $\Z$. 

For any finitely presented morphism
$T \rightarrow S$ from an algebraic space
$($resp. scheme$)$ $T$, there exists an $\alpha_0$ and 
a finitely presented map $T_{\alpha_0} \rightarrow S_{\alpha_0}$ from an
algebraic space $($resp. scheme$)$ 
$T_{\alpha_0}$ such that $T \simeq T_{\alpha_0} \times_{S_{\alpha_0}} S$
over $S$, so the natural map 
$$T \rightarrow \invlim_{\alpha \ge \alpha_0} T_{\alpha}$$
is an isomorphism with $T_{\alpha} := T_{\alpha_0} \times_{S_{\alpha_0}} S_{\alpha}$.
\end{proposition}

The hypothesis that some $S_{\alpha'}$ is affine over an algebraic space
of finite presentation over $\Z$ is always satisfied when the
$S_{\alpha}$ are schemes, by absolute noetherian approximation for
qcqs schemes \cite[Thm.~C.9]{tt}.  This hypothesis will be satisfied in general once absolute
noetherian approximation for qcqs algebraic spaces is proved in \S\ref{absnoetherian}.
The proof of this latter approximation result requires Proposition \ref{lemma7} only in the special
case that all $S_{\alpha}$ are of finite presentation over $\Z$, so
there is no circularity. 

Note also that  if $\{S_{\alpha}\}$ has schematically dominant transition maps and
$T \rightarrow S$ is flat then 
$\{T_{\alpha}\}$ has schematically dominant transition maps for sufficiently large indices 
because  the quasi-compact map 
$T_{\alpha} \rightarrow S_{\alpha}$ is flat for large $\alpha$
(as we see by working \'etale-locally
on $S_{\alpha_0}$ and $T_{\alpha_0}$ to reduce to the known scheme case).

\begin{proof}
Fix some $\alpha_0$ and choose an \'etale chart in schemes
$$R_0 := U_0 \times_{S_{\alpha_0}} U_0 \rightrightarrows U_0 \rightarrow S_{\alpha_0}$$
with $U_0$ qcqs \'etale over $S_{\alpha_0}$
(e.g., an affine scheme), 
so both $U_0$ and $R_0$ are qcqs schemes.  For any $\alpha \ge \alpha_0$ we define
$U_{\alpha} = U_0 \times_{S_{\alpha_0}} S_{\alpha}$
and $R_{\alpha} = R_0 \times_{S_{\alpha_0}} S_{\alpha}$, so 
$R_{\alpha} \rightrightarrows U_{\alpha}$ is naturally an \'etale scheme chart for $S_{\alpha}$
compatible with change in $\alpha$ and moreover the inverse systems
$\{U_{\alpha}\}$ and $\{R_{\alpha}\}$ have affine transition maps.  
The corresponding inverse limits $U$ and $R$ provide an \'etale
equivalence relation $R \rightrightarrows U$ over $S$ obtained by base change of
any $R_{\alpha} \rightrightarrows U_{\alpha}$ along $S \rightarrow S_{\alpha}$.
It is clear that $R \rightarrow U \times_S U$ is an isomorphism,
so $R \rightrightarrows U$ is an \'etale chart for the algebraic space $S$. 
Note that $U \rightarrow S$ is a qcqs \'etale (i.e., finitely presented \'etale) map. 

Now we turn our attention to the finitely presented map $T \rightarrow S$. 
First we consider the case when $T$ is a scheme, and then we bootstrap this to the
case when $T$ is an algebraic space.  
The pullback $U' = T \times_S U$ is quasi-compact and \'etale over $T$,
and it is a scheme since $T$ and $U$ are schemes.  Likewise,
$R' = T \times_S R$ is a quasi-compact scheme and $R' \rightrightarrows U'$
is an \'etale equivalence relation with $R' = U' \times_T U'$, so $U'/R' = T$. 
The natural \'etale $S$-map 
$q:U' \rightarrow U = \invlim U_{\alpha}$
is a qcqs map, hence it is of finite presentation, so since $U'$ and
the $U_{\alpha}$'s are schemes we may use standard limit results for schemes to increase $\alpha_0$ so that
there exists a qcqs \'etale map of schemes 
$U'_{\alpha_0} \rightarrow U_{\alpha_0}$ which descends $q$.
We use this descended map to define $U'_{\alpha_0} \rightarrow S_{\alpha_0}$,
and define $U'_{\alpha} = U'_{\alpha_0} \times_{S_{\alpha_0}} S_{\alpha}$
for all $\alpha \ge \alpha_0$.  Clearly $\{U'_{\alpha}\}_{\alpha \ge \alpha_0}$
has affine transition maps (so each $U'_{\alpha}$ is a scheme)
and $U' = \invlim U'_{\alpha}$ over $U = \invlim U_{\alpha}$.  
By increasing $\alpha_0$ if necessary, we may also descend
$R' \rightarrow R$ to a qcqs \'etale map of schemes $R'_{\alpha_0} \rightarrow R_{\alpha_0}$
and define the inverse system of schemes $R'_{\alpha} = R'_{\alpha_0} \times_{S_{\alpha_0}} S_{\alpha}$ with affine transition maps and limit $R'$.   

Consider the two \'etale $S$-maps
$p'_1, p'_2:R' \rightrightarrows U'$ over the pair of $S$-maps $p_1, p_2:R \rightrightarrows U$ between
schemes \'etale over $S$.  By Corollary \ref{mapdescend1}, if we increase
$\alpha_0$ then we can arrange that the maps $p'_1$ and $p'_2$
descend to a pair of $S_{\alpha_0}$-maps $p'_{1,\alpha_0}, p'_{2,\alpha_0}:R'_{\alpha_0}
\rightrightarrows U'_{\alpha_0}$ over the respective projections $p_{1,\alpha_0}$ and
$p_{2,\alpha_0}$.  For $\alpha \ge \alpha_0$ and $i \in \{1,2\}$, define the $S_{\alpha}$-map
$p'_{i,\alpha}:R'_{\alpha} \rightarrow U'_{\alpha}$ to be the base change of
$p'_{i,\alpha_0}$ along $S_{\alpha} \rightarrow S_{\alpha_0}$.   Now consider
only $\alpha \ge \alpha_0$.  
 In particular, each $p'_{i,\alpha}$ is a qcqs \'etale $S_{\alpha}$-map
and base change along $S_{\beta} \rightarrow S_{\alpha}$ for $\beta \ge \alpha$
carries $p'_{i,\alpha}$ to $p'_{i,\beta}$.  

We claim that $R'_{\alpha} \rightrightarrows U'_{\alpha}$
is an \'etale equivalence relation for sufficiently large $\alpha$.  In other words, for large $\alpha$ 
we claim that the maps $p'_{i,\alpha}$ are
\'etale and that the map $R'_{\alpha} \rightarrow U'_{\alpha} \times_{\Spec \Z} U'_{\alpha}$
is a monomorphism which functorially defines an equivalence relation.  This latter map factors
through the subfunctor $U'_{\alpha} \times_{S_{\alpha}} U'_{\alpha}$ which is also a scheme, so
for the equivalence relation condition it is equivalent to check that the map $\delta_{\alpha}:
R'_{\alpha} \rightarrow U'_{\alpha} \times_{S_{\alpha}}
U'_{\alpha}$ is a monomorphism defining a functorial equivalence relation for large $\alpha$.
Note that $\delta_{\alpha}$ is a finitely presented map since it is an $S_{\alpha}$-map
between schemes of finite presentation along over $S_{\alpha}$.  

Since $\{p'_{i,\alpha}\}$ is a compatible system of finitely presented maps
between inverse systems of qcqs {\em schemes} with affine transition maps 
and the limit map $p'_i:R' \rightarrow U'$ is \'etale, by
\cite[IV$_4$,~17.7.8(ii)]{ega} it follows that $p'_{1,\alpha}$ and $p'_{2,\alpha}$
are \'etale for all large $\alpha$.  
Likewise, $\{U'_{\alpha} \times_{S_{\alpha}} U'_{\alpha}\}$ is 
an inverse system of qcqs {\em schemes} with affine transition maps and limit $U' \times_S U'$,
and base change of the finitely presented
scheme map $\delta_{\alpha}$ along a transition map in this system (resp. 
along pullback to the limit) is identified with the base change along the corresponding
transition map in $\{S_{\alpha}\}$ (resp. along pullback from $S_{\alpha}$ to $S$). 
Since the limit map $\delta = (p'_1, p'_2):R' \rightarrow U' \times_S U'$ is an equivalence relation,
in particular it is a monomorphism.  Thus, by \cite[IV$_4$,~8.10.5(i bis)]{ega}
the $S_{\alpha}$-maps $\delta_{\alpha}$ are monomorphisms for sufficiently large $\alpha$.
Likewise, if $\alpha$ is large enough then the symmetry automorphism $R' \simeq R'$ over the flip on
the scheme $U' \times_S U'$ carrying $(p'_1, p'_2)$ to $(p'_2, p'_1)$ descends to such 
an automorphism $R'_{\alpha} \simeq R'_{\alpha}$ over
the flip on the scheme $U'_{\alpha} \times_{S_{\alpha}} U'_{\alpha}$ carrying
$(p'_{1,\alpha}, p'_{2,\alpha})$ to $(p'_{2,\alpha}, p'_{1,\alpha})$, 
so the functorial relation $R'_{\alpha}$
on $U'_{\alpha}$ is symmetric for sufficiently large $\alpha$.   Similarly, the transitivity morphism
of schemes $R' \times_{p'_2, U', p'_1} R' \rightarrow R'$ induced by the scheme map 
$q_{13}:U' \times_S U' \times_S U'
\rightarrow U' \times_S U'$ descends to the $\alpha$-level for sufficiently large $\alpha$, 
and the diagonal $\Delta_{U'_{\alpha}/S_{\alpha}}$ factors through $\delta_{\alpha}$
for sufficiently large $\alpha$.  Hence, for some large $\alpha_1$ we have that 
$R'_{\alpha}$ defines an equivalence relation on $U'_{\alpha}$ for all $\alpha \ge \alpha_1$.

The quotient $T_{\alpha} := U'_{\alpha}/R'_{\alpha}$ 
makes sense as a quasi-compact algebraic space locally
of finitely presentation over $U_{\alpha}/R_{\alpha} = 
S_{\alpha}$ for $\alpha \ge \alpha_1$, and it is quasi-separated 
since $R'_{\alpha} \rightarrow U'_{\alpha}
\times_{S_{\alpha}} U'_{\alpha}$ is quasi-compact and $S_{\alpha}$ is quasi-separated.
In particular, $\{T_{\alpha}\}$ is an inverse system of qcqs algebraic
spaces.  We clearly have $S_{\alpha}$-isomorphisms
$T_{\alpha} \simeq T_{\alpha_1} \times_{S_{\alpha_1}} S_{\alpha}$
compatibly with change in $\alpha \ge \alpha_1$, so $\{T_{\alpha}\}$ has
affine transition maps and there is an $S$-isomorphism 
$T \simeq T_{\alpha_1} \times_{S_{\alpha_1}} S$.  In particular, 
$T$ is identified with $\invlim T_{\alpha}$ over $S$.  

But we are considering the case when $T$ is a scheme, so
by Lemma \ref{scheme} it follows that $T_{\alpha}$ is a scheme for sufficiently large $\alpha$
provided that some $T_{\alpha'}$ is affine over an algebraic space
of finite presentation over $\Z$.   For example, we have have solved the case when $T$ is a scheme
(with the stronger conclusion that we can choose the $T_{\alpha}$ to be schemes) provided 
that the $S_{\alpha}$ are of finite presentation over $\Z$, as then all
$T_{\alpha}$ above are of finite presentation over $\Z$ as well (so Lemma \ref{scheme}
is applicable to the inverse system $\{T_{\alpha}\}$). 

Now suppose more generally that $T$ is an algebraic space, but assume that
all $S_{\alpha}$ are of finite presentation over $\Z$.  (The reason for this temporary extra assumption is
that the case when $T$ is a scheme has thus far been settled only in such cases,
though with the stronger conclusion that the $T_{\alpha}$ can be chosen to be schemes in such cases.)  
Fix a qcqs \'etale cover $U' \twoheadrightarrow T$ by a scheme, so
$R' = U' \times_T U'$ is also a scheme and both $U'$ and $R'$ are
qcqs over $S = \invlim S_{\alpha}$.  We may apply to $U'$
and $R'$ the settled case of schemes (with all $S_{\alpha}$ of
finite presentation over $\Z$), so by 
restricting to sufficiently large $\alpha$ we can construct cartesian limit presentations
$U' = \invlim U'_{\alpha}$ and $R' = \invlim R'_{\alpha}$
in finitely presented schemes (with affine transition maps) over the $S_{\alpha}$'s.  

By Corollary \ref{mapdescend1}, the $S$-map
$\delta':R' \rightarrow U' \times_S U'$ arises
from a compatible system of $S_{\alpha}$-maps
$\delta'_{\alpha}:R'_{\alpha} \rightarrow U'_{\alpha} \times_{S_{\alpha}} U'_{\alpha}$
for large $\alpha$.  For fixed $i \in \{1, 2\}$, the system of $S_{\alpha}$-maps
$p'_{i,\alpha}:R'_{\alpha} \rightarrow U'_{\alpha}$
between finitely presented $S_{\alpha}$-schemes is cartesian with respect to change in $\alpha$
and has limit $p'_i:R' \rightarrow U'$ that is \'etale.  
Exactly as we have just argued above in our treatment of the case when $T$ is a scheme, 
by restricting to large $\alpha$ we can arrange that each $p'_{i,\alpha}$ is \'etale
and that each $\delta'_{\alpha}$ is a monomorphism. 
A further application of the same method handles 
the symmetry and transitivity aspects, so $R'_{\alpha} \rightrightarrows
U'_{\alpha}$ is a cartesian inverse system of \'etale equivalence relations in finitely presented
schemes over $S_{\alpha}$ for large $\alpha$.  

The quotients
$T_{\alpha} = U'_{\alpha}/R'_{\alpha}$ form a cartesian inverse system 
of finitely presented algebraic spaces over $S_{\alpha}$,  so
$\{T_{\alpha}\}$ has affine transition maps.  The limit
$T' := \invlim T_{\alpha}$ therefore makes sense as an algebraic space over
$S$ and for a large $\alpha_0$ we have an $S$-isomorphism 
$$T' \simeq T_{\alpha_0} \times_{S_{\alpha_0}} S = (U'_{\alpha_0}/R'_{\alpha_0}) \times_{S_{\alpha_0}} S
\simeq U'/R' = T.$$
Thus, we have solved the general problem when the $S_{\alpha}$ are all of finite presentation
over $\Z$, including the refined claim that (in such cases) when $T$ is a scheme
we may choose the $T_{\alpha}$ to all be schemes. 

Finally, we consider the general case by relaxing the assumption that all $S_{\alpha}$ are of
finite presentation over $\Z$ to the assumption that some $S_{\alpha'}$ is affine 
over an algebraic space $\mathscr{S}$ of finite presentation over $\Z$.  
We only consider $\alpha \ge \alpha'$, so all $S_{\alpha}$ are affine over $\mathscr{S}$.  
Arguing via \cite[III,~Thm.~1.1,~Cor.~1.2]{knutson}
as in the proof of Lemma \ref{scheme}, we get another limit presentation
$S = \invlim S'_{\gamma}$ with 
$\{S'_{\gamma}\}$ an inverse system with affine transition maps
of algebraic spaces over $\mathscr{S}$ with each map 
$S'_{\gamma} \rightarrow \mathscr{S}$ affine
and finitely presented.  Thus, all $S'_{\gamma}$ are finitely presented over
$\Z$, so the settled cases imply that the finitely presented map
$T \rightarrow S$ descends to a finitely presented
map $T_{0} \rightarrow S'_{\gamma_0}$ for some $\gamma_0$,
and that we can choose $T_{0}$ to be a scheme when $T$ is a scheme. 

Apply Proposition \ref{fpres} to the finitely presented map $S'_{\gamma_0} \rightarrow \mathscr{S}$
and the inverse system $\{S_{\alpha}\}$ over $\mathscr{S}$ with limit $S$.
This gives that the canonical $\mathscr{S}$-map $S = \invlim S_{\alpha} \rightarrow S'_{\gamma_0}$
factors through an $\mathscr{S}$-map $S_{\alpha_0} \rightarrow S'_{\gamma_0}$
for some $\alpha_0$.  Note that this map is affine, since $S_{\alpha_0}$
and $S'_{\gamma_0}$ are affine over $\mathscr{S}$.
Define $T_{\alpha_0} = T_{0} \times_{S'_{\gamma_0}} S_{\alpha_0}$
and $$T_{\alpha} := T_{\alpha_0} \times_{S_{\alpha_0}} S_{\alpha} = T_0 \times_{S'_{\gamma_0}} S_{\alpha}$$
for all $\alpha \ge \alpha_0$.
Then $\{T_{\alpha}\}$ is an inverse system with affine transition maps of finitely presented
objects over $\{S_{\alpha}\}$ and
$$T_{\alpha_0} \times_{S_{\alpha_0}} S = T_{0} \times_{S'_{\gamma_0}} S = T$$
over $S$.  Also, if $T$ is a scheme then we can choose
$T_{0}$ to be a scheme, so then every $T_{\alpha}$ is a scheme since 
it is affine over $T_0$.  Thus, we are done.  
\end{proof}

\begin{corollary}\label{mapdescend2}
In the setup of Corollary {\rm{\ref{mapdescend1}}}, 
let $f_{\alpha}:X_{\alpha} \rightarrow Y_{\alpha}$
over $S_{\alpha}$ be a compatible system maps inducing $f:X \rightarrow Y$ over $S$
in the limit. Assume that some $S_{\alpha'}$ is affine over an algebraic space of
finite presentation over $\Z$.

The map $f$ satisfies $\P$ if and only if 
$f_{\alpha}$ does for sufficiently large $\alpha$, where
$\P$ is any of the properties of morphisms as in
{\rm{\cite[IV$_3$,~8.10.5(i)--(xii),~11.2.6(ii);~IV$_4$,~17.7.8(ii)]{ega}}}.
\end{corollary}

The hypothesis that some $S_{\alpha'}$ is affine over an algebraic space
of finite presentation over $\Z$ is needed solely because of its role in 
the previous two results (which are used in the following proof),
so as with those results we can eliminate this hypothesis once
absolute noetherian approximation is proved for qcqs algebraic spaces in 
\S\ref{absnoetherian}.   The proof of this approximation result uses
Corollary \ref{mapdescend2} only in cases with all $S_{\alpha}$ of finite presentation over $\Z$,
so there is no circularity. 

\begin{proof}
The descent of properties of morphisms which
are \'etale-local on the source and target does not require the hypothesis
on some $S_{\alpha'}$ and is easily reduced to the known case of
schemes treated in the given references in \cite[IV]{ega}.  The remaining properties
(all in \cite[IV$_3$,~8.10.5]{ega}) are \'etale-local on the base, so we
may assume that all $S_{\alpha}$, $S$, $Y_{\alpha}$, and $Y$ are affine schemes.
The key point is that if $f$ is representable in schemes (i.e., if $X$ is a scheme) then so is
$f_{\alpha}$ for large $\alpha$ (i.e., $X_{\alpha}$ is a scheme for large $\alpha$), 
as follows from Lemma \ref{scheme}.  This reduces everything
to the known case of schemes except for the properties of being
surjective, radiciel, quasi-finite, or proper.   The first three are fibral properties
cutting out a constructible locus in the base, and so the
proofs for these conditions in the scheme case carry over to the case of algebraic spaces.
Indeed, constructible loci
in qcqs algebraic spaces interact with limits
exactly as in the scheme case \cite[IV$_3$,~8.3.4]{ega}
(as one shows by working \'etale-locally
to reduce to the case of schemes), so the behavior of
fibral properties with respect to limits as in \cite[IV$_3$,~9.3.3]{ega} holds for
algebraic spaces too.  

As for properness, if $f$ is proper
then the finitely presented map $f_{\alpha}:X_{\alpha} \rightarrow Y_{\alpha}$
is at least separated for 
large $\alpha$.  Fixing such an $\alpha_0$, since
$Y_{\alpha_0}$ is affine 
we can express $Y_{\alpha_0}$ as the limit
of an inverse system of affine schemes of finite type over 
$\Z$ and then use Proposition \ref{lemma7}
and the settled descent of
separatedness to descend $f_{\alpha_0}$ to a finitely presented and separated 
map $X_0 \rightarrow Y_0$
from an algebraic space to a finite type $\Z$-scheme.  By Chow's
Lemma for separated maps of finite type between noetherian algebraic spaces, there is 
a surjective, proper, and finitely presented map $P_{\alpha_0} \rightarrow X_{\alpha_0}$
with $P_{\alpha_0}$ a scheme.   The induced 
$S_{\alpha}$-maps $P_{\alpha} \rightarrow Y_{\alpha}$ for $\alpha \ge \alpha_0$
only involve schemes, and the limit map $P \rightarrow Y$ over $S$ is proper
since $f:X \rightarrow Y$ is proper.  Hence, by the known scheme case
$P_{\alpha}$ is $Y_{\alpha}$-proper for large $\alpha$, so surjectivity of 
the $Y_{\alpha}$-map $P_{\alpha} \rightarrow X_{\alpha}$ 
forces the separated and finite type map $X_{\alpha} \rightarrow Y_{\alpha}$
to be proper for large $\alpha$.
\end{proof}

We end this section by recording a result on how inverse limit
presentations of quasi-compact open subschemes can be extended
to such presentations of an ambient qcqs scheme at the expense of passing to a
cofinal subsystem of the given inverse system.  

\begin{lemma}\label{lemma2}
Let $X$ be a qcqs scheme, and $V \subseteq X$ a quasi-compact open subscheme. 
Assume that $X$ is endowed with a structure of $\Lambda$-scheme for a noetherian
ring $\Lambda$ $($e.g., $\Lambda = \Z$$)$, and 
choose an inverse system $\{V_i\}_{i \in I}$ of finite type $\Lambda$-schemes 
with affine and schematically dominant transition maps such that $V \simeq \invlim V_i$.

This can be extended to a limit presentation of $X$ in the sense that 
there is an inverse system $\{X_j\}_{j \in J}$ of finite type
$\Lambda$-schemes with
affine and schematically dominant 
transition maps, a cofinal map of directed sets $\psi:J \rightarrow I$, and 
a compatible system of open immersions $V_{\psi(j)} \hookrightarrow X_j$ over $\Lambda$
such that the diagrams
$$\xymatrix{{V_{\psi(j')}} \ar[d] \ar[r] & {X_{j'}} \ar[d] \\ {V_{\psi(j)}} \ar[r] & {X_j}}$$
are cartesian for all $j' \ge j$ and the open immersion
$V = \invlim V_i \simeq \invlim V_{\psi(j)} \hookrightarrow \invlim X_j$
extends to an isomorphism $X \simeq \invlim X_j$.
\end{lemma}

Such an inverse system $\{V_i\}$ always exists for any $V \subseteq X$, 
 by \cite[Thm.~C.9]{tt} and an argument with scheme-theoretic images
as in the proof of Lemma \ref{zapprox}.  

\begin{proof}
Choose a non-empty finite collection of affine open subschemes $\{U_1, \dots, U_m\}$
of $X$ such that $V$ together with the $U_k$ cover $X$.  (For example, we
could take  $\{U_k\}$ to be an affine open covering of $X$.)  We shall argue by
induction on $m$.  The case $m = 1$ is established in the proof of \cite[Thm.~C.9]{tt},
and here it is used that $\Lambda$ is noetherian and the affine transition maps
in $\{V_i\}$ are schematically dominant.  
In general, if $m > 1$ then we may apply induction to $V' = V \cup (U_1 \cup \dots \cup
U_{m-1})$ in the role of $X$ and then apply the case $m=1$ to the inclusion of
$V'$ into $X$ (with $X = V' \cup U_m$). 
\end{proof}

\subsection{Square-zero thickenings of algebraic spaces}\label{sqzerosec}

In order to reduce the proof of Nagata's theorem to the reduced case, we need to relate
square-zero thickenings of  an algebraic space $X$ to (suitable) square-zero extensions of 
the sheaf of rings ${\ms O}_X$ on $X_{\et}$. It is convenient to express this in terms of an equivalence
of categories, as follows.

Let $\mathscr{C}'$ denote the category of pairs $(X', \mathscr{I})$
consisting of an algebraic space $X'$ and 
a quasi-coherent sheaf of ideals $\mathscr{I} \subseteq {\ms O}_{X'}$ on $X'_{\et}$ 
such that $\mathscr{I}^2 = 0$.  Let $\mathscr{C}$ denote the category of pairs 
$(X, \mathscr{A} \twoheadrightarrow {\ms O}_X)$ consisting of an algebraic space $X$
and a surjective map
$\theta:\mathscr{A} \twoheadrightarrow {\ms O}_X$  of sheaves of rings on $X_{\et}$ 
whose kernel $\mathscr{J} \subseteq \mathscr{A}$ 
is a square-zero ideal that is quasi-coherent as an 
${\ms O}_X$-module.

Fix an object $(X', \mathscr{I})$ in $\mathscr{C}'$.
Let $X$ denote the closed subspace of $X'$ cut out by $\mathscr{I}$.
Pullback along $X_{\et} \rightarrow X'_{\et}$ induces an equivalence
of topoi, by Proposition \ref{redequiv}, so for a sheaf of sets $\mathscr{F}'$ 
on $X'_{\et}$ we will therefore abuse notation 
(when the context makes it clear) by also writing $\mathscr{F}'$ to denote
the pullback sheaf on $X_{\et}$.
For example, we view ${\ms O}_{X'}$ as a sheaf of rings on $X_{\et}$ in this way,
and $\mathscr{I}$ as
a square-zero sheaf of ideals in ${\ms O}_{X'}$ on $X_{\et}$. There is an evident quotient map 
${\ms O}_{X'} \twoheadrightarrow
{\ms O}_X$ of sheaves on $X_{\et}$ with kernel
$\mathscr{I}$ whose induced module structure over ${\ms O}_{X'}/\mathscr{I} \simeq {\ms O}_X$
on $X_{\et}$ is clearly quasi-coherent.
Thus, $(X', \mathscr{I}) \rightsquigarrow (X, {\ms O}_{X'} \twoheadrightarrow {\ms O}_X)$
is a functor from $\mathscr{C}'$ to $\mathscr{C}$.

\begin{theorem}\label{sqzero} The above functor $\mathscr{C}' \rightarrow \mathscr{C}$
is an equivalence of categories. 
\end{theorem}

This theorem says that for any algebraic space $X$, to give a square-zero
thickening of $X$ as an algebraic space is functorially the same as to choose a quasi-coherent
sheaf $\mathscr{J}$ on $X_{\et}$ and a sheaf of rings $\mathscr{A}$ on $X_{\et}$
that is a square-zero extension of ${\ms O}_X$ by the ${\ms O}_X$-module $\mathscr{J}$.  This latter
point of view is expressed entirely in terms of the \'etale topos of $X$ and so
is well-suited to deformation-theoretic considerations (such as with the cotangent complex).  

\begin{remark}\label{remsqzero}
 Although we state and prove Theorem \ref{sqzero} without noetherian hypotheses, it is only
used in this paper the noetherian case (in \S\ref{redredcase}).   Aside from using 
a couple of self-contained 
results (Lemma \ref{lemma3} and Proposition \ref{redequiv}), our proof uses
Corollary \ref{corred}, which rests on Theorem \ref{absapprox}, 
whose proof in turn relies on almost everything in \S\ref{absnoetherian} and \S\ref{limprop}.
  However, 
the noetherian case of Corollary \ref{corred} is an old result of Knutson \cite[III,~Thm.~3.3]{knutson}.
Consequently, the reader only interested in \S\ref{exccase} may safely restrict
attention to the simpler noetherian case with quasi-coherence replaced with coherence in several
places (since a square-zero extension of a noetherian
ring $A_0$ by a finite $A_0$-module is automatically a noetherian ring).
\end{remark}

\begin{proof}[Proof of Theorem $\ref{sqzero}$]

{\bf Step 1}. (Faithfulness) For objects $(X'_1, \mathscr{I}_1)$
and $(X'_2, \mathscr{I}_2)$ in $\mathscr{C}'$ consider
two maps $f', h':X'_1 \rightrightarrows X'_2$ such that
(i) $f'$ and $h'$ carry $X_1$ into $X_2$,
(ii) the induced maps $f, h:X_1 \rightrightarrows X_2$ coincide,
and (iii) the resulting pairs of maps
$${\ms O}_{X'_2} \rightrightarrows f_{\ast}({\ms O}_{X'_1}) = h_{\ast}({\ms O}_{X'_2})$$
of sheaves on $(X_2)_{\et}$ coincide.  We wish to show that $f' = h'$.
In case $X'_1$ and $X'_2$ are schemes, so morphisms between them can be considered
within the category of ringed spaces (rather than ringed topoi), the desired equality is obvious.
In general we will use \'etale covers by schemes to reduce the faithfulness problem to the 
settled scheme case.

Pick an \'etale covering $\pi'_2:U'_2 \rightarrow X'_2$ by a scheme, 
and let $\pi_2:U_2 \rightarrow X_2$ denote the  pullback \'etale scheme covering. 
The resulting pullback \'etale covers
$$\xymatrix{
{{f'}^{-1}(U'_2)} \ar[dr]_-{p'} && {h'}^{-1}(U'_2) \ar[dl]^-{q'} \\
& {X'_1} & }$$
restrict over $X_1 \subseteq X'_1$ to algebraic space
\'etale covers $f^{-1}(U_2)$ and $h^{-1}(U_2)$ of $X_1$ that are naturally
identified since $f = h$.  
This $X_1$-isomorphism $\phi:{f'}^{-1}(U'_2)|_{X_1} \simeq {h'}^{-1}(U'_2)|_{X_1}$ 
between algebraic spaces \'etale over $X_1$ uniquely lifts to an $X'_1$-isomorphism
$\phi':{f'}^{-1}(U'_2) \simeq {h'}^{-1}(U'_2)$ between algebraic spaces \'etale over $X'_1$, 
due to the topological
invariance of the \'etale site of algebraic spaces (Proposition \ref{redequiv}).    
Since $p'$ and $q'$ are \'etale surjections and $q' \circ \phi' = p'$, 
to prove the equality $f' = h'$ it suffices to prove that 
the diagram
$$\xymatrix{
{{f'}^{-1}(U'_2)} \ar[rr]_-{\simeq}^-{\phi'} \ar[dr]_-{f' \circ p'} & & {{h'}^{-1}(U'_2)} 
\ar[dl]^-{h' \circ q'} \\ & {X'_2} &}$$
commutes.

Choose an \'etale scheme cover $\psi':U'_1 \rightarrow {f'}^{-1}(U'_2)$ and define the \'etale
scheme covering map $U'_1 \rightarrow {h'}^{-1}(U'_2)$ 
to be $\phi' \circ \psi'$, so 
these define a common \'etale covering map $\pi'_1:U'_1 \rightarrow X'_1$
since $\phi'$ is an $X'_1$-isomorphism.  
Composing $f' \circ p'$ and $h' \circ q'$ back to $U'_1$ via $\psi'$ and $\phi' \circ \psi'$
respectively recovers the pair of maps
$f' \circ \pi'_1, h' \circ \pi'_1:U'_1 \rightrightarrows X'_2$, so we get a co-commutative diagram
$$\xymatrix{
{U'_1} \ar@<2pt>[r]^-{\widetilde{f}'} \ar@<-2pt>[r]_-{\widetilde{h}'} \ar[d]_-{\pi'_1} & {U'_2} \ar[d]^-{\pi'_2} \\
{X'_1} \ar@<2pt>[r]^-{f'} \ar@<-2pt>[r]_-{h'} & {X'_2}}$$
in which $\widetilde{f}'$ is 
$U'_1 \stackrel{\psi'}{\rightarrow} {f'}^{-1}(U'_2) \rightarrow U'_2$ and
$\widetilde{h}'$ is $U'_1 \stackrel{\phi' \circ \psi'}{\longrightarrow}
{h'}^{-1}(U'_2) \rightarrow U'_2$.  
But the pair $(\widetilde{f}', \widetilde{h}')$ satisfies
the same initial hypotheses as the pair $(f', h')$, using
the quasi-coherent  square-zero ideals given by the 
pullbacks of $\mathscr{I}'_1$ to $U'_1$ and 
$\mathscr{I}'_2$ to $U'_2$.  Hence, by the settled scheme case we conclude that
$\widetilde{f}' = \widetilde{h}'$, so $f' = h'$.

{\bf Step 2}. (Fullness: reduction to schemes) For a pair of objects
$(X'_1, \mathscr{I}_1)$ and $(X'_2, \mathscr{I}_2)$ in $\mathscr{C}'$, suppose
we are given a map $f:X_1 \rightarrow X_2$ of algebraic spaces and 
a map $\theta':{\ms O}_{X'_2} \rightarrow f_{\ast}({\ms O}_{X'_1})$ of sheaves
of rings on $(X_2)_{\et}$ (using the equivalences $(X_1)_{\et} = (X'_1)_{\et}$
and $(X_2)_{\et} = (X'_2)_{\et}$ to make sense of $f_{\ast}({\ms O}_{X'_1})$) such that the diagram
\begin{equation}\label{cartgiven}
\xymatrix{
{{\ms O}_{X'_2}} \ar[r]^-{\theta'} \ar[d] & {f_{\ast}({\ms O}_{X'_1})} \ar[d] \\
{{\ms O}_{X_2}} \ar[r]_-{f^{\sharp}} & {f_{\ast}({\ms O}_{X_1})}}
\end{equation}
of sheaves on $(X_2)_{\et}$ commutes.  We seek to construct a map $f':X'_1 \rightarrow X'_2$ inducing 
the pair $(f, \theta')$.  

By Step 1 such an $f'$ is unique if it exists, so by descent
we can compose along an \'etale scheme covering of $X'_1$ to reduce to the
case when $X'_1$ (and hence $X_1$) is a scheme.
We will reduce the proof of existence of $f'$ to the case when $X'_2$
is also a scheme.  In the case that $X'_2$ is a scheme there is trivially such an $f'$ which
uniquely solves our problem for the Zariski topology since we can then work with
locally ringed spaces, but even in this scheme case there is more to do:  
we have to show that the scheme morphism which is a solution for the Zariski topology
is a solution for the \'etale topology.  Granting this additional property in the scheme case
for now, let us see how to solve the general case.

Choose an \'etale scheme cover $\pi'_2:U'_2 \rightarrow X'_2$
and form the cartesian square
$$\xymatrix{
{U_2} \ar[r] \ar@{->>}[d]_-{\pi_2} &{U'_2} \ar@{->>}[d]^-{\pi'_2} \\ {X_2} \ar[r] & {X'_2}}$$
in which the horizontal arrows are closed immersions defined by square-zero ideals.
In particular, the left side is an \'etale scheme covering. Pulling back $\pi_2$
along $f:X_1 \rightarrow X_2$ defines an \'etale cover
$$\pi_1:U_1 := X_1 \times_{X_2} U_2 \rightarrow X_1$$
with $U_1$ a scheme since $X_1$ and $U_2$ are schemes, and
by topological invariance of the \'etale site
of algebraic spaces (Proposition \ref{redequiv}) this fits into a unique cartesian square of
algebraic spaces
$$\xymatrix{
{U_1} \ar@{->>}[d]_-{\pi_1} \ar[r] & {U'_1} \ar@{->>}[d]^-{\pi'_1} \\
{X_1} \ar[r] & {X'_1}}$$
such that $\pi'_1$ is \'etale (and necessarily surjective, since $\pi_1$ is).  The underlying reduced space 
$(U'_1)_{\rm{red}} = (U_1)_{\rm{red}}$ is a scheme, so by 
\cite[III,~Thm.~3.3]{knutson} in the noetherian case and by Corollary \ref{corred} 
in the general case 
we conclude that 
$U'_1$ is actually a scheme.  (See Remark \ref{remsqzero}.)  

The natural map $\widetilde{f} = 
\pi_2^{\ast}(f):U_1 = X_1 \times_{X_2} U_2 \rightarrow U_2$ fits into the top row of the 
diagram of schemes
\begin{equation}\label{udiagram}
\xymatrix{
{U_1} \ar[r]^-{\widetilde{f}} \ar[d] & {U_2} \ar[d] \\
{U'_1} & {U'_2}}
\end{equation}
in which the vertical maps are closed immersions defined by the square-zero pullbacks of
$\mathscr{I}_1$ and $\mathscr{I}_2$ along $\pi'_1:U'_1 \rightarrow X'_1$ and 
$\pi'_2:U'_2 \rightarrow X'_2$ respectively.
The diagram (\ref{cartgiven}) of sheaves of rings on $(X_2)_{\et}$
restricts over the subcategory $(U_2)_{\et}$ to give a commutative diagram
\begin{equation}\label{ddig}
\xymatrix{
{{\ms O}_{U'_2}} \ar[r]^-{\theta'|_{U_2}} \ar[d] & {\widetilde{f}_{\ast}({\ms O}_{U'_1})} \ar[d] \\
{{\ms O}_{U_2}} \ar[r]_-{\widetilde{f}^{\sharp}} & {{\ms O}_{U_1}}}
\end{equation}
in which the vertical maps are the natural surjections.
Hence, the pair $(\widetilde{f}, \theta'|_{U_2})$ is a morphism 
\begin{equation}\label{cmap}
(U_1, {\ms O}_{U'_2} \twoheadrightarrow {\ms O}_{U_1}) \rightarrow 
(U_2, {\ms O}_{U'_2} \twoheadrightarrow {\ms O}_{U_2})
\end{equation}
in the category $\mathscr{C}$.  Since we are assuming that
the fullness problem is solved for schemes (with the \'etale topology), 
the morphism (\ref{cmap}) arises from a unique morphism
$\widetilde{f}':U'_1 \rightarrow U'_2$ of algebraic spaces fitting into the bottom side of (\ref{udiagram})
as a commutative diagram 
and compatible with $\theta'|_{U_2}$.

Consider the resulting diagram
\begin{equation}\label{diag}
\xymatrix{
{U'_1 \times_{X'_1} U'_1} \ar@<2pt>[r] \ar@<-2pt>[r] & {U'_1} \ar[r]^-{\pi'_1} \ar[d]_-{\widetilde{f}'} &
{X'_1} \ar@{-->}[d]_-{?}^-{f'} \\
& {U'_2} \ar[r]_-{\pi'_2} & {X'_2}}
\end{equation}
Composing the map $\pi'_2 \circ \widetilde{f}':U'_1 \rightarrow X'_2$ 
with the projections $U'_1 \times_{X'_1} U'_1 \rightrightarrows U'_1$
gives a pair of maps $U'_1 \times_{X'_1} U'_1 \rightrightarrows X'_2$ 
that coincide, due to Step 1 (using the square-zero quasi-coherent pullback
of $\mathscr{I}_1$ on $U'_1 \times_{X'_1} U'_1$ and
the commutativity of (\ref{cartgiven}) and (\ref{ddig})), so by descent 
we can uniquely fill in the arrow $f':X'_1 \rightarrow X'_2$ in (\ref{diag}) 
to make a commutative square.  It is easy to check (with the help of $\widetilde{f}'$) that 
$f'$ is the morphism from $(X'_1, \mathscr{I}_1)$ to
$(X'_2, \mathscr{I}_2)$ in $\mathscr{C}'$ that we sought to construct.
This completes the reduction of the fullness problem to the special case
of schemes (with the \'etale topology!).  We will address this in Step 3.

{\bf Step 3}.  (The case of schemes)  We now prove that the faithful functor 
$\mathscr{C}' \rightarrow \mathscr{C}$ is full on the full subcategory of pairs
$(X', \mathscr{I})$ for which $X'$ is a scheme (with the \'etale topology) 
and that its essential image contains
all pairs $(X, \mathscr{A} \twoheadrightarrow {\ms O}_X)$ for which $X$ is a scheme
(with the \'etale topology).   The fullness for schemes is trivial to check for the Zariski topology
using the viewpoint of locally ringed spaces.  Fullness for schemes with the \'etale topology
then follows formally from two facts:  (i) for a scheme $X'$, the stalks of
${\ms O}_{X'_{\et}}$ at geometric points of $X'_{\et}$ are the strict henselizations of
the Zariski local rings of ${\ms O}_{X'}$, and (ii) the uniqueness aspect of the universal mapping property for
the strict henselization of a local ring (relative to a specified separable closure 
of its residue field).  By Step 2, fullness has now been proved in general. 

Now consider a pair $(X, \mathscr{A} \twoheadrightarrow {\ms O}_X)$ for a scheme $X$ 
with the \'etale topology.   We seek to construct a pair $(X', \mathscr{I})$ in
$\mathscr{C}'$ giving rise to this, and will do so
with $X'$ a scheme.  By the settled full faithfulness, it suffices to solve
this problem Zariski-locally on $X$, so we may assume $X = \Spec A$ is an affine
scheme.  We then let $A' = \mathscr{A}(X)$, so $A' \rightarrow A$ is surjective
since $X$ has vanishing higher Zariski cohomology for quasi-coherent
${\ms O}_X$-modules.  The kernel $I = \ker(A' \twoheadrightarrow A)$
is a square-zero ideal in $A'$, and so $(\Spec A', \widetilde{I})$
with the \'etale topology is the natural candidate to consider for $(X', \mathscr{I})$.
Consideration of Zariski stalks shows that 
that the natural map ${\ms O}_{\Spec A'} \rightarrow \mathscr{A}|_{X_{\rm{Zar}}}$ is
an isomorphism.  Our problem is to show that this (necessarily uniquely) lifts to
an isomorphism ${\ms O}_{\Spec A', \et} \simeq \mathscr{A}$ of sheaves of rings on $X_{\et}$
respecting the identification ${\ms O}_{\Spec A} = {\ms O}_X$ of quotient rings on $X_{\et}$
and the identification of ideal sheaves $\widetilde{I}_{\et} \simeq \mathscr{I} := 
\ker(\mathscr{A} \twoheadrightarrow {\ms O}_X)$ on $X_{\et}$. 

Let $X' = \Spec A'$.  Consider geometric points $\overline{x}$ of $X'_{\et} = 
X_{\et}$ that are algebraic over
their physical image points $x \in X$. 
Observe that $\mathscr{A}_{\overline{x}}$ 
is a strictly henselian local ring with residue field $k(\overline{x})$ since 
its quotient ${\ms O}_{X_{\et}}$ by a square-zero ideal sheaf has strictly henselian 
local $\overline{x}$-stalk with residue field $k(\overline{x})$.  
Letting $\pi:X_{\et} \rightarrow X_{\rm{Zar}}$ be the natural map of sites, the
given isomorphism ${\ms O}_{X'_{\rm{Zar}}} \simeq \pi_{\ast}\mathscr{A}$ corresponds
to a map $\pi^{-1}{\ms O}_{X'_{\rm{Zar}}} \rightarrow \mathscr{A}$ on $X_{\et}$
that on $\overline{x}$-stalks is the natural local map
$\theta_{\overline{x}}:{\ms O}_{X',x} \rightarrow \mathscr{A}_{\overline{x}}$ inducing
$k(x) \rightarrow k(\overline{x})$ on residue fields.  But the map on $\overline{x}$-stalks
induced by the natural map
$\pi^{-1}{\ms O}_{X'_{\rm{Zar}}} \rightarrow {\ms O}_{X'_{\et}}$ is 
identified with the unique local map ${\ms O}_{X',x} \rightarrow {\ms O}_{X',\overline{x}}^{\rm{sh}}$
lifting $k(x) \rightarrow k(\overline{x})$ on residue fields, so
by the universal property of strict henselization there is at most one way to fill in the 
dotted arrow in the commutative square 
\begin{equation}\label{diagmap}
\xymatrix{
{\pi^{-1}{\ms O}_{X'_{\rm{Zar}}}} \ar[r] \ar[d] & {{\ms O}_{X'_{\et}}} \ar@{-->}[dl]^-{?}  \ar@{->>}[d]^-{\phi} \\
{\mathscr{A}} \ar@{->>}[r]_-{\psi} & {{\ms O}_{X_{\et}}}}
\end{equation}
to make commutative triangles of sheaves of rings.  (Such a diagonal map is necessarily local on stalks
and an isomorphism on residue fields since $\phi$ and $\psi$ are quotients
by square-zero ideal sheaves.)  The lower triangle commutes
if the upper one does, by the universal property of strict henselizations.
Moreover, both ideal sheaves $\ker \phi$ and $\ker \psi$ are quasi-coherent on $X_{\et}$, 
and a map between quasi-coherent sheaves for the \'etale topology of a scheme is
an isomorphism if and only if it is an isomorphism for the Zariski topology.  Thus, since 
a diagonal map as in (\ref{diagmap}) must be an isomorphism on $X'_{\rm{Zar}} = X_{\rm{Zar}}$
if it exists (as the left and top sides are isomorphisms over $X'_{\rm{Zar}}$), 
 it must induce an isomorphism
between the quasi-coherent $\ker \phi$ and $\ker \psi$ on $X_{\et}$ and hence be an isomorphism
if it exists.  

We are now reduced to constructing a map of sheaves of $\pi^{-1}{\ms O}_{X_{\rm{Zar}}}$-algebras
${\ms O}_{X'_{\et}} \rightarrow \mathscr{A}$ over $X_{\et} = X'_{\et}$.  We will do this functorially
on sections over \'etale maps $h:U \rightarrow X$
with $U$ an affine scheme.  Via the equivalence of $X_{\et}$ and $X'_{\et}$, there is a unique
cartesian diagram of schemes 
\begin{equation}\label{uuxxx}
\xymatrix{
{U} \ar[d]_-{h} \ar[r] & {U'} \ar[d]^-{h'} \\
{X} \ar[r] & {X'}}
\end{equation}
in which $h'$ is \'etale, and
$U'$ is affine since the square-zero closed subscheme $U$ is affine.  The natural map
$\mathscr{A}(U) \rightarrow {\ms O}_{X_{\et}}(U) = {\ms O}_U(U)$ is surjective
since $U$ is affine and $\mathscr{J} := \ker(\mathscr{A} \twoheadrightarrow {\ms O}_X)$ is quasi-coherent
on $X_{\et}$, so the natural map $i:U \rightarrow \Spec \mathscr{A}(U)$
is a closed immersion defined by the square-zero ideal $\mathscr{J}(U)$. 

By functoriality the diagram
$$\xymatrix{
{U} \ar[d]_-{h} \ar[r]^-{i} & {\Spec \mathscr{A}(U)} \ar[d] \\
{X} \ar[r] & {\Spec \mathscr{A}(X)}}$$
commutes, so 
the natural map $\Spec \mathscr{A}(U) \rightarrow \Spec \mathscr{A}(X) =: X'$
has restriction to the square-zero closed subscheme $U$ (via $i$) that 
factors through the \'etale map $U' \rightarrow X'$ due to (\ref{uuxxx}).
By the functorial property of \'etale maps of schemes, there is a unique
$X'$-map $\Spec \mathscr{A}(U) \rightarrow U'$ lifting the identity on
the common square-zero closed subscheme $U$ over $X$.  This provides
an $A'$-algebra map ${\ms O}_{X'_{\et}}(U') = {\ms O}_{U'}(U') \rightarrow \mathscr{A}(U)$,
and by the uniqueness of this construction it is easily checked to be functorial in
$U \rightarrow X$.  Hence, we have constructed a $\pi^{-1}({\ms O}_{X'_{\rm{Zar}}})$-algebra map 
${\ms O}_{X'} \rightarrow \mathscr{A}$ over $X_{\et} = X'_{\et}$, so we are done.

{\bf Step 4}. (Essential surjectivity: general case)
By Step 3, essential surjectivity is solved for objects $(X, \mathscr{A}
\twoheadrightarrow {\ms O}_X)$ in $\mathscr{C}$ such that 
the algebraic space $(X,{\ms O}_X)$ is a scheme (with its \'etale topology),
as well as full faithfulness in general. 
In particular, if an object $(X', \mathscr{I})$ in $\mathscr{C'}$ is carried
to $(X, \mathscr{A} \twoheadrightarrow {\ms O}_X)$ for a scheme $X$ then
$X'$ is a scheme by construction.
We will now prove essential surjectivity in general.

Choose an \'etale scheme covering $U \twoheadrightarrow X$, and consider
the object $(U, \mathscr{A}|_{U_{\et}} \twoheadrightarrow {\ms O}_U)$ in
$\mathscr{C}$; for ease of notation we will denote this as $(U, \mathscr{A}|_{U_{\et}})$.
Since $U$ is a scheme, by essential surjectivity that we are assuming
in the scheme case (with the \'etale topology) there is an object $(U', \mathscr{J})$ in $\mathscr{C}'$
that is carried to $(U, \mathscr{A}|_{U_{\et}})$,
and $U'$ is necessarily a scheme.  Likewise, the fiber product
$R := U \times_X U$ is a scheme that is \'etale over $U$ in two ways
and is \'etale over $X$ in a canonical way, so $(R, \mathscr{A}|_{R_{\et}})$ 
is an object in $\mathscr{C}$ with two
natural maps to $(U, \mathscr{A}|_{U_{\et}})$.
Thus, $(R, \mathscr{A}|_{R_{\et}})$ arises from a unique (up to unique isomorphism)
object $(R', \mathscr{K})$
in $\mathscr{C}'$ with $R'$ a scheme such that $R$ is identified with the zero scheme of
$\mathscr{K}$ in $R'$ and the equivalence
$R'_{\et} \simeq R_{\et}$ carries ${\ms O}_{R'}$ to
$\mathscr{A}|_{R_{\et}}$.  Moreover, for each $i \in \{1,2\}$
there is a unique map of schemes 
$p'_i:R' \rightarrow U'$ lifting the \'etale map $p_i:R \rightarrow U$
and identifying ${\ms O}_{U'} \rightarrow p'_{i \ast}({\ms O}_{R'})$ on $U_{\et}$ with $\calA|_{U_{\et}} 
\rightarrow p_{i\ast}(\calA|_{R_{\et}})$.

To prove that each $p'_i$ is \'etale  (and in particular locally of finite
presentation) we give an alternative construction of 
$p'_i$ as follows.  The \'etale map $p_i:R \rightarrow U$ uniquely lifts
to an \'etale map of schemes  $q'_i:R'_i \rightarrow U'$, so
${\ms O}_{R'_i}$ viewed on $R_{\et}$ is identified with ${\ms O}_{U'}|_{R_{\et}} =
\mathscr{A}|_{R_{\et}}$.   Hence, by the unique characterization of
$(R', \mathscr{K})$ we get unique $U'$-isomorphisms $(R', p'_i) \simeq (R'_i,q'_i)$.
This shows that each $p'_i$ is \'etale.

Since the diagram of schemes 
$$\xymatrix{
{R} \ar[r]^-{(p_1,p_2)} \ar[d]  & {U \times_{\Spec \Z} U} \ar[d] \\
{R'} \ar[r]_-{(p'_1,p'_2)} & {U' \times_{\Spec \Z} U'}}$$
is cartesian (due to the constructions of $R'$ and the $p'_i$) and 
the top side is an equivalence relation, the easy full faithfulness for schemes with the Zariski topology
can be used to construct a canonical groupoid structure on $R' \rightarrow U' \times U'$
lifting the one on $R \rightarrow U \times U$.  In other words, $R' \rightrightarrows U'$
is an \'etale equivalence relation in schemes, so 
the quotient $X' = U'/R'$ makes sense as an algebraic space
and Lemma \ref{lemma3} ensures that the commutative diagram
$$\xymatrix{{U} \ar[r] \ar[d] & {U'} \ar[d] \\ {U/R} \ar[r] & {U'/R'}}$$
is cartesian.  But the vertical maps are \'etale coverings and the top side is
a square-zero thickening, so $X' = U'/R'$ is a square-zero thickening of $X = U/R$.  It
is easy to check that this square-zero thickening $X'$ of $X$ solves the original
essential surjectivity problem for the object $(X, \mathscr{A} \twoheadrightarrow {\ms O}_X)$.
\end{proof}

\subsection{Flattening for finite type maps}\label{flatsec}

Our approximation arguments made essential use of results of Raynaud and Gruson from \cite[I,~\S5.7]{rg}.
In this section, we wish to
address a minor issue in the formulation of \cite[I,~5.7.10--5.7.14]{rg}:  these
results are stated for finite type maps $f:X \rightarrow S$
between qcqs algebraic spaces, but the proofs rest on a result
\cite[I,~5.7.9]{rg} which is proved under finite presentation hypotheses.

This does not affect anything in the present paper, as we only use 
\cite[I,~5.7.10--5.7.14]{rg} 
in the proof of Lemma \ref{lemma11}, where such results are applied to noetherian algebraic spaces, 
in which case the finite presentation requirement is satisfied.
Hence, there is circularity in now 
using our approximation results to show that the results in \cite[I,~5.7.10--5.7.14]{rg}
are true as stated with ``finite type'' rather than ``finite presentation''.

Once the two results \cite[I,~5.7.10,~5.7.11]{rg} are proved
as stated in the finite type case, the arguments of Raynaud and Gruson for their subsequent results 
work as written.  We explain below how to make the proofs of \cite[I,~5.7.10,~5.7.11]{rg} work 
in the finite type case by applying the preliminary result
\cite[I,~5.7.9]{rg} (which has finite presentation hypotheses) to
an auxiliary finitely presented map that is provided by Theorem \ref{fpresred}.

\begin{lemma}\label{dimrg} Let $h:V \rightarrow S$ be a flat and locally
finitely presented map of algebraic spaces.  If there exists a schematically
dense open subspace  $U \subseteq V$  with $\dim(U/S) := 
\sup_{s \in |S|} \dim U_s \le n$ then $\dim(V/S) \le n$.  
\end{lemma}

\begin{proof}
By working \'etale-locally on $S$ and then on $V$, we may assume that $S$ and $V$
are affine schemes (so $h$ is finitely presented).  
It therefore suffices to show that if $S = \Spec(A)$ with $A$ a local ring then
$\dim(V_0) \le n$, where $V_0$ is the special fiber.  Letting $\eta$ denote a generic point
of $S$, by flatness and finite presentation for $h$ we have $\dim(V_0) = \dim(V_{\eta})$.  
We can therefore replace $A$ with $A_{\eta}$ so that $A$ is a 0-dimensional local ring. 
In this case, since $U$ is schematically dense in $V$ it is topologically dense,
so $U_{\rm{red}}$ is a topologically dense open subscheme of
the scheme $V_{\rm{red}}$ of finite type over the field $A_{\rm{red}}$.  Hence
$\dim(V) = \dim(V_{\rm{red}}) = \dim(U_{\rm{red}}) \le n$. 
\end{proof}

\begin{proposition}[{\cite[I,~5.7.10]{rg}}]\label{5.7.10}
Let $f:X \rightarrow S$ be a finite type map of qcqs algebraic spaces, and let
$U \subseteq S$ be a quasi-compact open subspace such that $f^{-1}(U) \rightarrow U$
has all fibers with dimension $\le n$, with $n$ a fixed integer.  Then there
exists a $U$-admissible blow-up $g:S' \rightarrow S$ such that the strict transform
$X' \rightarrow S'$ of $f$ with respect to $g$ has all fibers with dimension $\le n$. 
\end{proposition}

\begin{proof}
The case of finitely presented $f$ is what is proved by Raynaud and Gruson,
and we will use Theorem \ref{fpresred} to make their method of proof work in the general 
finite type case.

By Theorem \ref{fpresred}, we may choose a closed immersion $i:X \hookrightarrow
\overline{X}$ over $S$ into an algebraic space $\overline{X}$ of finite presentation over $S$.
Let $\overline{f}:\overline{X} \rightarrow S$ denote the structure map.  By Chevalley's
semi-continuity theorem for fiber dimension of locally finite type maps
\cite[IV$_3$,~13.1.3]{ega} (which works for algebraic spaces by using \'etale localization
to reduce to the case of schemes), the locus
$$W = \{\overline{x} \in \overline{X}\,|\,\dim_{\overline{x}}(\overline{f}^{-1}(\overline{f}(\overline{x})))
\le n\}$$
is open in $\overline{X}$.  

By \cite[I,~5.7.8]{rg}, the quasi-coherent ideal
sheaf $\mathscr{I}$ of $X$ in $\overline{X}$ is the direct limit of its
finite type quasi-coherent ${\ms O}_{\overline{X}}$-submodules
$\mathscr{I}_{\lambda}$.  Hence, $X = \invlim X_{\lambda}$
where $\{X_{\lambda}\}$ is the inverse system (with affine transition maps)
of finitely presented closed subspaces of $\overline{X}$ that contain $X$.  
If $f_{\lambda}:X_{\lambda} \rightarrow S$ denotes the finitely presented structure map
then $\{f_{\lambda}^{-1}(U)\}$ is an inverse system of finitely presented
closed subspaces of $\overline{f}^{-1}(U)$ with $\invlim f_{\lambda}^{-1}(U) = f^{-1}(U)$.
By hypothesis, the open subspace $W$ in the qcqs $\overline{X}$ contains
$f^{-1}(U)$. Hence, if $Z$ denotes $\overline{X} - W$ with its reduced structure
then $\{Z \cap f_{\lambda}^{-1}(U)\}$ is an inverse system of
qcqs algebraic spaces with limit $Z \cap f^{-1}(U) = \emptyset$.
It follows that $Z \cap f_{\lambda}^{-1}(U)$ is empty for sufficiently 
large $\lambda$, so $f_{\lambda}^{-1}(U) \subseteq W$ for sufficiently large $\lambda$. 
For any $\overline{x} \in f_{\lambda}^{-1}(U)$, the fiber 
$f_{\lambda}^{-1}(f_{\lambda}(\overline{x})) = X_{\lambda} \cap \overline{f}^{-1}(\overline{f}(
\overline{x}))$ has dimension $\le n$ at $\overline{x}$
since $\overline{x} \in W$.   In other words, by replacing
$\overline{X}$ with such an $X_{\lambda}$ we get to the situation in which
all fibers of $\overline{f}^{-1}(U) \rightarrow U$ have dimension $\le n$.

Let $\mathscr{M} = {\ms O}_{\overline{X}}$.  Since $\overline{f}$ is finitely presented
and the open set $V = \emptyset$ in $\overline{f}^{-1}(U)$ has complement
whose fibers over $U$ have dimension $< n+1$, we can say that
$\mathscr{M}|_{\overline{f}^{-1}(U)}$ is $U$-flat in dimensions $\ge n+1$
(in the sense of \cite[I,~5.2.1]{rg}).  Thus, by \cite[I,~5.7.9]{rg} there is
a $U$-admissible blow-up $g:S' \rightarrow S$ such that
the strict transform $\overline{f}':\overline{X}' \rightarrow S'$ of $\overline{f}$ relative to $g$ is
finitely presented over $S'$ and there is an $S'$-{\em{flat}} quasi-compact open
$\overline{V}' \subseteq \overline{X}'$ such that $\dim((\overline{X}' - \overline{V}')/S') \le n$
(that is, all fibers have dimension at most $n$).   By \cite[I,~5.1.2(v),~5.1.4]{rg} (whose proofs
also work for algebraic spaces instead of schemes), we can make a further 
$U$-admissible blow-up on $S'$ to get to the situation in which $\overline{U}' := {\overline{f}'}^{-1}(U')$
is a schematically dense open subspace of $\overline{X}'$, where $U'$ denotes 
the preimage of $U$ in $S'$. Note that $\dim(\overline{U}'/S') \le n$
since $S' \rightarrow S$ is a $U$-admissible blow-up and $\dim(\overline{f}^{-1}(U)/U) \le n$. 

We claim that $\dim(\overline{X}'/S') \le n$, so since
the strict transform $X' \rightarrow S'$ of $f$ relative to $g$ is a closed subspace
of $\overline{X}'$ over $S'$ it would follow that $\dim(X'/S') \le n$, as desired.
Since $\dim((\overline{X}' - \overline{V}')/S') \le n$, we just have to check that
$\dim(\overline{V}'/S') \le n$.  The overlap $\overline{V}' \cap \overline{U}'$
is a schematically dense open subspace of $\overline{V}'$ with fiber dimensions at
most $n$ (since $\dim(\overline{U}'/S') \le n$), so we may apply Lemma \ref{dimrg} to conclude. 
\end{proof}

The proof of  \cite[I,~5.7.11]{rg} begins with a finite type map $f:X \rightarrow S$
and asserts at the start (using \cite[I,~5.7.9]{rg}) that there exists a 
blow-up $S' \rightarrow S$ (along a finite-type quasi-coherent ideal)
with respect to which the strict transform of $X$ 
is flat and finitely presented over $S'$.  Since \cite[I,~5.7.9]{rg} has finite presentation
hypotheses, an additional argument is needed to ensure the existence of such a blow-up when $f$ is 
of finite type rather than of finite presentation.  We now prove this via a further application of Theorem 
\ref{fpresred}. 

\begin{proposition}\label{5.7.11} 
Let $f:X \rightarrow S$ be a map of finite type between qcqs algebraic spaces,
and let $U \subseteq S$ be a quasi-compact open subspace such that 
$f^{-1}(U) \rightarrow U$ is flat and finitely presented
$($e.g., an open immersion$)$.   There exists a $U$-admissible blow-up
$g:S' \rightarrow S$ such that the strict transform $f':X' \rightarrow S'$ of $f$ relative to $g$
is flat and finitely presented. 
\end{proposition}

\begin{proof}
By Theorem \ref{fpresred},
we may choose a closed
immersion $i:X \hookrightarrow \overline{X}$ over $S$ into an algebraic
space of finite presentation over $S$.  Let $\overline{f}:\overline{X} \rightarrow S$
denote the structure map, and $\mathscr{M} = i_{\ast}({\ms O}_X)$.  Since
$f^{-1}(U) \rightarrow U$ is finitely presented and flat, $\mathscr{M}|_{\overline{f}^{-1}(U)}$
is finitely presented over ${\ms O}_{\overline{f}^{-1}(U)}$ and $U$-flat.   In other words,
$\mathscr{M}|_{\overline{f}^{-1}(U)}$ is $U$-flat in dimension $\ge 0$
(in the sense of \cite[I,~5.2.1]{rg}).   Thus, by the finite presentedness
of $\overline{f}$ we can apply \cite[I,~5.7.9]{rg} to construct 
a $U$-admissible blow-up $g:S' \rightarrow S$ such that 
the strict transform $\overline{X}' \rightarrow S'$ of $\overline{f}$
relative to $g$ is finitely presented and
the strict transform $\mathscr{M}'$ of $\mathscr{M}$ is
finitely presented over ${\ms O}_{\overline{X}'}$ and $S'$-flat.   But
$\mathscr{M}'$ is the quotient of ${\ms O}_{\overline{X}'}$
corresponding to the closed immersion $X' \hookrightarrow \overline{X}'$
of strict transforms relative to $g$.  Hence, $X'$ is $S'$-flat and its closed immersion
into $\overline{X}'$ is finitely presented, so $X' \rightarrow S'$ is also
finitely presented (as $\overline{X}' \rightarrow S'$ is finitely presented).
\end{proof}

\end{appendix}

\end{document}